\documentclass[article]{amsart}
\usepackage{cases}
\usepackage{amsmath, amsfonts, pifont, amssymb}
\usepackage{verbatim}
\usepackage{geometry}
\newtheorem{theorem}{Theorem}[section]
\newtheorem{lemma}[theorem]{Lemma}

\newtheorem{corollary}[theorem]{Corollary}

\newcommand{\beq}{\begin{equation}}
\newcommand{\eeq}{\end{equation}}
\newcommand{\beqq}{\begin{equation*}}
\newcommand{\eeqq}{\end{equation*}}

\theoremstyle{definition}
\newtheorem{definition}[theorem]{Definition}

\theoremstyle{remark}

\newtheorem{conjecture}[theorem]{Conjecture}
\numberwithin{equation}{section}

\numberwithin{equation}{section}

\begin{document}

\title[Scattering for cubic NLS on $\mathbb{R}^2$ $\times$ $\mathbb{T}^2$]{Global well-posedness and scattering for the defocusing cubic Schr{\"o}dinger equation on waveguide $\mathbb{R}^2$ $\times$ $\mathbb{T}^2$}

\author{Zehua Zhao}

\maketitle

\begin{abstract}
In this paper, we consider the problem of large data scattering for the defocusing cubic nonlinear Schr{\"o}dinger equation on $\mathbb{R}^2$ $\times$ $\mathbb{T}^2$. This equation is critical both at the level of energy and mass. The key ingredients are global-in-time Stricharz estimate, resonant system approximation, profile decomposition and energy induction method. Assuming the large data scattering for the 2d cubic resonant system, we prove the large data scattering for this problem. This problem is the cubic analogue of [13].
\end{abstract}
\bigskip

\bigskip
\section{Introduction}
\noindent We consider the defocusing cubic nonlinear Schr{\"o}dinger equation on $\mathbb{R}^2\times \mathbb{T}^2$,

\begin{equation}\label{equation}
\aligned
(i\partial_t+ \Delta_{\mathbb{R}^{2} \times \mathbb{T}^{2}}) u &= F(u) = |u|^{2} u, \\
u(0,x) &= u_{0} \in H^{1}(\mathbb{R}^{2} \times \mathbb{T}^{2})
\endaligned
\end{equation}
where $\Delta_{\mathbb{R}^{2} \times \mathbb{T}^{2}}$ is the Laplace-Beltrami operator on $\mathbb{R}^2\times \mathbb{T}^2$ and $u:\mathbb{R}\times \mathbb{R}^2\times \mathbb{T}^2\rightarrow \mathbb{C}$ is a complex-valued function. We are interested in the long-time behavior of solutions of initial value problem (1.1) with large data.\vspace{3mm}

\noindent First, the equation (1.1) has the following conserved quantities: 
\begin{enumerate}
\item Energy $E$: $E(u)=\frac{1}{2}||\nabla u||_{L^2(\mathbb{R}^2\times \mathbb{T}^2)}^2+\frac{1}{4}||u||_{L^4(\mathbb{R}^2\times \mathbb{T}^2)}^4$.
\item Mass $M$: $M(u)=||u||_{L^2(\mathbb{R}^2\times \mathbb{T}^2)}^2$.
\item Momentum: $P(u)=Im \int_{}  \bar{u}(x,y,t)\nabla u(x,y,t)dxdy$.
\end{enumerate}
We define ``full energy'' $L$ as follows: $L(u)=\frac{1}{2}M(u)+E(u)$, which is obviously conserved.\vspace{2mm}

\noindent Generally, the equation (1.1) is a special case of the general defocusing nonlinear Schr{\"o}dinger equation on the waveguide $\mathbb{R}^m \times \mathbb{T}^n$:
\begin{equation}\label{equation}
\aligned
(i\partial_t+ \Delta_{\mathbb{R}^{m} \times \mathbb{T}^{n}}) u &= F(u) = |u|^{p-1} u, \\
u(0,x) &= u_{0} \in H^{1}(\mathbb{R}^{m} \times \mathbb{T}^{n}).
\endaligned
\end{equation}

\noindent As for the background, it is well known that there are many existing results regarding NLS problems on Euclidean space. In this paper, the nonlinear Schr{\"o}dinger equation is discussed on a semiperiodic space, i.e. $\mathbb{R}^2\times \mathbb{T}^2$. The motivation is to better understand the broad question of the effect of the geometry of the domain on the asymptotic behavior of large solutions to nonlinear dispersive equations. The study of solutions of the nonlinear Schr{\"o}dinger equation on compact or partially compact domains has been the subject of many works. Such equations have also been studied in applied sciences and various background. Those product spaces ($\mathbb{R}^{m} \times \mathbb{T}^{n}$) are also called waveguide manifolds. In paper [25, 26], there is a good overview of the main results for waveguide manifolds.\vspace{3mm}

\noindent The studies on global well-posedness for energy critical and subcritical equations seem to point to the absence of any geometric obstruction to global existence. Moreover, it is clear that the geometry influences the asymptotic dynamics of solutions. Thus, it is meaningful to explore when one can obtain the simplest asymptotic behavior, i.e. scattering, which means that all nonlinear solutions asymptotically resemble linear solutions. Based on the theory of NLS on Euclidean space, i.e. $\mathbb{R}^d$, the equation 
\begin{equation}\label{equation}
\aligned
(i\partial_t+ \Delta_{\mathbb{R}^{d}}) u &= F(u) = |u|^{p-1} u, \\
u(0,x) &= u_{0} \in H^{1}(\mathbb{R}^{d})
\endaligned
\end{equation}
\noindent would scatter in the range $1+\frac{4}{d} \leq p\leq 1+\frac{4}{d-2}$. When $p=1+\frac{4}{d}$, the equation (1.3) is mass critical; when $p=1+\frac{4}{d-2}$, the equation (1.3) is energy critical; when  $p<1+\frac{4}{d}$, the equation is mass subcritical; when  $p>1+\frac{4}{d-2}$, the equation (1.3) is energy supercritical; when $1+\frac{4}{d} < p< 1+\frac{4}{d-2}$, the equation is mass supercritical and energy subcritical. \vspace{3mm}

\noindent Naturally, we are also interested in the range of the nonlinearity index $p$ (in (1.2)) for well-posedness and scattering of the NLS on  $\mathbb{R}^m\times \mathbb{T}^n$. Based on the existing results and theories, we expect that the solution of (1.2) globally exists and scatters in the range $1+\frac{4}{m}\leq p \leq 1+\frac{4}{m+n-2}$. And fortunately the index ($p=3$) in equation (1.1) lies in the range (exactly at the endpoints of the interval), so it is reasonable for us to consider this problem.\vspace{3mm}

\noindent Moreover, we are mainly inspired by a related result [13] which studies the defocusing quintic NLS on space $\mathbb{R}$ $\times$ $\mathbb{T}^2$ by Zaher Hani and Benoit Pausader. For that problem, the defocusing NLS equation is also critical both at the level of energy and mass. Moreover, [15, 16, 17, 18, 19] contain some of the important related results. Also, the introduction of [7] has a summary of the main known results about NLS problems on waveguides, i.e. $\mathbb{R}^m\times \mathbb{T}^n$.\vspace{3mm}

\noindent As in [13], we need to assume the large data scattering for a cubic resonant system (1.6) as well. The next task for us is to deal with the regarding resonant system, i.e. to prove this conjecture 1.2. An inspiring thing is that K. Yang and L. Zhao ([35]) have proved the large data scattering for a similar resonant system recently. Conjecture 1.2 seems to be a reachable problem to work on by using scattering result for 2d mass critical NLS in B. Dodson [9], and we leave it for a later work.\vspace{3mm} 

\noindent \textbf{One word on cubic NLS problems on 4d waveguides:} While scattering holds for the cubic equation on $\mathbb{R}^4$ (see [24]), $\mathbb{R}^3\times \mathbb{T}$ (see an oncoming paper which is our second project of a series work) and $\mathbb{R}^2\times \mathbb{T}^2$ (see this paper), it is not expected to hold on $\mathbb{T}^4$. The situation on $\mathbb{R}^2\times \mathbb{T}^2$ seems to be a borderline case for this question, i.e. defocusing cubic NLS equation on four dimensional waveguides. \vspace{3mm}

\noindent The first result asserts that the small data leads to solutions that are global and scattering.
\begin{theorem}\label{theorem} There exists $\delta > 0$ such that any initial data $u_0 \in H^1(\mathbb{R}^2 \times \mathbb{T}^2)$ satisfying \vspace{5mm}
\[ ||u_0||_{H^1(\mathbb{R}^2 \times \mathbb{T}^2)} \leq \delta 
\]
\noindent leads to a unique global solution $u \in X^1_{c}(\mathbb{R})$ which scatters in the sense that there exists $v^{\pm \infty} \in H^1(\mathbb{R}^2 \times \mathbb{T}^2)$ such that
\begin{equation}\label{equation}
||u(t)-e^{it\Delta_{\mathbb{R}^2 \times \mathbb{T}^2}}v^{\pm \infty}||_{H^1(\mathbb{R}^2 \times \mathbb{T}^2)} \rightarrow 0 \quad as \quad t \rightarrow \pm \infty .
\end{equation}
\end{theorem}
\noindent The uniqueness space $X^1_c\subset C_t(\mathbb{R}:H^1(\mathbb{R}^2\times \mathbb{T}^2))$ was essentially introduced by Herr-Tataru-Tzvetkov [15]. In order to extend our analysis to large data, we use a method formalized in [20, 21]. One key ingredient is a linear and nonlinear profile decomposition for solutions with bounded energy. The so-called profiles correspond to sequences of solutions exhibiting an extreme behavior. It is there that the ``energy critical" and ``mass critical" nature of our equation become manifest.\vspace{3mm}

\noindent For this problem, in view of the scaling-invariant of the initial value problem (1.1) under 
\[\mathbb{R}_x^2 \times \mathbb{T}_y^2 \rightarrow M_{\lambda}:=\mathbb{R}_x^2 \times (\lambda^{-1}\mathbb{T}^2)_{y},\quad u \rightarrow \tilde{u}(x,y,t)=\lambda u(\lambda x,\lambda y,\lambda^2 t).\]
\noindent There are two situations:\vspace{3mm}

\noindent When $\lambda \rightarrow 0$, the manifolds $M_\lambda$ will be similar to $\mathbb{R}^4$ and we can use the four dimensional energy critical result [24] by (E. Ryckman and M. Visan). The appearance is a manifestation of the energy-critical nature of the nonlinearity. This corresponds in $M_1$ to solutions with initial data (This behavior corresponds to Euclidean profiles which we will define more precisely in section 5 and section 6.)
\[ u^{\lambda}(x,y,0)=\lambda^{-1}\phi(\lambda^{-1} (x,y)) \quad \phi \in C^{\infty}_0(\mathbb{R}^4),\lambda \rightarrow 0.
\]

\noindent When $\lambda \rightarrow \infty$, the manifolds $M_\lambda$ become thinner and thinner and resembles $\mathbb{R}^2$. The problem will become similar to the cubic mass critical NLS problem on $\mathbb{R}^2$ :
\begin{equation}\label{equation}
(i\partial_t+\Delta_x )u=|u|^2u ,\quad u(0) \in H^1(\mathbb{R}^2) .
\end{equation} 
\noindent Those solutions on $M_{\lambda}$ correspond to solutions on $M_1$ with initial data
\begin{equation}
u^{\lambda}(x,y,0)=\lambda^{-1}\phi(\lambda^{-1}x,y) \quad \phi \in C^{\infty}_0(\mathbb{R}^2\times \mathbb{T}^2),\lambda \rightarrow \infty.
\end{equation}
\noindent The second situation is closely related to the following cubic resonant system.\vspace{3mm}

\noindent \textbf{The cubic resonant system:} We consider the cubic resonant system, 
\begin{equation} \label{equation}
\aligned
&(i\partial_t+\Delta_{x})u_j=\sum_{(j_1,j_2,j_3)\in R(j)} u_{j_{1}}\bar{u}_{j_{2}}u_{j_{3}} , \\
&R(j)=\{(j_1,j_2,j_3) \in (\mathbb{Z}^2)^3:j_1-j_2+j_3=j \quad and \quad  |j_1|^2-|j_2|^2+|j_3|^2=|j|^2 \}
\endaligned
\end{equation}
\noindent with unknown $\vec{u}=\{u_j \}_{j \in \mathbb{Z}^2}$, where $u_j :\mathbb{R}_x^2 \times \mathbb{R}_t \rightarrow \mathbb{C}$. \vspace{3mm}

\noindent In the special case when $u_j=0$ for $j \neq 0$, it is exactly equation (1.5). Similar systems of nonlinear Schr{\"o}dinger equations arise in the study of nonlinear optics in waveguides.\vspace{3mm}

\noindent As we show in Section 8, the system (1.6) is Hamiltonian and it has a nice local theory and retains many properties of (1.5). In view of this and of the 2d mass critical NLS result [9] (by B. Dodson), it seems reasonable to assume the following conjecture. Another reason for us to assume this is that K. Yang and L. Zhao [35] have proved the large data scattering for a similar resonant system when index $j\in \mathbb{Z}$ (1d case).

\begin{conjecture} Let $E \in (0,\infty)$. For any smooth initial data $\vec{u}_0$ satisfying:
\[ E_{ls}(\vec{u}_0)=\frac{1}{2}\sum\limits_{j \in \mathbb{Z}^2} \langle j \rangle^2 ||u_{0,j}||^2_{L^2_{x}(\mathbb{R}^2)}\leq E.
\]
\noindent There exists a global solution $\vec{u}(t),\vec{u}(t=0)=\vec{u}_0 $ with conserved $E_{ls}(\vec{u}(t))$ satisfying: 
\begin{equation}\label{equation}
||\vec{u}||^{2}_{\vec{W}}:=\sum\limits_{j \in \mathbb{Z}^2}\langle j \rangle^2 ||u_{j}||^2_{L^4_{t,x}(\mathbb{R}_t \times \mathbb{R}_x^2)} \leq \Lambda_{ls}(E_{ls}(\vec{u}_{0}))
\end{equation}
\noindent for some finite non-decreasing function $\Lambda_{ls}(E)$.
\end{conjecture}
\noindent  \emph{Remark.} As for a more general case, when $n=2$ and $p=1+\frac{4}{m}$, the initial value problem (1.2) is also both critical at the level of mass and energy. If $m>2$, the index $p$ would no longer be an integer, which may cause some trouble for us to use the resonant system approximation.\vspace{3mm}

\noindent We now give the main result of this paper which asserts the large data scattering for (1.1) conditioned on Conjecture 1.1.
\begin{theorem}\label{theorem}
Assume that $\textit{Conjecture 1.1}$ holds for all $E \leq E_{max}^{ls}$, then any initial data $u_0 \in H^{1}(\mathbb{R}^2\times \mathbb{T}^2)$ satisfying 
\[ L(u_0)=\int_{\mathbb{R}^2 \times \mathbb{T}^2} (\frac{1}{2} |\nabla u_0|^2+\frac{1}{2}  |u_0|^2 +\frac{1}{4}|u_0|^{4} )dx \leq E_{max}^{ls}
\]leads to a solution $u \in X^1_{c}(\mathbb{R})$ which is global, and scatters (in the sense of (1.3)). In particular if $E_{max}^{ls}=+\infty$, then all solutions of (1.1) with finite energy and mass scatter.
\end{theorem}
\noindent As a consequence of the local theory for the system (1.6), Conjecture 1.1 holds below a nonzero threshold $E_{max}^{ls}>0$, so Theorem 1.2 is non-empty and indeed strengthens Theorem 1.1. Another point worth mentioning is that while Theorem 1.2 is stated as an implication of Conjecture 1.1, it is actually an equivalence as it is easy to see that one can reverse the analysis needed to understand the behavior of large-scale profile initial data for (1.1) in order to control general solutions of (1.6) and prove Conjecture 1.1 assuming that Theorem 1.2 holds (cf. Section 8). At last, the ``scattering threshold" for (1.1) and the resonant system (1.6) are same. \vspace{3mm}

\noindent The proof of the Theorem 1.2 follows from a standard skeleton based on the Kenig-Merle machinery ([20, 21]) for proving global well-posedness and scattering problem. Mainly there are several important points of the proof: global Strichartz estimates, large-scale profile and the resonant cubic system, profile decomposition and energy induction method.\vspace{3mm}

\noindent \textbf{The organization of this paper:} in Section 2, we introduce some basic notations and important function spaces; in Section 3, we prove the global-in-time Strichartz estimate; in Section 4, we prove the local well-posedness and small data scattering of (1.1); in Section 5, we describe the Euclidean profiles and large-scale profiles and construct Euclidean approximation and resonant system approximation in Section 6, we obtain a good linear profile decomposition that leads us to analyze the large data case; in Section 7, we prove the contradiction argument leading to Theorem 1.2 (Main Theorem); in Section 8, we prove the local theory for the cubic resonant system (1.6) and also give a proof for Lemma 3.3 (local-in-time $L^p$ estimate) in Section 3 for completeness.\vspace{3mm}

\section{Notations and function spaces}
\noindent About the notation, we write $A \lesssim B$ to say that there is a constant $C$ such that $A\leq CB$. We use $A \simeq B$ when $A \lesssim B \lesssim A $. Particularly, we write $A \lesssim_u B$ to express that $A\leq C(u)B$ for some constant $C(u)$ depending on $u$.\vspace{3mm}

\noindent In addition to the usual isotropic Sobolev spaces $H^s(\mathbb{R}^2 \times \mathbb{T}^2)$,  we have non-isotropic versions. For $s_1,s_2 \in \mathbb{R}$ we define:
\begin{equation}\label{equation} 
H^{s_1,s_2}(\mathbb{R}^2 \times \mathbb{T}^2)=\{ u:\mathbb{R}^2 \times \mathbb{T}^2 \rightarrow \mathbb{C}: \langle \xi \rangle^{s_1}   \langle n \rangle^{s_2}\hat{u}(\xi,n)  \in L^2_{\xi,n}(\mathbb{R}^2 \times \mathbb{Z}^2)\}. 
\end{equation}
Particularly $H^{0,1}(\mathbb{R}^2 \times \mathbb{T}^2)$ is a Hilbert space with inner product: 
\[ \langle u,v \rangle_{H^{0,1}}=\langle u,v \rangle_{L^2}+\langle \nabla_y u, \nabla_y v \rangle_{L^2}.
\]
\noindent We can also define a discrete analogue. For $\vec{\phi}=\{\phi_p \}_{p \in \mathbb{Z}^2}$ a sequence of real-variable functions, we let
\begin{equation}\label{equation}
h^{s_1}H^{s_2}:=\{\vec{\phi}=\{ \phi_{p} \} : ||\vec{\phi}||_{h^{s_1}H^{s_2}}^2=\sum\limits_{p \in \mathbb{Z}^2} \langle p \rangle^{2s_1}||\phi_p||_{H^{s_2}}^{2} < +\infty \}. 
\end{equation}
\noindent We can naturally identify $ H^{0,1}(\mathbb{R}^2 \times \mathbb{T}^2)$ and $h^1L^2$ by via the Fourier transform in the periodic variable $y$.\vspace{3mm}

\noindent \textbf{Function spaces.} In this paper, we will use some function spaces. For example, the $X^1$ space was essentially introduced by Herr-Tataru-Tzvetkov [15] (see also [13]). \vspace{3mm}

\noindent For $C=[-\frac{1}{2},\frac{1}{2})^4 \in \mathbb{R}^4$ and $z\in \mathbb{R}^4$, we denote by $C_z=z+C$ the translate by $z$ and define the sharp projection operator $P_{C_z}$ as follows: ($\mathcal{F}$ is the Fourier transform): 
\[
\mathcal{F}(P_{C_z} f)=\chi_{C_z}(\xi) \mathcal{F} (f)  (\xi).
\]
\noindent We use the same modifications of the atomic and variation space norms that were employed in some other papers [15, 16]. Namely, for $s\in \mathbb{R}$, we define:
\[ \|u\|_{X_0^s(\mathbb{R})}^2=\sum_{z\in \mathbb{Z}^4} \langle z \rangle^{2s} \|P_{C_z} u\|_{U_{\Delta}^2(\mathbb{R};L^2)}^2
\]
\noindent and similarly we have,  
\[\|u\|_{Y^s(\mathbb{R})}^2=\sum_{z\in \mathbb{Z}^4} \langle z \rangle^{2s} \|P_{C_z} u\|_{V_{\Delta}^2(\mathbb{R};L^2)}^2\]          
\noindent where the $U_{\Delta}^p$ and $V_{\Delta}^p$ are the atomic and variation spaces respectively of functions on $\mathbb{R}$ taking values in $L^2(\mathbb{R}^2\times \mathbb{T}^2)$. There are some nice properties of those spaces. We refer to [15, 16] for the description and properties. For convenience, we also give the some definitions here. 

\begin{definition}Let $1\leq p < \infty$, and $H$ be a complex Hilbert space. A $U^p$-atom is a piecewise defined function, $a:\mathbb{R} \rightarrow H$，
\[ a=\sum_{k=1}^{K}\chi_{[t_{k-1},t_k)}\phi_{k-1}
\]
\noindent where $\{t_k\}_{k=0}^{K} \in \mathcal{Z}$ and $\{\phi_k\}_{k=0}^{K-1} \subset H$ with $\sum_{k=0}^{K}||\phi_k||^p_H=1$. Here we let $\mathcal{Z}$ be the set of finite partitions $-\infty<t_0<t_1<...<t_K\leq \infty$ of the real line.\vspace{3mm}

\noindent The atomic space $U^p(\mathbb{R};H)$ consists of all functions $u:\mathbb{R}\rightarrow H$ such that

\noindent                           $u=\sum_{j=1}^{\infty}\lambda_j a_j$ for $U^p$-atoms $a_j$, $\{\lambda_j\} \in l^1$, 

\noindent with norm 
 \[||u||_{U^p}:=inf\{\sum^{\infty}_{j=1}|\lambda_j|:u=\sum_{j=1}^{\infty}\lambda_j a_j,\lambda_j\in \mathbb{C}, a_j \quad U^p\textmd{-atom}\}.\]
\end{definition}
\begin{definition}Let $1\leq p < \infty$, and $H$ be a complex Hilbert space. We define $V^p(\mathbb{R},H)$ as the space of all functions $v:\mathbb{R} \rightarrow H$ such that
\[ ||u||_{V^p}:=\sup\limits_{\{t_k\}^K_{k=0} \in \mathcal{Z}}(\sum_{k=1}^{K}||v(t_k)-v(t_{k-1})||^p_{H})^{\frac{1}{p}} \leq +\infty,
\]
\noindent where we use the convention $v(\infty)=0$.
\noindent Also, we denote the closed subspace of all right-continuous functions $v:\mathbb{R}\rightarrow H$ such that $\lim\limits_{t\rightarrow -\infty}v(t)=0$ by $V^p_{rc}(\mathbb{R},H)$. 
\end{definition}
\begin{definition}For $s\in \mathbb{R}$, we let $U^p_{\Delta}H^s$ resp. $V^p_{\Delta}H^s$ be the spaces of all functions such that $e^{-it\Delta}u(t)$ is in $U^p(\mathbb{R},H^s)$ resp. $V^p_{rc}(\mathbb{R},H)$, with norms
\[||u||_{U^p_{\Delta}H^s}=||e^{-it\Delta}u||_{U^p(\mathbb{R},H^s)}, \quad ||u||_{V^p_{\Delta}H^s}=||e^{-it\Delta}u||_{V^p(\mathbb{R},H^s)}.
\]
\end{definition}
\noindent For this problem, we choose $H$ to be $L^2(\mathbb{R}^2\times \mathbb{T}^2)$. Norms $X_0^s$ and $Y^s$ are both stronger than the $L^{\infty}(\mathbb{R};H^s)$ norm and weaker than the norm $U^2_{\Delta}(\mathbb{R}:H^s)$. Moreover, they satisfy the following property (for $p>2$):
\[U^2_{\Delta}(\mathbb{R}:H^s) \hookrightarrow X_0^s\hookrightarrow  Y^s \hookrightarrow V^2_{\Delta}(\mathbb{R}:H^s)\hookrightarrow  U^p_{\Delta}(\mathbb{R}:H^s) \hookrightarrow L^{\infty}(\mathbb{R};H^s).
\]
\noindent For an interval $I \subset \mathbb{R}$, we can also define the restriction norms $X^s_0(I)$ and $Y^s(I)$ in the natural way:
\noindent                                               $||u||_{X_0^s(I)}=$ inf $\{||v||_{X_0^s(\mathbb{R})}:v\in X_0^s(\mathbb{R})$ satisfying $v_{|I}=u_{|I}\}$.

\noindent And similarly for $Y^s(I)$. \vspace{3mm}

\noindent A modification for to $X_0^s(\mathbb{R})$: \vspace{3mm}

\noindent $X^s(\mathbb{R}):=\{u:\phi_{-\infty}=\lim\limits_{t\rightarrow -\infty}e^{-it\Delta}u(t) \textmd{ exists in }  H^s,u(t)-e^{it\Delta}\phi_{-\infty} \in X^s_0(\mathbb{R})\} $ equipped with the norm:
\begin{equation}
||u||^2_{X^s(\mathbb{R})}=||\phi_{-\infty}||^2_{H^s(\mathbb{R}^2\times \mathbb{T}^2)}+||u-e^{it\Delta}\phi_{-\infty}||^2_{X_0^s(\mathbb{R})}.
\end{equation}\noindent Our basic space to control solutions is $X^1_c(I)=X^1(I)\cap C(I:H^1).$  Also we use $X^1_{c,loc}(I)$ to express the set of all solutions in $C_{loc}(I:H^1)$ whose $X^1(J)$-norm is finite for any compact subset $J \subset I$.\vspace{3mm}

\noindent In order to control the nonlinearity on interval $I$, we need to define `$N$ -Norm' as follows, on an interval $I=(a,b)$ we have:
\begin{equation}\label{equation}
\| h\|_{N^s(I)}=\|\int_{a}^{t} e^{i(t-s)\Delta} h(s) ds \|_{X^s(I)} .
\end{equation}
\noindent And then we can define the following spacetime norm, i.e. `Z-norm' by
\[
\|u\|_{Z(I)}=(\sum_{N\geq 1}N^{2} \|1_I(t) P_N u\|_{L_{x,y,t}^4(\mathbb{R}^2\times \mathbb{T}^2 \times \mathbb{R})}^4)^\frac{1}{4}.
\]
\noindent $Z$ is a weaker norm than $X^1$, in fact:
\[||u||_{Z(I)} \lesssim ||u||_{X^1(I)},
\]\noindent which follows from Strichartz estimate (see Section 3).\vspace{3mm}

\noindent We also need the following theorem which has analogues in [13, 15, 16].
\begin{theorem}\label{theorem}[[13, 15, 16]] If $f\in L^1_t(I,H^1(\mathbb{R}^2\times \mathbb{T}^2))$, then 
\[||f||_{N(I)} \lesssim \sup_{v\in Y^{-1}(I),||v||_{Y^{-1}(I)}\leq 1} \int_{I \times (\mathbb{R}^2 \times \mathbb{T}^2)} f(x,t)\overline{v(x,t)}dxdt.
\]
Also, we have the following estimate holds for any smooth function $g$ on an interval $I=[a,b]$: 
\[ ||g||_{X^1(I)}\lesssim ||g(0)||_{H^1(\mathbb{R}^2\times \mathbb{T}^2)}+(\sum_N ||P_N(i\partial_t+\Delta)g||^2_{L^1_t(I,H^1(\mathbb{R}^2\times \mathbb{T}^2))})^{\frac{1}{2}}.
\]
\end{theorem}
\noindent \emph{Proof:} The proof follows as in [15, Proposition 2.11] and [16, Proposition 2.10].
\section{Global Strichartz estimate}
\begin{theorem}\label{theorem}Then we can prove the following Strichartz Estimate:
\begin{equation}\label{equation}
\|e^{it\Delta_{\mathbb{R}^2\times\mathbb{T}^2}}P_{\leq N}u_0\|_{{l_{\gamma}^q L_{x,y,t}^p}(\mathbb{R}^2\times \mathbb{T}^2\times[2\pi\gamma,2\pi(\gamma+1)])}\lesssim N^{2-\frac{6}{p}}\|u_0\|_{L^2( \mathbb{R}^2 \times \mathbb{T}^2)} 
\end{equation}
\noindent whenever 
\begin{equation}\label{equation} 
p>\frac{10}{3} \quad and \quad \frac{1}{q}+\frac{1}{p}=\frac{1}{2}.
\end{equation}
\end{theorem}
\emph{Remark.} We explain the norm $l_{\gamma}^q L_{x,y,t}^p$. First, we decompose $\mathbb{R}$ into $\mathbb{R}=\cup_{\gamma\in \mathbb{Z}}I_{\gamma}=\cup_{\gamma\in \mathbb{Z}}2\pi[\gamma,\gamma+1)$, where $I_{\gamma}=[2\pi\gamma,2\pi(\gamma+1)]$. Moreover, we take the $L^p$-timespace norm on each $I_{\gamma}$, and then take $l^q$-sequence norm for the sequence $\{L_{x,y,t}^p(\mathbb{R}^2\times \mathbb{T}^2\times I_{\gamma})\}_{\gamma}$. \vspace{2mm}

\noindent \emph{Proof:} The main idea of the proof is similar to [13, Theorem 3.1], i.e. using $T-T^*$ argument, a partition of unity and then estimating the diagonal part and non-diagonal part separately. One remarkable difference is that in the diagonal estimate part, we can not use Bourgain' s $L^p$ estimate on $\mathbb{T}^2$ ([2]) directly as in [13] since we need a Stricharz estimate with a threshold less than 4 (precisely it is $\frac{10}{3}$) to do the interpolation later. And we use Hardy-Littlewood circle method as in [19, Proposition 2.1] to obtain the local-in-time $L^p$ estimate.\vspace{3mm}

\noindent First, we prove a more precise conclusion and we can get the estimate by duality:

\begin{lemma}
For any $h\in C_c^{\infty}(\mathbb{R}_x^2\times \mathbb{T}_y^2\times \mathbb{R}_t)$, the inequality
\begin{equation}\label{equation}
\aligned
&\|\int_{s\in \mathbb{R}} e^{-is\Delta_{\mathbb{R}^2\times \mathbb{T}^2}} P_{\leq N} h(x,y,s) ds\|_{L_{x,y}^{2}(\mathbb{R}_x^2\times \mathbb{T}_y^2)} \\ &\lesssim N^{2-\frac{6}{p}} \|h\|_{l_{\gamma}^2 L_{x,y,t}^{p^{'}}(\mathbb{R}^2 \times \mathbb{T}^2\times [2\pi\gamma,2\pi(\gamma+1)])}+N^{1-\frac{8}{3p}} \|h\|_{l_{\gamma}^{q^{'}} L_{x,y,t}^{p^{'}}(\mathbb{R}^2\times \mathbb{T}^2\times [2\pi\gamma,2\pi(\gamma+1)])}  
\endaligned
\end{equation}
\noindent holds for any (p,q) satisfies (3.2).
\end{lemma}
\noindent \emph{Proof:} In order to distinguish between the large and small time scales, we consider a smooth partition of unity $1=\sum\limits_{\gamma\in \mathbb{Z}} \chi(t-2\pi\gamma)$ with $\chi$ supported in $[-2\pi,2\pi]$. We also denote by $h_{\alpha}(t)=\chi(t) h(2\pi\alpha+t)$. By using the semigroup property and the unitarity of $e^{it\Delta_{\mathbb{R}^2\times \mathbb{T}^2}}$ we can obtain:

\begin{align*}
&\|\int_{s\in \mathbb{R}} e^{-is\Delta_{\mathbb{R}^2\times \mathbb{T}^2}} P_{\leq N} h(x,y,s)ds\|^2_{L_{x,y}^2(\mathbb{R}^2\times \mathbb{T}^2)}\\
&=\int_{s,t\in \mathbb{R}} \langle e^{-is\Delta_{\mathbb{R}^2\times \mathbb{T}^2}} P_{\leq N} h(s),e^{-it\Delta_{\mathbb{R}^2\times \mathbb{T}^2}} P_{\leq N} h(t) \rangle_{L_{x,y}^2(\mathbb{R}^2\times \mathbb{T}^2)\times L_{x,y}^2(\mathbb{R}^2\times \mathbb{T}^2)}dsdt  \\
&=\sum_{\alpha,\beta} \int_{s,t\in \mathbb{R}} \langle \chi(s-2\pi\alpha) e^{-is\Delta_{\mathbb{R}^2\times \mathbb{T}^2}} P_{\leq N} h(s),\chi(s-2\pi\beta) e^{-is\Delta_{\mathbb{R}^2\times \mathbb{T}^2}} P_{\leq N} h(t)\rangle_{L_{x,y}^2(\mathbb{R}^2\times \mathbb{T}^2)\times L_{x,y}^2(\mathbb{R}^2\times \mathbb{T}^2)}dsdt \\
&=\sum_{\alpha,\beta} \int_{s,t\in [-2\pi,2\pi]}\langle e^{-i(2\pi(\alpha-\beta)+s))\Delta_{\mathbb{R}^2\times \mathbb{T}^2}} P_{\leq N}h_{\alpha}(s),e^{-it\Delta_{\mathbb{R}^2\times \mathbb{T}^2}} P_{\leq N}h_{\beta}(t)\rangle_{L_{x,y}^2 \times L_{x,y}^2} dsdt \\
&=\sigma_{d}+\sigma_{nd}. 
\end{align*}
\noindent Here we have,
\[
\sigma_{d}=\sum_{\alpha\in \mathbb{Z},| \gamma| \leq 9} \int_{s,t\in \mathbb{R}} \langle e^{-i(s-2\pi\gamma)\Delta_{\mathbb{R}^2\times \mathbb{T}^2}} P_{\leq N}h_{\alpha}(s),e^{-it\Delta_{\mathbb{R}^2\times \mathbb{T}^2}} P_{\leq N}h_{\alpha+\gamma}(t) \rangle_{L_{x,y}^2 \times L_{x,y}^2} dsdt.
\]
\[
\sigma_{nd}=\sum_{\alpha,\gamma \in \mathbb{Z},| \gamma| > 10} \int_{s,t\in \mathbb{R}} \langle e^{-i(s-2\pi\gamma)\Delta_{\mathbb{R}^2\times \mathbb{T}^2}} P_{\leq N}h_{\alpha}(s),e^{-it\Delta_{\mathbb{R}^2\times \mathbb{T}^2}} P_{\leq N}h_{\alpha+\gamma}(t) \rangle _{L_{x,y}^2 \times L_{x,y}^2} dsdt.
\]
\noindent Here `$d$' is short for `diagonal' and `$nd$' is short for `non-diagonal'. Now, we will estimate the diagonal part and the non-diagonal part separately by using different methods as follows. \vspace{3mm}

\noindent \emph{For the diagonal part:} First we need a local-in-time $L^p$ estimate as follows:
\begin{lemma} 
\noindent Let $p_1=\frac{10}{3}$, then for any $p>p_1$, $N\geq 1$,and $f\in L^2(\mathbb{R}^2\times \mathbb{T}^2)$,
\begin{equation}
||e^{it\Delta}P_N f||_{L^p(\mathbb{R}^2\times \mathbb{T}^2\times [0,2\pi])}\lesssim_p N^{2-\frac{6}{p}}||f||_{L^2(\mathbb{R}^2\times \mathbb{T}^2)}.
\end{equation}
\end{lemma}

\noindent We will give the proof of Lemma 3.3 in the Appendix (Section 8).\vspace{3mm}

\noindent According to the estimate (3.4) above, by duality we have
\begin{equation}
\|\int_{s\in \mathbb{R}} e^{-is\Delta_{\mathbb{R}^2\times \mathbb{T}^2}} P_{\leq N} h(s) ds\|_{L_{x,y}^2(\mathbb{R}^2\times \mathbb{T}^2)} \lesssim N^{2-\frac{6}{p}} \|h\|_{L_{x,y,t}^{p^{'}}(\mathbb{R}^2\times \mathbb{T}^2\times [-2\pi,2\pi])} 
\end{equation}
\noindent where $h$ is supported in $[-2\pi,2\pi]$. And consequently,

\begin{equation}\label{equation}
\aligned
\sigma_{d} &=\sum_{\alpha \in \mathbb{Z} ,|\gamma| \leq 9} \int_{s,t\in \mathbb{R}} \langle e^{-i(s-2\pi\gamma)\Delta_{\mathbb{R}^2\times \mathbb{T}^2}} P_{\leq N}h_{\alpha}(s),e^{-it\Delta_{\mathbb{R}^2\times \mathbb{T}^2}} P_{\leq N}h_{\alpha+\gamma}(t)\rangle_{L_{x,y}^2 \times L_{x,y}^2} dsdt  \\
&\leq \sum_{\alpha\in \mathbb{Z} ,|\gamma| \leq 9} \|\int_{s\in \mathbb{R}} e^{-is\Delta_{\mathbb{R}^2\times \mathbb{T}^2}} P_{\leq N} h_{\alpha}(2\pi\gamma+s)ds\|_{L_{x,y}^2} \|\int_{s\in \mathbb{R}} e^{-is\Delta_{\mathbb{R}^2\times \mathbb{T}^2}} P_{\leq N} h_{\alpha+\gamma}(s)ds\|_{L_{x,y}^2} \\
&\lesssim N^{2(2-\frac{6}{p})} \sum_{\alpha} \|h_{\alpha}\|_{L_{x,y,t}^{p^{'}}(\mathbb{R}^2\times \mathbb{T}^2 \times [-2\pi,2\pi])}^2.  
\endaligned
\end{equation}

\noindent This finishes the estimate for the diagonal part.\vspace{3mm}

\noindent \emph{For the non-diagonal part:} We need a lemma (Lemma 3.4) that we will prove shortly and we can apply it to estimate the non-diagonal part by using H$\ddot{o}$lder's inequality and the discrete Hardy-Sobolev inequality as below:
\begin{equation}\label{equation}
\aligned
\sigma_{nd}&=\sum_{\alpha,\gamma \in \mathbb{Z},| \gamma| > 10} \int_{t \in \mathbb{R}} \langle \int_{s \in \mathbb{R}}e^{-i(s-2\pi\gamma)\Delta_{\mathbb{R}^2\times \mathbb{T}^2}} P_{\leq N}h_{\alpha}(s)ds,e^{-it\Delta_{\mathbb{R}^2\times \mathbb{T}^2}} P_{\leq N}h_{\alpha+\gamma}(t) \rangle_{L_{x,y}^2 \times L_{x,y}^2} dt \\
&\lesssim N^{2-\frac{16}{3p}} \sum_{\alpha,\gamma \in \mathbb{Z},| \gamma| > 3} |\gamma|^{\frac{2}{p}-1} \|h_{\alpha}\|_{L_{x,y,t}^{p^{'}}}\|h_{\alpha+\gamma}\|_{L_{x,y,t}^{p^{'}}}    \\
&\lesssim N^{2-\frac{16}{3p}}\|h_{\alpha}\|_{l_{\alpha}^{q^{'}} L_{x,y,t}^{p^{'}} (\mathbb{R}^2\times \mathbb{T}^2\times[-2\pi,2\pi])}^2 .
\endaligned
\end{equation}
\begin{lemma} Suppose $\gamma\in \mathbb{Z}$ satisfies $|\gamma|\geq 3$ and that $p> \frac{10}{3}$. For any function $h\in L_{x,y,s}^{p^{'}}(\mathbb{R}^2\times \mathbb{T}^2 \times [-2\pi,2\pi])$, the following inequality holds:
\[
\|\int_{s\in \mathbb{R}} \chi(s) e^{i(t-s+2\pi\gamma)\Delta_{\mathbb{R}^2\times \mathbb{T}^2}} P_{\leq N} h(s)ds\|_{L_{x,y,t}^{p}(\mathbb{R}^2\times \mathbb{T}^2 \times [-2\pi,2\pi])} \lesssim |\gamma|^{\frac{2}{p}-1} N^{2-\frac{16}{3p}} \|h\|_{L_{x,y,t}^{p^{'}}(\mathbb{R}^2\times \mathbb{T}^2 \times [-2\pi,2\pi])}.
\]\end{lemma}

\noindent \emph{Proof:}  The proof of this lemma is similar to [14, Lemma 3.3] by using Hardy-Littlewood circle method. The main idea of the proof is to study the Kernel $K_{N,\gamma}$, use a partition and decompose the corresponding index set into three parts and estimate over the three parts separately. \vspace{3mm}

\noindent Without loss of generality, we assume that:
\[h=\chi(s)P_{\leq N}h,\quad ||h||_{L^{p^{'}}(\mathbb{R}^2\times \mathbb{T}^2\times [-2\pi,2\pi])}=1
\]
\noindent and we define:
\[g(x,y,s)=\int_{s\in \mathbb{R}}e^{i(t-s+2\pi \gamma)\Delta_{\mathbb{R}^2\times \mathbb{T}^2}} h(x,y,s) ds.
\]
\noindent Also the Kernel is defined as follows: 
\begin{equation}\label{equation}
\aligned
K_N(x,y,t)&=\sum_{k \in \mathbb{Z}^2}\int_{\mathbb{R}^2_{\xi}}[\eta^1_{\leq N}(\xi_1)\eta^1_{\leq N}(\xi_2)]^2[\eta^1_{\leq N}(k_1)\eta^1_{\leq N}(k_2) ]^2e^{i[x\cdot \xi+y\cdot k+t(|k|^2+|\xi|^2)]} d\xi\\
&=[\int_{\mathbb{R}^2_{\xi}}[\eta^1_{\leq N}(\xi_1) \eta^1_{\leq N}(\xi_2)]^2e^{i[x\cdot \xi+t|\xi|^2]}]\cdot[\sum_{k \in \mathbb{Z}^2}[\eta^1_{\leq N}(k_1)\eta^1_{\leq N}(k_2) ]^2e^{i[y\cdot k+t|k|^2]}]\\
&=K^{\mathbb{R}^2}_N(x,t)\bigotimes K^{\mathbb{T}^2}_N(x,t).
\endaligned
\end{equation}
\noindent And we define $K_{N,\gamma}(x,y,t):=K_N(x,y,2\pi\gamma+t)$ so that we have $g(x,y,t)=K_{N,\gamma}*h$. Notice that a remarkable difference is the non-stationary phase estimate because of the dimension, we have:
\[
\|K_{N,\gamma}\|_{L_{x,y,t}^{\infty}} \lesssim |\gamma|^{-1} N^2
\]
\noindent instead of
\[
\|K_{N,\gamma}\|_{L_{x,y,t}^{\infty}} \lesssim |\gamma|^{-\frac{1}{2}} N^2.
\]
\noindent And
\[||\mathcal{F}_{x,y,t} K_{N,\gamma}||_{L_{\xi,k,\tau}^{\infty}} \lesssim 1
\]
still holds.\vspace{3mm}

\noindent For $\alpha$ a dyadic number, we define $g^{\alpha}(x,y,t)=\alpha^{-1}g(x,y,t)1_{\{\frac{\alpha}{2} \leq |g|\leq \alpha\}}$ which has modulus in $[\frac{1}{2},1]$. We define $h^{\beta}$ similarly for $\beta\in 2^{\mathbb{Z}}$. And we have the following decomposition:
\begin{equation}\label{equation}
\aligned
||g||^p_{L^p_{x,y,t}}&=\langle |g|^{p-2}g,g \rangle \\
&=\sum_{\alpha,\beta}\alpha^{p-1} \beta \langle |g^{\alpha}|^{p-2}g^{\alpha},K_{N,\gamma}*h^{\beta} \rangle \\
&=[\sum_{\mathcal{S}_1}+\sum_{\mathcal{S}_2}+\sum_{\mathcal{S}_3}]\alpha^{p-1} \beta \langle |g^{\alpha}|^{p-2}g^{\alpha},K_{N,\gamma}*h^{\beta} \rangle\\
&=\sum_{1}+\sum_{2}+\sum_{3},
\endaligned
\end{equation}
\noindent where $\mathcal{S}_1$, $\mathcal{S}_2$, $\mathcal{S}_3$ are three index sets. And similarly in this case we have the following decomposition:
\begin{enumerate}
\item $\mathcal{S}_1=\{(\alpha,\beta):C|\gamma|^{-1}N^2\leq \alpha (\beta)^{{p^{'}}-1}\}$, 
\item $\mathcal{S}_2=\{(\alpha,\beta):\alpha (\beta)^{{p^{'}}-1}\leq CN|\gamma|^{-1}\}$,
\item $\mathcal{S}_3=\{(\alpha,\beta):CN|\gamma|^{-1} \leq \alpha (\beta)^{{p^{'}}-1} \leq C|\gamma|^{-1}N^2\}$
\end{enumerate}
\noindent for $C$ a large constant to be decided later. For fixed $\alpha,\beta$, we will decompose $K_{N,\gamma}=K^1_{N,\gamma;\alpha,\beta}+K^2_{N,\gamma;\alpha,\beta}$ and estimate them as follows.
\begin{equation}
\langle |g^{\alpha}|^{p-2}g^{\alpha},K^{1}_{N,\gamma}*h^{\beta} \rangle \lesssim ||K^1_{N,\gamma;\alpha,\beta}||_{L^{\infty}_{x,y,t}}||g^{\alpha}||_{L^1}||h^{\beta}||_{L^1},
\end{equation}
\begin{equation}
\langle |g^{\alpha}|^{p-2}g^{\alpha},K^{2}_{N,\gamma}*h^{\beta} \rangle \lesssim ||\mathcal{F}_{x,y,t} K^2_{N,\gamma;\alpha,\beta}||_{L^{\infty}_{\xi,k,\tau}}||g^{\alpha}||_{L^2}||h^{\beta}||_{L^2}.
\end{equation}
\noindent Then we can estimate the three parts as in [13]. Eventually for the conclusion, there is one difference as follows:
\[||g||^p_{L^p_{x,y,t}}\lesssim_C ||g||^{\frac{p}{2}}_{L^p_{x,y,t}}|\gamma|^{\frac{2-p}{2}}\textmd{max}(N^{p-4+\epsilon},N^{\frac{p-2}{2}})\lesssim ||g||^{\frac{p}{2}}_{L^p_{x,y,t}}|\gamma|^{\frac{2-p}{2}}N^{p-\frac{8}{3}}
\]
\noindent if $p>\frac{10}{3}$. The rest follows as in [13] so we omitted.\vspace{3mm}

\noindent This finishes the estimate for the non-diagonal part.

\section{Local well-posedness and small-data scattering}
\noindent Recall ``Z-norm" (scattering norm) is as follows
\[
\|u\|_{Z(I)}=(\sum_{N\geq 1}N^{2} \|1_I(t) P_N u\|_{L_{x,y,t}^4(\mathbb{R}^2\times \mathbb{T}^2 \times \mathbb{R})}^4)^\frac{1}{4}.
\]
\noindent Now for convenience, we define ``$Z^{'}$-norm" which is a mixture of $Z$-norm and $X^{1}$-norm as follows
\begin{equation}\label{equation}
\|u\|_{Z^{'}(I)}=\|u\|_{Z(I)}^{\frac{3}{4}} \|u\|_{X^1(I)}^{\frac{1}{4}}.
\end{equation}
\begin{lemma}[Bilinear Estimate] Suppose that $ u_i=P_{N_i}u$, for $i=1,2$ satisfying $N_1\geq N_2$. There exists $\delta$ such that the following estimate holds for any interval $I\in \mathbb{R}$:   
\begin{equation}\label{equation}
\|u_1 u_2\|_{L_{x,t}^2(\mathbb{R}^2\times \mathbb{T}^2\times I)} \lesssim (\frac{N_2}{N_1}+\frac{1}{N_2})^\delta \|u_1\|_{Y^0(I)} \|u_2\|_{Z^{'}(I)} .
\end{equation}
\end{lemma}
\noindent \emph{Proof:} Without loss of generality, we can assume that $I=\mathbb{R}$. On one hand, we need the following estimate which follows from [16, Proposition 2.8],
\begin{equation}\label{equation}
\|u_1 u_2\|_{L^2(\mathbb{R}^2\times \mathbb{T}^2\times \mathbb{R})} \lesssim N_2(\frac{N_2}{N_1}+\frac{1}{N_2})^\delta \|u_1\|_{Y^{0}(\mathbb{R})} \|u_2\|_{Y^{0}(\mathbb{R})} .
\end{equation}  
\noindent And it suffices to prove the following estimate, if it is true then we can just combine the two inequalities (noticing the definition of $Z^{'}$-norm) and we will get the lemma proved.
\begin{equation}\label{equation}
\|u_1 u_2\|_{L^2(\mathbb{R}^2\times \mathbb{T}^2\times \mathbb{R})} \lesssim  \|u_1\|_{Y^{0}(\mathbb{R})} \|u_2\|_{Z(\mathbb{R})} .
\end{equation}  
\noindent We first notice that, by orthogonality considerations, we may replace $u_1$ by $P_C u_1$ where $C$ is a cube of dimension $N_2$. By using H$\ddot{o}$lder's inequality, we have,
\begin{align*}
\|(P_{C}u_{1})u_2\|_{L^2_{x,y,t}} &\lesssim \|P_{C}u_{1}\|_{L_{x,y,t}^4(\mathbb{R}^2\times \mathbb{T}^2 \times \mathbb{R})} \|u_2\|_{L_{x,y,t}^4(\mathbb{R}^2\times \mathbb{T}^2\times \mathbb{R})}      \\
&\lesssim N_2^{\frac{1}{2}}\|P_{C}u_{1}\|_{U_{\Delta}^4} \|u_2\|_{L_{x,y,t}^4(\mathbb{R}^2\times \mathbb{T}^2\times \mathbb{R})} \\
&\lesssim \|P_{C}u_1\|_{Y^0} \|u_2\|_{Z(\mathbb{R})}.
\end{align*}

\noindent Here we have used some properties of the function spaces and another form of Strichartz inequality.\vspace{3mm}

\noindent The Strichartz inequality in another form: for $p>\frac{10}{3}$ and $q$ as in Theorem 3.1, the following estimate holds for any time interval $I\subset \mathbb{R}$ and every cube $Q\subset \mathbb{R}^4$ of size $N$:
\begin{equation}\label{equation}
\|1_{I}(t) P_{Q}u\|_{l_{\gamma}^q L_{x,y,t}^p} \lesssim N^{2-\frac{6}{p}} \|u\|_{U_{\Delta}^{min(p,q)}}(I;L^2(\mathbb{R}^2\times \mathbb{T}^2)). 
\end{equation}

\noindent By using the properties of atomic spaces, it follows from the Strichartz inequality straightly.
\begin{lemma} [Nonlinear Estimate] For $u_i \in X^1(I)$, $i=1,2,3$. There holds that
\begin{equation}\label{equation}
\|\widetilde{u}_1 \widetilde{u}_2 \widetilde{u}_3 \|_{N(I)} \leq \sum_{(i,j,k)=(1,2,3) } \|u_i\|_{X^{1}(I)} \|u_j\|_{Z^{'}(I)} \|u_k\|_{Z^{'}(I)} 
\end{equation} where $\tilde{u}_i$ is either $u_i$ or $\bar{u}_i$.
\end{lemma}
\noindent \emph{Proof:} It suffices to prove the following estimate: (Without loss of generality, let $I = \mathbb{R}$)
\begin{equation}\label{equation}
\|\sum_{K\geq 1}P_{K} u_1 \prod_{i=2,3}P_{\leq CK} \tilde{u}_{i}\|_{N(\mathbb{R})}\lesssim_C \|u_1\|_{X^{1}(\mathbb{R})} \|u_2\|_{Z^{'}(\mathbb{R})} \|u_3\|_{Z^{'}(\mathbb{R})} .
\end{equation}
\noindent It suffices to prove for any $u_0\in Y^{-1}$ and $||u_0||_{Y^{-1}}\leq 1$ (By using Theorem 2.1)
\begin{equation}\label{equation}
\sum_{N_{1}} |\int_{\mathbb{R}^2\times \mathbb{T}^2\times \mathbb{R}} \bar{u}_0 P_{N_1}u_1 \prod_{i=2,3}(P_{\leq CN_{1}} \tilde{u}_{i}) dx dy dt| \leq \|u_0\|_{Y^{-1}} \|u_1\|_{X^1(\mathbb{R})} \|u_2\|_{Z^{'}(\mathbb{R})} \|u_3\|_{Z^{'}(\mathbb{R})}  .
\end{equation}
\noindent Now we split them as follows, let $u_i=\sum\limits_{N_i\geq 1}P_{N_i} u_i$, $i=0,1,2,3$, denoting $u_j^{N_j}= P_{N_j} u_j$ and then the estimate would follow from the following bound:
\begin{equation}\label{equation}
\sum_{S(N_0,N_1,N_2,N_3)} |\int u_0^{N_0} u_1^{N_1} u_2^{N_2} u_3^{N_3}dxdydt| \lesssim \|u_0\|_{Y^{-1}}\|u_1\|_{X^{1}}\|u_2\|_{Z^{'}}\|u_3\|_{Z^{'}} .
\end{equation}
\noindent Here we have set of index $S$ to be $\{(N_0,N_1,N_2,N_3) :  N_1 \sim max(N_2,N_0)\geq N_2\geq N_3 \}$ and we split $S$ into the disjoint union of $S_1$ and $S_2$ and $S_1$ is for the elements in $S$ that satisfy $N_1 \sim N_0$ and $S_2$ is for the elements in $S$ that satisfy $N_1 \sim N_2$. And we will estimate $S_1$ and $S_2$ separately. We omit the  $S_2$ part since the estimate is similar.\vspace{3mm}

\noindent By using bilinear estimate (4.2), some basic inequalities and the properties of function spaces, we have, for term in $S_1$:
\begin{align*}
|\int u_0^{N_0} u_1^{N_1} u_2^{N_2} u_3^{N_3}dxdydt|&\leq \|u_0^{N_0} u_2^{N_2}\|_{L^2} \|u_1^{N_1} u_3^{N_3}\|_{L^2} \\
&\leq (\frac{N_2}{N_0}+\frac{1}{N_2})^{\delta} (\frac{N_3}{N_1}+\frac{1}{N_3})^{\delta} \|u_0^{N_0}\|_{Y^{0}(\mathbb{R})} \|u_1^{N_1}\|_{Y^{0}(\mathbb{R})} \|u_2^{N_2}\|_{Z^{'}(\mathbb{R})} \|u_3^{N_3}\|_{Z^{'}(\mathbb{R})}.
\end{align*}
\noindent By using Cauchy-Schwarz inequality, the sum of the terms in $S_1$:
\begin{align*}
S_1&\lesssim \sum_{N_1\sim N_0} (\frac{N_2}{N_0}+\frac{1}{N_2})^{\delta} (\frac{N_3}{N_1}+\frac{1}{N_3})^{\delta} \|u_0^{N_0}\|_{Y^{0}(\mathbb{R})} \|u_1^{N_1}\|_{Y^{0}(\mathbb{R})} \|u_2^{N_2}\|_{Z^{'}(\mathbb{R})} \|u_3^{N_3}\|_{Z^{'}(\mathbb{R})} \\
&\lesssim (\sum_{N_1\sim N_0} \frac{N_0}{N_1} \|u_0^{N_0}\|_{Y^{-1}(\mathbb{R})} \|u_1^{N_1}\|_{Y^{1}(\mathbb{R})}) \|u_2\|_{Z^{'}(\mathbb{R})} \|u_3\|_{Z^{'}(\mathbb{R})} \\
&\lesssim \|u_0\|_{Y^{-1}(\mathbb{R})} \|u_1\|_{X^{1}(\mathbb{R})} \|u_2\|_{Z^{'}(\mathbb{R})} \|u_3\|_{Z^{'}(\mathbb{R})}.
\end{align*}
\noindent This finishes Lemma 4.2.
\begin{theorem}\label{}[Local Well-posedness] Let $E > 0$ and $\|u_0\|_{ H^1(\mathbb{R}^2 \times \mathbb{T}^2) }<E$,
then there exists $\delta_0=\delta_0(E)>0$ such that if 
\[\|e^{it\Delta}u_0\|_{Z(I)}< \delta \] 
\noindent for some $ \delta \leq \delta_0 $, $0\in I$. Then there exists a unique strong solution $u\in X_{c}^{1}(I)$ satisfying $u(0)=u_0$ and we can get an estimate,
\begin{equation}\label{equation}
\|u(t)-e^{it\Delta_{\mathbb{R}^2\times \mathbb{T}^2}}u_0\|_{X^1(I)}\leq (E\delta)^{\frac{3}{2}} .
\end{equation}  
\end{theorem}
\noindent \emph{Remark.} Observe that if $u \in X^1_c(\mathbb{R})$, then $u$ scatters as $t \rightarrow \pm\infty $ as in (1.3). Also, if $E$ is small enough, $I$ can be taken to be $\mathbb{R}$ which proves Theorem 1.1.\vspace{3mm}

\noindent \emph{Proof:} First, we consider a mapping defined as follows,
\[
\Phi(u)=e^{it\Delta} u_0-\int_{0}^{t}{e^{i(t-s)\Delta} |u(s)|^2 u(s)ds}.
\]
\noindent And we define a set $ B=\{u\in X_{c}^{1}(I) :   \|u\|_{X^1(I)}\leq 2E  \textmd{ and }  \|u\|_{Z(I)}\leq 2\delta \}$. Now we will check two properties of $\Phi$: 1. $\Phi$ maps $B$ to $B$. 2. $\Phi$ is a contraction mapping. \vspace{3mm}

\noindent 1. For $ u \in B$, we can use the nonlinear estimate in Lemma 4.2 and let $\delta \leq 1$ and small enough to make $E^3 \delta$ small enough, we have:                
\[
\|\Phi(u)\|_{X^1(I)} \leq \|e^{it\Delta} u_0\|_{X^1(I)}+\||u|^2 u\|_{N(I)} \leq E+CE^{\frac{3}{2}} \delta^{\frac{3}{2}} \leq 2E,
\]
\[
\|\Phi(u)\|_{Z(I)} \leq \|e^{it\Delta} u_0\|_{Z(I)}+\||u|^2 u\|_{N(I)} \leq \delta+CE^{\frac{3}{2}}\delta^{\frac{3}{2}} \leq 2\delta.
\]
\noindent 2.\begin{align*}
\|\Phi(u)-\Phi(v)\|_{X^{1}(I)} &\lesssim \|u-v\|_{X^1(I)} (\|u\|_{X^1(I)}+\|v\|_{X^1(I)})
(\|u\|_{Z^{'}(I)}+\|v\|_{Z^{'}(I)})    \\
&\leq  C\|u-v\|_{X^1(I)} E^{\frac{5}{4}} \delta^{\frac{3}{4}}  \\
&\leq  C\|u-v\|_{X^1(I)} [(E^3\delta)^\frac{5}{12} \delta^\frac{1}{3}] \\
&\leq C\frac{1}{2}\|u-v\|_{X^1(I)}.
\end{align*}

\noindent Thus the result now follows from the Picard's fixed point argument.
\begin{theorem}\label{theorem}
[Controlling Norm] Let $u\in X_{c,loc}^{1}(I)$ be a strong solution on $I\in \mathbb{R} $ satisfying 
\begin{equation}\label{equation}
\|u\|_{Z(I)}< \infty .
\end{equation}
\noindent Then we have two conclusions, \vspace{3mm}

\noindent (1) If $I$ is finite, then u can be extended as a strong solution in $X_{c,loc}^{1}(I^{'})$ on a strictly larger interval $I^{'}$, $I\subsetneq I^{'}\subset \mathbb{R}$. In particular, if u blows up in finite time, then the $Z$-norm of $u$ has to blow up.  \vspace{3mm}

\noindent (2) If $I$ is infinite, then $ u\in X_{c}^{1}(I)$.     
\end{theorem}
\noindent \emph{Proof:} Without loss of generality, for the finite case we can assume $I=[0,T)$ and we want to extend it to $[0,T+v)$ for some $v>0$. Denoting $E=\sup\limits_{I}\|u(t)\|_{H^{1}(\mathbb{R}^2 \times \mathbb{T}^2)}$ and using the time-divisibility of `$Z$-norm', there exists $T_1$ such that $T-1<T_1<T$ such that
\[\|u\|_{Z([T_{1},T])} \leq \epsilon,\]
\noindent where $\epsilon$ is to be decided. This allows to conclude:
\[\|u(t)-e^{i(t-T_{1})\Delta } u(T_1) \|_{X^1([T_{1},T])} \lesssim \|u\|_{X^1([T_1,T))}^{\frac{3}{2}} \|u\|_{Z([T_1,T))}^{\frac{3}{2}}\leq C\epsilon^{\frac{3}{2}}  \|u\|_{X^1([T_1,T))}^{\frac{3}{2}}.\]
\noindent By bootstrap argument, we get,

\[ \|u\|_{X^1([T_1,T))}^{\frac{3}{2}} \lesssim E.  \]

\noindent If $\epsilon$ is small enough and, making $\epsilon$ possibly smaller, we have,
\begin{align*}
\|e^{i(t-T_1)\Delta} u(T_1)\|_{Z([T_1,T))} &\leq \|u\|_{Z([T_1,T))} +\|e^{i(t-T_1)\Delta} u(T_1)-u(t)\|_{Z([T_1,T))}   \\
&\leq \epsilon+\|e^{i(t-T_1)\Delta} u(T_1)-u(t)\|_{X^{1}([T_1,T))}   \\
&\leq \epsilon+C^{'}\epsilon^{\frac{3}{2}} E^{\frac{3}{2}}    \\
&\leq \frac{3}{4} \delta_{0}(E).
\end{align*}

\noindent Notice that we can let $\epsilon$ small enough s.t. $\epsilon<\frac{1}{4}$ and $\epsilon E<(\frac{1}{2} \delta_{0}(E))^{\frac{2}{3}} .$\vspace{3mm}

\noindent This allows to find an interval $[T_1,T+v]$ for which:

\[ \|e^{i(t-T_1)\Delta}u(T_1)\|_{Z([T_1,T+v])}<\delta_0 .  \]
\noindent That finishes the proof by using the Theorem 4.3 (local well-posedness). Using the symmetries of the equation, the above argument also covers the case when $I$ is an arbitrary bounded interval.\vspace{3mm}

\noindent For the infinite case, without loss of generality, it is enough to consider the case $I=(a,\infty)$. Choosing $T$ to be large enough so that
\[||u||_{Z([T,\infty))} \leq \epsilon
\]
\noindent we get that for any $T^{'}>T$ : 
\[\|u(t)-e^{i(t-T)\Delta } u(T) \|_{X^1([T,T^{'}))} \lesssim \|u\|_{X^1([T,T^{'}))}^{\frac{3}{2}} \|u\|_{Z([T_1,T))}^{\frac{3}{2}}\leq C\epsilon^{\frac{3}{2}}  \|u\|_{X^1([T,T^{'}))}^{\frac{3}{2}}\]

\noindent which gives that $||u||_{X^1([T,T^{'}))} \lesssim E$ for any $T^{'}>T$ and we have
\[||e^{i(t-T)\Delta}u(T)||_{Z([T,\infty))} \leq 2\epsilon \leq \delta_0(E)
\]
\noindent if $\epsilon$ small enough. Now the result follows from Theorem 4.3. 
\begin{theorem}\label{theorem}[Stability Theory] Let $ I\in \mathbb{R}$ be an interval, and let $\widetilde{u} \in X^1(I) $ solve the approximate solution,
\begin{equation}\label{equation}
(i\partial_{t}+\Delta_{\mathbb{R}^2 \times \mathbb{T}^2} )\widetilde{u}=\rho |\widetilde{u}|^2\widetilde{u}+e \textmd{ and }  \rho\in [0,1]. 
\end{equation}
\noindent Assume that:
\begin{equation}\label{equation}
\|\widetilde{u}\|_{Z(I)}+\|\widetilde{u}\|_{L_{t}^{\infty}(I,H^{1}(\mathbb{R}^2\times \mathbb{T}^2))} \leq M.
\end{equation}
\noindent There exists $\epsilon_{0}=\epsilon_{0}(M) \in (0,1] $ such that if for some $t_{0} \in I$:
\begin{equation}\label{equation}
\|\widetilde{u}(t_{0})-u_{0}\|_{H^{1}(\mathbb{R}^2\times \mathbb{T}^2)}+\|e\|_{N(I)} \leq \epsilon < \epsilon_{0},
\end{equation}
\noindent then there exists a solution $u(t)$ to the exact equation:
\begin{equation}\label{equation}
(i\partial_{t}+\Delta_{\mathbb{R}^2 \times \mathbb{T}^2} )u=|u|^2u \end{equation}
\noindent with initial data $u_0$ that satisfies 
\begin{equation}\label{equation} \|u\|_{X^1(I)}+\|\widetilde{u}\|_{X^1(I)} \leq C(M),   \quad  \quad  \|u-\widetilde{u}\|_{X^1(I)} \leq C(M)\epsilon. \end{equation}
\end{theorem}
\noindent \emph{Proof:} The proof is similar to the proof of [13, Proposition 4.7]. The proof relies heavily on the estimate of Lemma 4.2 (nonlinear estimate) and the trick of division of the intervals.\vspace{3mm}

\noindent First, we consider for an interval $J\in I$ s.t. $\|\widetilde{u}\|_{Z(J)} \leq \epsilon $ (That is the additional smallness assumption, $\epsilon$ is to be decided). We will prove the theorem under the assumption first.\vspace{3mm}

\noindent Then by local existing argument for the approximate equation, there exists $ \delta_{1}(M)$ such that if 
\[ \|e^{i(t-t_{*})\Delta} \widetilde{u}(t_{*})\|_{Z(J)}+\|e\|_{N(J)} \leq \delta_{1} \]
for some $t_{*}\in J$, then $\widetilde{u}\in X^1(J)$ is unique and satisfies:
\[ \|\widetilde{u}- e^{i(t-t_{*})\Delta} \widetilde{u}(t_{*})\|_{X^1(J)} \leq C\|\widetilde{u}\|_{X^1(J)}^{\frac{3}{2}} \|\widetilde{u}\|_{Z(J)}^{\frac{3}{2}}+\|e\|_{N(J)}. \]

\noindent We can conclude
\[ \|\widetilde{u}\|_{X^1(J)} \lesssim M+1\quad \textmd{and} \quad \|e^{i(t-t_{*})\Delta} \widetilde{u}(t_{*})\|_{Z(J)} \lesssim \epsilon.\]\noindent if $\epsilon<\epsilon_1(M)$ is small enough.\vspace{3mm}

\noindent Second, let us estimate the difference of the solutions. Consider solution $u$ with initial data $u_{*}$ satisfying $\|u_{*}-\widetilde{u}(t_{*})\|_{H^{1}} \leq \epsilon $ and living on an interval $J_{u} \in J$ containing $t_{*}$. We want to prove the following estimate for some constant $C$ independent of $J_u$ to be specified later:
\begin{equation}\label{equation}
\|u-\widetilde{u} \|_{X^1(J_u)}\leq C\epsilon. 
\end{equation}

\noindent Let $w=u-\widetilde{u} $, then we know that $w$ satisfies: 
\[
(i\partial_t +\Delta)w=\rho(|\widetilde{u}+w|^2 (\widetilde{u}+w)-|\widetilde{u}|^2 \widetilde{u} )-e.
\]
\noindent Adopting the bootstrap hypothesis: 
\[ \|w\|_{X^1(J_u\cap[t_{*}-t,t_{*}+t]) }\leq2C\epsilon. \]
\noindent For convenience, we denote $J_u\cap[t_{*}-t,t_{*}+t) $ by $J_t$, by using nonlinear estimate, we compute:
\begin{align*}
\|w\|_{X^1(J_t)}&\lesssim \|u(t_{*})-\widetilde{u}(t_{*})\|_{H^1(\mathbb{R}^2\times \mathbb{T}^2)}+\|w\|_{X^1(J_t)} \| \widetilde{u}\|_{X^1(J_t)}  \| \widetilde{u}\|_{Z^{'}(J_t)}+\|e\|_{N(J_t)}     \\
&\lesssim \epsilon+\|w\|_{X^1(J_t)} \|\widetilde{u}\|_{X^1(J_t)}^{\frac{5}{4}}   \|\widetilde{u}\|_{Z(J_t)}^{\frac{3}{4}}   \\
&\leq C_1 \epsilon+C_1 M^{\frac{5}{4}} \epsilon^{\frac{3}{4}} \|w\|_{X^1(J_t)}.    
\end{align*}
\noindent As a result, if $\epsilon<\epsilon_1(M)$ with $\epsilon_1(M)$ small enough in terms of $M$, we conclude that $||u-\tilde{u}||_{X^1(J_t)}\leq 2C_1\epsilon$, which close the the bootstrap argument with $C=2C_1$. This finishes the proof under the additional assumption ($\|\widetilde{u}\|_{Z(J)} \leq \epsilon $). \vspace{3mm}

\noindent Now, to generalize the argument to the whole interval $I$, we split $I$ into $N=C(M,\epsilon_1(M))$ intervals $I_k=[T_k,T_{k+1})$ such that:
\[||u||_{Z(I_k)}\leq \frac{\epsilon_1(M)}{100} \quad \textmd{and} \quad ||e||_{N(I_k)}\leq \frac{\epsilon_1(M)}{100}.
\]
\noindent If $\epsilon_0(M)$ is chosen sufficiently small in terms of $N$, $M$ and $\epsilon_1(M)$, we can iterate the first part of the proof on each interval $I_k$ while keeping the condition
\[||u(T_k)-\tilde{u}(T_k)||_{H^1(\mathbb{R}^2\times \mathbb{T}^2)}+||e||_{N(I_k)}+||u||_{Z(I_k)}<\epsilon_1(M)
\]
\noindent satisfied for each $k$. This finishes the proof.

\section{Nonlinear Analysis of the Profiles}
\noindent In this section, we describe and analyze the profiles that appear in our linear and nonlinear profile decomposition. \vspace{3mm}

\noindent \textbf{5.1 Euclidean Profiles.} The Euclidean profiles define a regime where we can compare solutions of cubic NLS on $\mathbb{R}^4$ with those on $\mathbb{R}^2\times \mathbb{T}^2$. \vspace{3mm}

\noindent \emph{Remark.} We refer to [13, Section 5] and [19, Section 4] for more information. For those problems, Euclidean profiles also appear in the analysis of profile decomposition according to the structures of the corresponding equations.\vspace{3mm}

\noindent Here is the precise characterization. We fix a spherically symmetric function $\eta \in C_{0}^{\infty}(\mathbb{R}^4)$ supported in the ball of radius 2 and equal to 1 in the ball of radius 1. Given $\phi \in H^1(\mathbb{R}^4)$ and a real number $N \geq 1$, we define:
\[
Q_{N} \phi \in H^1(\mathbb{R}^4)   \quad  \quad  (Q_{N} \phi)  (x)=\eta(\frac{x}{N^{\frac{1}{2}}}) \phi(x)
\]
\begin{equation}\label{euqation}
\phi_{N}\in H^1(\mathbb{R}^4)  \quad  \quad \phi_{N}(x)=N(Q_N \phi)(Nx)
\end{equation}
\[
f_N\in H^1(\mathbb{R}^2\times \mathbb{T}^2) \quad \quad f_N(y)=\phi_N(\Psi^{-1}(y))
\]
\noindent where $\Psi$ is the identity map from the unit ball of $\mathbb{R}^4$ to $\mathbb{R}^2\times \mathbb{T}^2$. Thus $ Q_N \phi $ is a compactly supported modification of the profile $\phi$, $\phi_{N}$ is an $\dot{H^{1}}$-invariant rescaling of $Q_N \phi $, and $f_{N}$ is the function obtained by transferring $\phi_{N}$ to a neighborhood of $0$ in $\mathbb{R}^2 \times \mathbb{T}^2$. Notice that 
\[
\|f_{N}\|_{H^1(\mathbb{R}^2\times \mathbb{T}^2)}  \lesssim \|\phi\|_{\dot{H^{1}}(\mathbb{R}^4)}.
\]
\noindent Then we use the main theorem of [24] (by E. Ryckman and M. Visan) in the following form:
\begin{theorem}\label{theorem}
Assume $\psi \in \dot{H^{1}}(\mathbb{R}^4)$, then there is a unique global solution $v \in C(\mathbb{R}:\dot{H^{1}}(\mathbb{R}^4))$ of the initial-value problem 
\begin{equation}\label{euqation}
(i\partial_{t} +\Delta_{\mathbb{R}^4})v=v|v|^{2}, \quad v(0)=\psi,
\end{equation}
\noindent and
\begin{equation}\label{euqation}
\||\nabla_{\mathbb{R}^4}v|\|_{(L_{t}^{\infty}L_{x}^{2} \cap L_{t}^{2}L_{x}^{4})(\mathbb{R}^4\times \mathbb{R})} \leq \tilde{C}(E_{\mathbb{R}^4}(\psi)).
\end{equation}
\noindent Moreover this solution scatters in the sense that there exists $\psi^{\pm \infty} \in \dot{H^{1}}(\mathbb{R}^4)$ such that 
\begin{equation}\label{euqation}
\|v(t)-e^{it\Delta} \psi^{\pm \infty}\|_{\dot{H^{1}}(\mathbb{R}^4)} \to 0
\end{equation}
as $t \to \pm \infty$. Besides if $\psi \in H^{5}(\mathbb{R}^4)$, then $v \in C(\mathbb{R}:H^5(\mathbb{R}^4))$ and 
\[
\sup\limits_{t\in \mathbb{R}} \|v(t)\|_{H^{5}(\mathbb{R}^4)}  \lesssim_{\|\psi\|_{H^{5}(\mathbb{R}^4)}}   1.
\] \end{theorem}
\noindent Based on the existing result, we have:
\begin{theorem}\label{theorem}
\noindent Assume $\phi \in \dot{H^{1}}(\mathbb{R}^4)$, $T_{0} \in (0,\infty)$, and $\rho \in \{0,1\}$ are given, and define $f_{N}$ as before. Then the following conclusions hold:\vspace{3mm} 

\noindent $(1)$ There is $N_0=N_0 (\phi,T_0)$ sufficiently large such that for any $N\geq N_0$, there is a unique solution $U_{N} \in C((-T_0 N^{-2},T_0 N^{-2}); H^1(\mathbb{R}^2\times \mathbb{T}^2))$ of the initial-value problem 
\begin{equation}\label{euqation}
(i\partial_t+\Delta)U_{N}=\rho U_{N} |U_{N}|^{2} ,\quad and \quad  U_{N}(0)=f_{N}.
\end{equation}
\noindent Moreover, for any $N \geq N_{0}$, 
\begin{equation}\label{euqation}
\|U_N\|_{X^{1}(-T_0 N^{-2},T_0 N^{-2})} \lesssim_{E_{\mathbb{R}^4}(\phi)} 1.
\end{equation}
$(2)$ Assume $\epsilon_{1} \in (0,1]$ is sufficiently small (depending only on $E_{\mathbb{R}^4}(\phi)$), $\phi ^{'} \in H^{5}(\mathbb{R}^4)$, and $\|\phi-\phi^{'}\|_{\dot{H^{1}}(\mathbb{R}^4)} \leq \epsilon_{1}$. Let $v^{'} \in C(\mathbb{R}:H^5)$ denote the solution of the initial-value problem 
\[
(i\partial_t+\Delta_{\mathbb{R}^4}) v^{'}=\rho v^{'} |v^{'}|^{2}, \quad v^{'} (0)=\phi^{'}.
\] 
\noindent For $R \geq 1$ and $N \geq 10R$, we define 
\[ v_{R}^{'}(x,t)=\eta(\frac{x}{R})v^{'}(x,t)    \quad \quad (x,t) \in \mathbb{R}^{4}\times (-T_0,T_0)
\]
\begin{equation}\label{euqation}
v_{R,N}^{'}(x,t)=N v_{R}^{'}(Nx,N^2 t) \quad \quad (x,t) \in \mathbb{R}^{4}\times (-T_0 N^{-2},T_0 N^{-2})
\end{equation}
\[
V_{R,N}(y,t)=v_{R,N}^{'}(\Psi^{-1}(y),t) \quad \quad (y,t) \in \mathbb{R}^2\times \mathbb{T}^2 \times (-T_0 N^{-2},T_0 N^{-2})
\]
\noindent Then there is $R_0 \geq 1$ (depending on $T_0$ and $\phi^{'}$ and $\epsilon_1$) such that, for any $R \geq R_{0}$ and $N \geq 10R$,
\begin{equation}\label{euqation}
\limsup\limits_{N \to \infty} \|U_N-V_{R,N}\|_{X^{1}(-T_0 N^{-2},T_0 N^{-2})}  \lesssim_{E_{\mathbb{R}^4}(\phi)} \epsilon_1.
\end{equation} 
\end{theorem}
\noindent \emph{Proof:} It suffices to prove part (2). All implicit constants are allowed to depend on $\|\phi\|_{\dot{H^{1}}(\mathbb{R}^4)}$. The idea of the proof is to show that with $ R_{0}$ chosen large enough, $V_{R,N}$ is an approximate solution.   \vspace{3mm}

\noindent First, we define: 
\[ e_{R}(x,t) := (i\partial_t +\Delta_{\mathbb{R}^4})v_{R}^{'}-\rho {|v_{R}^{'}|}^{2} v_{R}^{'}.
\]
\noindent Using the fact that $\sup\limits_{t} \|v^{'}(t)\|_{H^5} \lesssim_{\|\phi^{'}\|_{H^5}} 1$, we get that:  
\[
|e_{R}(t,x)|+|\nabla_{\mathbb{R}^4} e_{R}(t,x)| \lesssim 1_{[R/2,4R]} (|v^{'}(t,x)|+|\nabla_{\mathbb{R}^4}v^{'}(t,x)|+|\Delta_{\mathbb{R}^4}v^{'}(t,x)|)
\]
\noindent which directly gives that there exists $R_0 \geq 1$ such that for all $R > R_{0}$   
\[
\lim\limits_{R\to \infty} \| |e_{R}|+|\nabla_{\mathbb{R}^4} e_{R}| \|_{L_{t}^{1} L_{x}^{2}(\mathbb{R}^4 \times (-T,T))} = 0.
\]
\noindent Letting 
\[
e_{R,N}(x,t):=(i\partial_{t}+\Delta_{\mathbb{R}^4})v_{R,N}^{'}-\rho|v_{R,N}^{'}|^{2}v_{R,N}^{'},
\]
\noindent we have that for any $R>R_{0}$ and $N \geq 1$: 
\begin{equation}\label{euqation}
 \| |e_{R,N}|+|\nabla_{\mathbb{R}^4} e_{R,N}| \|_{L_{t}^{1} L_{x}^{2}(\mathbb{R}^4 \times (-TN^{-2},TN^{-2}))} \leq  2\epsilon_{1}
\end{equation}
\noindent with $V_{R,N} $ defined on $\mathbb{R}^2 \times \mathbb{T}^2 \times (-TN^{-2},TN^{-2})$, we let 
\[ E_{R,N}(y,t)=(i\partial_{t}+\Delta_{\mathbb{R}^4})V_{R,N}-\rho|V_{R,N}|^{2}V_{R,N}=e_{R,N}(\Psi^{-1}(y),t). 
\]
\noindent For $R >R_0$ and $N \geq 10R$:
\[  \| |E_{R,N}|+|\nabla_{\mathbb{R}^4} E_{R,N}| \|_{L_{t}^{1} L_{x}^{2}(\mathbb{R}^4 \times (-TN^{-2},TN^{-2}))} \lesssim  \epsilon_{1}
\]
\noindent from which  it follows (using Theorem 2.1) that:
\[ ||E_{R,N}||_{N(-TN^{-2},TN^{-2})} \lesssim \epsilon_{1}.
\]\noindent To verify the requirements of Theorem 4.5, we use (5.3) to conclude that:
\[ ||V_{R,N}||_{L_t^{\infty}H^1(\mathbb{R}^2 \times \mathbb{T}^2 \times(-TN^{-2},TN^{-2}))} \lesssim 1 .
\]
\noindent As for the $Z$-norm control, we choose $N$ to be big enough so that $TN^{-2} \leq \frac{1}{2}$ which makes all summations in the  $Z$-norm consist of at most two terms, after which we estimate the $Z$-norm by using Littlewood-Paley theory and Sobolev embedding theorem as follows: 
\[  ||K^{\frac{1}{2}}||P_k V_{R,N}||_{L^4_{t,x}(\mathbb{R}^2\times \mathbb{T}^2 \times (-TN^{-2},TN^{-2}))}||_{l_k^4} \lesssim ||(1-\Delta)^{\frac{1}{4}}V_{R,N}||_{L^4_{t,x} } \lesssim  ||(1-\Delta)^{\frac{1}{2}}v^{'}_{R,N}||_{L^4_{t}L^{\frac{8}{3}}_{x}} \lesssim_{E(\phi)} 1.
\]

\noindent At last, we know for $R_0$ big enough and $R >R_0$, $N \geq 10R$,
\[ ||f_N-V_{R,N}(0)||_{H^1(\mathbb{R}^2 \times \mathbb{T}^2)} \lesssim ||Q_N\phi-\phi||_{\dot{H}^1(\mathbb{R}^4)}+||\phi^{'}-\phi||_{\dot{H}^1(\mathbb{R}^4)}+||\phi^{'}-V^{'}_{R}(0)||_{\dot{H}^1(\mathbb{R}^4)} \lesssim \epsilon_1.
\]
\noindent This completes the verification of the requirements of Theorem 4.5 which concludes the proof by using perturbation theory.
\begin{lemma}[Extinction Lemma]Suppose that $\phi \in \dot{H^1}(\mathbb{R}^4), \epsilon > 0$, and  $I \subset \mathbb{R} $ is an interval. Assume that
\begin{equation}\label{euqation}
||\phi||_{\dot{H^1} (\mathbb{R}^4)} \leq 1, \quad ||\nabla_{x} e^{it\Delta}\phi||_{L_t^2 L_x^4(\mathbb{R}^4 \times I)} \leq \epsilon.
\end{equation}
\noindent For $N \geq 1$, we define as before: 
\[Q_N \phi=\eta(N^{-1/2}x)\phi(x),\quad \phi_N=N(Q_N\phi)(Nx),\quad f_N(y)=\phi_N(\Psi^{-1}(y)).
\]
\noindent Then there exists $N_0=N_0(\phi,\epsilon)$ such that for any $N \geq N_0$,
\[ ||e^{it\Delta} f_{N}||_{Z(N^{-2}I)} \lesssim \epsilon.
\]
\end{lemma}
\noindent \emph{Proof:} It suffices to prove that there exists $T_0$ such that for any $N>1$:
\begin{equation}\label{euqation}
||e^{it\Delta} f_N||_{Z(\mathbb{R} \setminus (-N^{-2}T_0,N^{-2}T_0))} \lesssim \epsilon
\end{equation}
\noindent as the rest follows from Lemma 5.2 (with $\rho=0$). Without loss of generality, by limiting arguments, we may assume that $\phi \in C^{\infty}_{0}(\mathbb{R}^4)$. We have, for any $p$, 
\[ f_{N,p}(x)=\frac{1}{(2\pi)^2} \int_{\mathbb{T}^2} \phi_{N}(x,y)e^{-i\langle y,p \rangle} dy= \frac{N}{(2\pi)^2} \int_{\mathbb{R}^2} e^{-i\langle y,p \rangle} \phi(Nx,Ny) dy.
\]
\noindent And using dispersive estimate and unitarity, we have 
\begin{equation}\label{euqation} ||e^{it\Delta}P_M f_N(t)||_{L^{\infty}_{x,y}(\mathbb{R}^2\times \mathbb{T}^2)} \lesssim \sup_{x \in \mathbb{R}^2}\sum_{|p| \leq M}|e^{it\Delta_{x}}f_{N,p}(x)| \lesssim \frac{M^2}{|t|}||f_N||_{L^1_{x,y}} \lesssim \frac{M^2N^{-3}}{|t|} 
\end{equation}
\noindent and 
\begin{equation}\label{euqation} ||e^{it\Delta}P_M f_N(t)||_{L^{2}_{x,y}(\mathbb{R}^2\times \mathbb{T}^2)}= ||P_M f_N(t)||_{L^{2}_{x,y}(\mathbb{R}^2\times \mathbb{T}^2)} \lesssim M^{-l} ||(1-\Delta)^{\frac{l}{2}} \phi_N||_{L^2(\mathbb{R}^4)} \lesssim M^{-l}N^{l-1}.
\end{equation}
\noindent Then by interpolation we have (choose $l=0, 10000$):
\begin{equation}\label{equation}
\aligned
||e^{it\Delta}P_M f_N(t)||_{L^{p}_{x,y}(\mathbb{R}^2\times \mathbb{T}^2)}  &\lesssim \frac{N^{-1}}{|t|^{1-\frac{2}{p}}} [(\frac{M}{N})^{2-\frac{4}{p}-\frac{2l}{p}}] \\
&\lesssim \frac{N^{-1}}{|t|^{1-\frac{2}{p}}} min[(\frac{M}{N})^{2-\frac{4}{p}},(\frac{N}{M})^{100}].
\endaligned
\end{equation}
\noindent As a result,
\begin{equation}\label{euqation} (\sum\limits_{M}M^{2}||e^{it\Delta}P_M f_N(t)||^4_{l^{4}L^{4}_{x,y}(\mathbb{R}^2\times \mathbb{T}^2 \times (|t|\geq TN^{-2}))})^{\frac{1}{4}}  \lesssim T^{-\frac{1}{4}}.
\end{equation}

\noindent According to the definition of $Z$-norm, we finish the proof.\vspace{3mm}

\noindent We can now describe the nonlinear solutions of (1.1) corresponding to data concentrating at a point. Let $\tilde{\mathcal{F}}_e$ denote the set of renormalized Euclidean frames as follows:\vspace{3mm}

\noindent $\tilde{\mathcal{F}}_e := \{(N_k,t_k,x_k)_{k\geq 1}:N_k \in [1,\infty),x_k \in \mathbb{R}^2 \times \mathbb{T}^2 ,N_k \rightarrow \infty$, and either $t_k=0$ for any $k\geq 1$ or $\lim\limits_{k\rightarrow \infty}N_k^2 |t_k|=\infty \}$.\vspace{3mm}

\noindent Given $f \in L^2(\mathbb{R}^2 \times \mathbb{T}^2)$, $t_0 \in \mathbb{R}$, and $x_0 \in \mathbb{R}^2\times \mathbb{T}^2$, we define:
\begin{equation}\label{euqation} 
\pi_{x_0}f=f(x-x_0), \quad \Pi_{(t_0,x_0)}f=(e^{-it_0\Delta_{\mathbb{R}^2 \times \mathbb{T}^2}}f)(x-x_0)=\pi_{x_0}e^{it_0\Delta_{\mathbb{R}^2 \times \mathbb{T}^2}}f. 
\end{equation}
\noindent Also for $\phi \in \dot{H}^1(\mathbb{R}^4)$ and $N \geq 1$, we denote the function obtained in (5.1) by:
\begin{equation}
T^e_{N_k}:=N\tilde{\phi}(N\Psi^{-1}(x)) \quad where \quad \tilde{\phi}(y):=\eta(\frac{y}{N^{\frac{1}{2}}}) \phi(y)
\end{equation}
\noindent and as before observe that $T_N^e:\dot{H}^1(\mathbb{R}^4) \rightarrow H^1(\mathbb{R}^2 \times \mathbb{T}^2)$ with $||T_N^e\phi||_{H^1(\mathbb{R}^2 \times \mathbb{T}^2)} \lesssim ||\phi||_{\dot{H}^1(\mathbb{R}^4)}$.
\begin{theorem}\label{theorem}
\noindent Assume that $\mathcal{O}=(N_k,t_k,x_k)_{k} \in \tilde{\mathcal{F}}_e$, $\phi \in \dot{H^1}(\mathbb{R}^4)$, and let $U_k(0)=\Pi_{t_k,x_k}(T^e_{N_k}\phi)$: \vspace{3mm}

\noindent (1) For $k$ large enough, there is a nonlinear solution $U_k \in X^1(\mathbb{R})$ of the equation (1.1) satisfying: 
\begin{equation}\label{euqation}
||U_k||_{X^1(\mathbb{R})} \lesssim_{E_{\mathbb{R}^4}(\phi)} 1.
\end{equation}
\noindent (2) There exists a Euclidean solution $u\in C(\mathbb{R}:\dot{H}^1(\mathbb{R}^4))$ of 
\begin{equation}\label{euqation} (i\partial_t+\Delta_{\mathbb{R}^4})u=|u|^2u
\end{equation}
\noindent with scattering data $\phi^{\pm  \infty}$ defined as in (5.4) such that up to a subsequence: for any $\epsilon > 0$, there exists $T(\phi,\epsilon)$ such that for all $T \geq T(\phi,\epsilon)$ there exists $R(\phi,\epsilon,T)$ such that for all $R \geq R(\phi,\epsilon,T)$, there holds that 
\begin{equation}\label{euqation} ||U_k-\tilde{u}_k||_{X^1(\{ |t-t_k| \leq TN_k^{-2}\})} \leq \epsilon,
\end{equation}
\noindent for $k$ large enough, where 
\begin{equation}\label{equation}
(\pi_{-x_k}\tilde{u})(x,t)=N_k \eta(N_k \Psi^{-1}(x)/R)u(N_k \Psi^{-1}(x),N_k^2(t-t_k)).
\end{equation}
\noindent In addition, up to a subsequence, 
\begin{equation} ||U_k(t)-\Pi_{(t_k-t,x_k)}T^e_{N_k} \phi^{\pm  \infty}||_{X^1(\{ |t-t_k| \leq TN_k^{-2}\})} \leq \epsilon
\end{equation}
\noindent for $k$ large enough (depending on ($\phi,\epsilon,T,R$)).
\end{theorem}
\noindent \emph{Proof:} Without loss of generality, we may assume $x_k=0$ for all $k$. We first consider the case when $t_k=0$ for all $k$.\vspace{3mm}

\noindent If $u$ is the solution with initial data $\phi$, there exists a time $T_0=T_0(\phi,\epsilon)$ such that: 
\begin{equation}\label{euqation}
||\nabla u||_{L^2_t L^4_x(\mathbb{R}^4 \times \{|t| \geq T\})} \ll_{E(\phi)} \epsilon.
\end{equation}
\noindent Theorem 5.2 tells us that for any $T \geq T_0$, there is $R_0=R_0(T,\phi,\epsilon)$ such that for all $R \geq R_0$ and $N_k \geq 10R$, it holds that: 
\[||U_k-\tilde{u}_k||_{X^1(\{|t|\leq TN_k^{-2}\})} \leq \epsilon.
\]\noindent Equation (5.23) along with Lemma 5.3 imply that $e^{it\Delta}U_k(\pm TN_k^{-2})$ is sufficiently small in $Z(\{\pm t\geq TN_k^{-2}\})$ thus guaranteeing that $U_k$ extends to a global solution in $X^1(\mathbb{R})$ that satisfies (5.22).\vspace{3mm}

\noindent We now turn to the other case when $ N_k^2|t_k| \rightarrow  \infty$. For definiteness, we assume that $N^2_kt_k \rightarrow \infty$ and $u$ be the solution to (5.20) satisfying:
\[||\nabla (u(t)-e^{it\Delta} \phi)||_{L^2(\mathbb{R}^4)} \rightarrow 0
\]
\noindent as $t \rightarrow -\infty$. Let $\tilde{\phi}=u(0)$ and let $V_k$ be the solution of (1.1) with initial data $T^e_{N_k} \tilde{\phi}$.\vspace{3mm}

\noindent Applying the first case of the proof to the frame $(N_k,0,0)$ and the family $V_k$ we conclude that:
\[ ||V_k(-t_k)-\Pi_{t_k,0}T^e_{N_k} \phi||_{H^1(\mathbb{R}^2 \times \mathbb{T}^2)} \rightarrow 0
\]
\noindent as $k \rightarrow \infty$.\vspace{3mm}

\noindent The conclusion of the proof now follows from Theorem 4.5 and by noticing the behavior of $V_k$.\vspace{3mm}

\noindent \textbf{Large Scale Profiles.} Also, we need to analyze the large-scale profiles that appear in the profile decomposition in the next section. We need some notation: given $\psi \in H^{0,1}(\mathbb{R}^2 \times \mathbb{T}^2)$ and $M \leq 1$, we define the large-scale rescaling 
\begin{equation}\label{euqation} T^{ls}_{M} \psi(x,y)=M\tilde{\psi^{*}}(Mx,y) \quad \textmd{where} \quad  \tilde{\psi^{*}}(x,y)= P^{x}_{\leq M^{-1/100}} \psi(x,y).
\end{equation}
\noindent It is not hard to see that,\vspace{3mm}

\noindent $T^{ls}_{M} :H^{0,1}(\mathbb{R}^2 \times \mathbb{T}^2) \rightarrow H^1(\mathbb{R}^2 \times \mathbb{T}^2)$ is a linear bounded operator.\vspace{3mm}

\noindent It is crucial to study the behavior of nonlinear solutions $u_k$ of (1.1) with initial data as above. And these solutions are closely related to the solutions of the following cubic resonant systems.\vspace{3mm}

\noindent \textbf{The cubic resonant system.} We consider the cubic resonant system:
\[ (i\partial_t+\Delta_{x})u_j=\sum_{(j_1,j_2,j_3)\in R(j)} u_{j_{1}}\bar{u}_{j_{2}}u_{j_{3}},
\]
\[ R(j)=\{(j_1,j_2,j_3) \in (\mathbb{Z}^2)^3:j_1-j_2+j_3=j \quad and \quad  |j_1|^2-|j_2|^2+|j_3|^2=|j|^2 \}
\]
\noindent with initial data $\vec{u}(0)=\{u_j(0)\} \in h^1L^2$, which is defined in (2.2). The following energy 
\begin{equation}\label{euqation} E_{ls}(\vec{u})= \sum_{p \in \mathbb{Z}^2} (1+|p|^2) ||u_{p}||^2_{L^2_{x}(\mathbb{R}^2)}
\end{equation}
is conserved and so is the $h^1L^2$ norm of any solution of (1.6). 
\begin{conjecture} In addition to Conjecture 1.1, any initial data $\vec{u}_0$ of finite $E_{ls}$ energy leads to a global solution of (1.6) satisfying
\[||\vec{u}||_{\vec{W}}^2:=\sum_{p\in \mathbb{Z}^2}[1+|p|^2]||u_p||^2_{L^4_{x,t}(\mathbb{R}^2_x\times \mathbb{T}_t)}\leq S(E_{ls}(\vec{u}))
\]
\noindent where $S$ is some nondecreasing finite function. Also, this solution scatters in the sense that there exists $\vec{v}^{\pm \infty} \in h^1_pL^2$ such that
\begin{equation}\label{euqation}  \sum_{p \in \mathbb{Z}^2}[1+|p|^2]||u_p(t)-e^{it\Delta_x}v_p^{\pm \infty}||_{L^2_{x}(\mathbb{R}^2)} \rightarrow 0 \quad as \quad t \rightarrow \pm \infty.
\end{equation}
\end{conjecture}
\noindent As we show in the Section 8, this conjecture is also implied by the Theorem 1.2. In addition, by using the local well-posedness theory for (1.6), this conjecture is true under the smallness hypothesis $E_{ls}(\vec{u})<\delta$ for some $\delta>0$. Finally the result of [9] implies that the conjecture is true if one adds the additional condition that $\vec{u}(0)$ is a scalar.\vspace{3mm}

\noindent By using the conjecture above and the persistence of regularity part of Theorem 8.1, we have:
\begin{theorem}\label{theorem}
Assume that Conjecture 5.1 hold true. Suppose that $\vec{u}_0 \in h^1L^2$ and that $\vec{u} \in C(\mathbb{R}:h^1L^2)$ is the solution of (1.6) with initial data $\vec{u}_0$ given by Conjecture 5.1. Suppose also that $\vec{v}_0 \in h^3H^2$ satisfies
\[ ||\vec{u}_0-\vec{v}_0||_{h^1L^2} \lesssim \epsilon ,
\]
\noindent and that $\vec{v}(t)$ is the solution to (1.6) with initial data $\vec{v}(0)=\vec{v}_0.$ Then, it holds that:
\[||(1-\Delta_x)\{(1+|p|^2)\vec{v}\}_p||_{L^{\infty}_t(h^1L^2) \cap \vec{W}(\mathbb{R})} \lesssim_{||u_0||_{h^1L^2}} ||v_0||_{h^3H^2},
\]
\[||\vec{u}-\vec{v}||_{L^{\infty}_t(h^1L^2) \cap \vec{W}(\mathbb{R})} \lesssim_{||u_0||_{h^1L^2}} \epsilon
\]
\noindent and there exists $\vec{\omega}^{\pm} \in h^3H^2$ such that
\[ \sum_{p \in \mathbb{Z}^2}[1+|p|^2]||v_p(t)-e^{it\Delta_x}\vec{\omega}_p^{\pm}||^2_{L^2_x(\mathbb{R}^2)} \rightarrow 0 \quad as \quad t\rightarrow \pm \infty.
\]\end{theorem}
\begin{lemma}Assume that Conjecture 5.1 holds true. Let $\psi \in H^{0,1}(\mathbb{R}^2\times \mathbb{T}^2)$, $T_0 \in (0,\infty)$, and $\rho \in \{0,1\}$ be given, and define $f_M=T_M^{ls}\psi(x,y)$. The following conclusions hold:\vspace{3mm}

\noindent (1) There is $M_0=M_0(\phi,T_0)$ sufficiently small such that for all $M \leq M_0$, there is a unique solution $U_M \in C((-T_0M^2,T_0M^2);H^1(\mathbb{R}^2 \times \mathbb{T}^2))$ of the initial-value problem
\begin{equation}\label{euqation} (i\partial_t+\Delta_{\mathbb{R}^2\times \mathbb{T}^2})U_M=\rho|U_M|^2U_M, \quad U_M(0)=f_M.
\end{equation}
\noindent Moreover, for any $M \leq M_0$,
\[||U_M||_{X^1(-T_0M^{-2},T_0M^{-2})} \lesssim_{E_{ls}(\psi)} 1.
\]
\noindent (2) Assume $\epsilon_1 \in (0,1]$ is sufficiently small (depending only on $E_{ls}(\psi)$), $v_0\in h^3H^2$, and $||\vec{\psi}-\vec{v}_0||_{h^1L^2} \leq \epsilon_1$. Let $\vec{v} \in C(\mathbb{R}:h^3H^2)$ denote the solution of the initial-value problem
\[(i\partial_t+\Delta_x)v_j=\rho \sum_{(p_1,p_2,p_3)\in R(j)} v_{p_{1}}\bar{v}_{p_{2}}v_{p_{3}},\quad v_j(0)=v_{0,j},j \in \mathbb{Z}^2.
\]
For $M \geq 1$ we define
\[v_{j,M}(x,t)=M v_j(Mx,M^2 t),\qquad (x,t)\in \mathbb{R}^2\times (-T_0M^{-2},T_0M^{-2}),
\]
\begin{equation}\label{euqation} V_M(x,y,t)=\sum_{q\in \mathbb{Z}^2} e^{-it|q|^2} e^{i\langle y,q \rangle} v_{q,M}(x,t), \qquad (x,y,t)\in \mathbb{R}^2 \times \mathbb{T}^2 \times (-T_0M^{-2},T_0M^{-2}),
\end{equation}
\noindent then
\begin{equation}\label{euqation} \limsup\limits_{M\rightarrow 0} ||U_M-V_M||_{X^1(-T_0M^{-2},T_0M^{-2})} \lesssim_{E_{ls}(\vec{\psi})} \epsilon_1. 
\end{equation}
\end{lemma}
\noindent \emph{Proof:} When $\rho=0$, it is trivial. It suffices to prove (2). \vspace{3mm}

\noindent First, by using Stricharz estimate on $\mathbb{R}^2$, we have that
\begin{equation}\label{euqation} \sum_{q\in \mathbb{Z}^2}\langle q \rangle^2 ||v_q||^2_{L_t^3L_x^{6}}\lesssim E_{ls}(\vec{\psi}),\quad \sum_{q\in \mathbb{Z}^2}\langle q \rangle^6 ||v_q||^2_{L_t^3W_x^{2,6}}\lesssim_{E_{ls}(\vec{\psi})} ||\vec{v}||^2_{h^3H^2}.
\end{equation}
\noindent We now want to show that $V_M$ is an approximate solution to (5.28) in the sense of Theorem 4.5.
\begin{align*}
(i\partial_t+\Delta_{\mathbb{R}^2 \times \mathbb{T}^2})V_M-|V_M|^2V_M=-\sum_{q\in \mathbb{Z}^2} e^{-it|q|^2}e^{i\langle y,q \rangle} \sum_{\vec{p}\in NR(q)} v_{p_1,M} \overline{v_{p_2,M}} v_{p_3,M}=RHS
\end{align*}
\[\Phi_{q,\vec{p}}=|p_1|^2-|p_2|^2+|p_3|^2-|q|^2
\]
\noindent where
\[ \mathcal{NR}(q)=\{\vec{p}=(p_1,p_2,p_3):p_1-p_2+p_3-q=0;\Phi_{q,\vec{p}} \neq 0\}.
\]
\noindent Naturally, $\mathcal{NR}$ is short for non-resonant. In addition, we now claim that
\begin{equation}\label{euqation}
||RHS||_{N^1(I)}\lesssim_{||\vec{v}_0||_{h^3H^2}} M.
\end{equation}

\noindent As in [13], We will estimate the high frequency part and low frequency part separately.
\[ RHS=P^x_{>2^{-10}} RHS+P^x_{<2^{-10}} RHS=P_{high} RHS+P_{low} RHS.
\]\noindent For the high frequency part,
\begin{align*}
||P^x_{>2^{-10}} RHS||^2_{N^1(0,S)} &\lesssim ||P^x_{>2^{-10}} \partial_x RHS||^2_{N^1(0,S)} =\sum_{q\in \mathbb{Z}^2} ||e^{i\langle y,q \rangle} P^x_{>2^{-10}} \mathcal{F}_y(RHS)(q)||^2_{N^1(0,S)}\\
&\lesssim \sum_{q\in \mathbb{Z}^2} \langle q \rangle ^{-2} [\langle q \rangle ^{2} \sum_{\vec{p} \in NR(q)}||\partial_x\{v_{p_1,M} \overline{v_{p_2,M}} v_{p_3,M} \}||_{L^1_t H^1_x}]^2 .
\end{align*}
\noindent Since $q \leq max\{p_1,p_2,p_3\}$, we see that, for any $q$,
\begin{align*}
&\langle q \rangle ^{2} \sum_{\vec{p} \in  \mathcal{NR}(q)} \partial_x\{v_{p_1,M} \overline{v_{p_2,M}} v_{p_3,M} \}||_{L^1_t H^1_x}\\
&\lesssim \sum_{\vec{p} \in  \mathcal{NR}(q)} \langle p_1 \rangle^2 ||\partial_x(v_{p_1,M})||_{L_t^3W_x^{1,6}} \Pi_{j=2}^{3}\{\langle p_j \rangle^2 ||v_{p_j,M}||_{L_t^3W_x^{1,6}} \} 	.
\end{align*}
\noindent Thus the high frequency part can be handled.\vspace{3mm}

\noindent For the low frequency part, we can use integration by parts,
\begin{align*}
&\int_0^S e^{i(S-\sigma)\Delta_{\mathbb{R}^2\times \mathbb{T}^2}} P_{low}RHS(\sigma)d\sigma \\
&=-\sum_{q\in \mathbb{Z}^2}\sum_{\vec{p} \in  \mathcal{NR}(q)}e^{-iS[|q|^2+\Phi_{q,\vec{p}}]} e^{i\langle y,q\rangle} \int_0^S e^{i(S-\sigma)[\Delta_x+\Phi_{q,\vec{p}}]} P_{low}(v_{p_1,M}\overline{v_{p_2,M}}v_{p_3,M})d\sigma\\
&=\sum_{q\in \mathbb{Z}^2}\sum_{\vec{p} \in  \mathcal{NR}(q)} e^{-iS[|q|^2+\Phi_{q,\vec{p}}]} e^{i\langle y,q\rangle} \times\\
&\{[ie^{i(S-\sigma)[\Delta_x+\Phi_{q,\vec{p}}]}(\Delta_x+\Phi_{q,\vec{p}})^{-1}P_{low}(v_{p_1,M}\overline{v_{p_2,M}}v_{p_3,M})]^S_0\\
&-i\int_0^S e^{i(S-\sigma)[\Delta_x+\Phi_{q,\vec{p}}]} (\Delta_x+\Phi_{q,\vec{p}})^{-1}P_{low} \partial_{\sigma}\{v_{p_1,M}\overline{v_{p_2,M}}v_{p_3,M}\}d\sigma\}.
\end{align*}
\noindent Now we need to estimate three terms: `$S$-boundary', `$0$-boundary', and the other term. We will estimate them separately. \vspace{3mm}

\noindent One important thing we should be clear is that $(\Delta_{x}+\Phi_{q,\vec{p}})^{-1}P_{low}$ is bounded on $L^2_x(\mathbb{R}^2)$ noticing the non-resonant condition and the low frequency cutoff as in [7, 13]. \vspace{3mm}

\noindent For the `$S$-boundary':
\begin{align*}
&||\sum_{q\in \mathbb{Z}^2}\sum_{\vec{p} \in  \mathcal{NR}(q)} e^{-iS[|q|^2+\Phi_{q,\vec{p}}]} e^{i\langle y,q\rangle}(\Delta_x+\Phi_{q,\vec{p}})^{-1} P_{low}(v_{p_1,M}\overline{v_{p_2,M}}v_{p_3,M}) ||^2_{X^1(I)} \\
&\lesssim \sum_{q\in \mathbb{Z}^2} \langle q \rangle^2 ||\sum_{\vec{p} \in  \mathcal{NR}(q)} (v_{p_1,M}\overline{v_{p_2,M}}v_{p_3,M})(0)||^2_{H^1_x(\mathbb{R}^2)} \\
&+\sum_{q\in \mathbb{Z}^2} \langle q \rangle^2 ||(i\partial_t+\Delta_x)\sum_{\vec{p} \in  \mathcal{NR}(q)} v_{p_1,M}\overline{v_{p_2,M}}v_{p_3,M}||^2_{L^1_t H^1_x}\\
&\lesssim_{||\vec{v}_0||_{h^3H^2}} M^4.
\end{align*}
\noindent For the  `$0$-boundary':
\begin{align*}
&||\sum_{q\in \mathbb{Z}^2}\sum_{\vec{p} \in  \mathcal{NR}(q)} e^{-iS\Delta_{\mathbb{R}\times \mathbb{T}^2}} e^{i\langle y,q\rangle}  (\Delta_x+\Phi_{q,\vec{p}})^{-1} P_{low}[v_{p_1,M}\overline{v_{p_2,M}}v_{p_3,M}](0)||^2_{X^1(I)} \\
&\lesssim \sum_{q\in \mathbb{Z}^2} \langle q \rangle^2 ||\sum_{\vec{p} \in  \mathcal{NR}(q)} (v_{p_1,M}\overline{v_{p_2,M}}v_{p_3,M})(0)||^2_{H^1_x(\mathbb{R}^2)}\lesssim_{||\vec{v}_0||_{h^3H^2}} M^4.
\end{align*}
\noindent For the other term:
\begin{align*}
&||\sum_{q\in \mathbb{Z}^2}\sum_{\vec{p} \in  \mathcal{NR}(q)} e^{-iS[|q|^2+\Phi_{q,\vec{p}}]} e^{i\langle y,q\rangle} \int_0^S e^{i(S-\sigma)[\Delta_x+\Phi_{q,\vec{p}}]} (\Delta_x+\Phi_{q,\vec{p}})^{-1}P_{low} \partial_{\sigma}\{v_{p_1,M}\overline{v_{p_2,M}}v_{p_3,M}\}d\sigma||^2_{X^1(I)} \\
&\lesssim \sum_{q\in \mathbb{Z}^2} \langle q \rangle^2 ||\sum_{\vec{p} \in  \mathcal{NR}(q)} \partial_{\sigma}\{v_{p_1,M}\overline{v_{p_2,M}}v_{p_3,M}\}||^2_{L^1_t H^1_x}\lesssim_{||\vec{v}_0||_{h^3H^2}} M^4.
\end{align*}
\noindent This finishes the proof of (5.23).\vspace{3mm}

\noindent We also have that
\[||V_M||^2_{L_t^{\infty}H^1_{x,y}(\mathbb{R}^2\times \mathbb{T}^2\times I)}\leq \sum_{q\in\mathbb{R}^2}\langle q \rangle^2||v_{q,M}||^2_{L_t^{\infty} H^1_{x}} \leq C||\vec{u}(0)||^2_{h^1L^2}+C(M)||\vec{v}||^2_{h^1H^1}
\]
\noindent and that
\begin{equation}\label{euqation} ||V_M||_{X^1(I)} \lesssim_{||\vec{v}(0)||_{h^1L^2}} 1+C(||\vec{v}||_{h^3H^2})M.
\end{equation}

\noindent Moreover, we have (Using Lemma 8.2)
\begin{align*}
||(i\partial_t+\Delta_{\mathbb{R}^2\times \mathbb{T}^2})V_M||^2_{N^1(I)}&\lesssim \sum_{q\in \mathbb{Z}^2}\langle q \rangle^2 ||\sum_{\vec{p}\in  \mathcal{R}(q)} v_{p_1,M}\overline{v_{p_2,M}}v_{p_3,M}||^2_{L^1_t H^1_x} \\
&\lesssim \sum_{q\in \mathbb{Z}^2}[\sum_{\vec{p}\in  \mathcal{R}(q)}\Pi_{k=1}^3 \langle p_k \rangle ||v_{p_k,M}||_{L_t^3W_x^{1,6}} \times \langle q \rangle \Pi_{k=1}^3 \langle p_k \rangle^{-1}]^2 \\
&\lesssim [\sum_{p\in \mathbb{Z}^2}\langle p \rangle^2 ||v_{p,M}||^2_{L_t^3W_x^{1,6}}]^3.
\end{align*}
\noindent which justifies (5.33).\vspace{3mm}

\noindent By using Theorem 4.5. we conclude that, for $M$ small enough (depending on $\vec{v}_0$), the solution $U_M$ of (1.1) with initial data $V_M(0)$ exists on $I$ and that
\[||U_M-V_M||_{} \lesssim \epsilon_1 +C(||\vec{v}_0||_{h^3H^2})M,
\]

\noindent which ends the proof.
\begin{lemma}For any $\psi \in H^{0,1}(\mathbb{R}^2 \times \mathbb{T}^2)$ and any $\epsilon > 0$, there exists $T_0=T(\psi,\epsilon)$ and $M_0=M(\psi,\epsilon)$such that for any $T\geq T_0$ and any $M \leq M_0$,
\[
||e^{it\Delta_{\mathbb{R}^2 \times \mathbb{T}^2}}T^{ls}_{M} \psi||_{Z(\{M^2|t| \geq T_0\})}  \lesssim \epsilon.
\]\end{lemma}
\noindent \emph{Proof:} By Stricharz estimate on $\mathbb{R}^2$, and dominated convergence theorem, there exists $T_0$ such that 
\begin{equation}\label{euqation}
\sum_{p\in \mathbb{Z}^2} \langle p \rangle^2 ||e^{it\Delta_x}\psi_p||^2_{L^4_{x,t}(\mathbb{R}^2\times \{|t| \geq T_0\})} \leq \epsilon^{1000}.
\end{equation}
\noindent Let $I= \{|t| \geq T_0 \}  $ and $I_M=\{M^2|t| \geq T_0 \}$. We have that 
\begin{align*}
e^{it\Delta_{\mathbb{R}^2\times \mathbb{T}^2}} T_M^{ls} \psi &=\sum_{q\in \mathbb{Z}^2}e^{i(\langle q,y \rangle-|q|^2t)} (Me^{iM^2t\Delta_x}{\tilde{\psi}_q}^{M}(Mx))\\
&=\sum_{q\in \mathbb{Z}^2} e^{i(\langle q,y \rangle-|q|^2t)} v_{q,M}(t,x)
\end{align*}
\noindent where we denoted by:
\[ v_{p,M}(x,t)=Me^{iM^2t\Delta_x} {\tilde{\psi}_p}^{M}(Mx).
\]\noindent Noticing that $e^{i\langle q,y \rangle} v_{q,M}(x,t)$ is supported in Fourier space in the box centered at $q$ of radius 2 and Bernstein's inequality in $y$, we can estimate:
\[||P_N e^{it\Delta_{\mathbb{R}^2\times \mathbb{T}^2}} T_M^{ls} \psi||_{L^4_{x,y,t}} \lesssim N^{\frac{1}{2}}||(\sum_{|q| \sim N}|v_{q,M}(x,t)|^2)^{\frac{1}{2}}||_{L^4_{x,t}} \lesssim N^{\frac{1}{2}} (\sum_{|q| \sim N}||v_{q,M}||^2_{L^4_{x,t}})^{\frac{1}{2}}.
\]
\noindent Thus we know:
\begin{equation} 
\aligned
\sum_{N \geq 1} N^2||P_N e^{i\Delta_{\mathbb{R}^2\times \mathbb{T}^2}} T_M^{ls} \psi||_{L^4_{x,y,t}(\mathbb{R}^2\times \mathbb{T}^2 \times I_M)}^4 &\lesssim \sum_{N \geq 1}N^4 (\sum_{|q| \sim N}||v_{q,M}||^2_{L^4_{x,t}(\mathbb{R}^2\times I_M)})^2\\
&\lesssim (\sum_{N \geq 1}N^2 \sum_{|q| \sim N}||v_{q,M}||^2_{L^4_{x,t}(\mathbb{R}^2\times I_M)})^2 \\
&\lesssim (\sum_{q\in \mathbb{Z}^2}  \langle q \rangle^2 ||e^{it\Delta_x}{\tilde{\psi}_q}^{M}||^2_{L^4_{x,t}(\mathbb{R}^2\times I)})^2 \lesssim \epsilon^{2000}. 
\endaligned
\end{equation}

\noindent According to the definition of $Z$-norm, we finish the proof of this lemma.\vspace{3mm}

\noindent Now we can describe the nonlinear solutions of the Initial Value Problem (1.1) corresponding to large-scale profile. In view of the profile analysis in the next section, we need to consider the renormalized large-scale frames by:\vspace{3mm}

\noindent $\tilde{\mathcal{F}}_{ls} := \{(M_k,t_k,p_k,\xi_k)_{k}:M_k \leq 1,M_k \rightarrow 0,p_k=(x_0,0) \in \mathbb{R}^2 \times \mathbb{T}^2 $, and $\xi_k \in \mathbb{R}$ ,with $\xi_k \rightarrow \xi_{\infty} \in \mathbb{R}$ and either $t_k=0$ or $M_k^{2}t_k \rightarrow \pm \infty$ and either $\xi_k=0$ or  $ M_k^{-1} \xi_k \rightarrow \pm \infty\}$.
\begin{theorem}\label{theorem} Assume Conjecture 5.1 holds true. Fix a renormalized large-scale frame $(M_k,t_k,(x_k,0),\xi_k)_k \in \tilde{\mathcal{F}}_{ls}$ and $\psi \in H^{0,1}(\mathbb{R}^2 \times \mathbb{T}^2)$ let
\[U_k(0)=\Pi_{t_k,x_k} e^{i\xi_k x}T^{ls}_{M_k} \psi.
\]
\noindent (1) For $k$ large enough (depending on $\psi$, $S$), there is a nonlinear solution $U_k \in X^1_c(\mathbb{R})$ of the equation (1.1) satisfying:
\begin{equation}\label{euqation}
||U_k||_{X^1(\mathbb{R})} \lesssim_{E_{ls}(\psi)} 1.
\end{equation}
\noindent (2) There exists a solution $\vec{v} \in C(\mathbb{R}:h^1L^2)$ of (1.6) with scattering data $\vec{v}_0^{\pm \infty}$ such that the following holds, up to a subsequence: for $\epsilon >0$, there exists $T(\psi,\epsilon)$ such that for all $T\geq T(\psi,\epsilon)$, there holds that 
\begin{equation}\label{euqation}
||U_k-W_k||_{X^1(\{|t-t_k| \leq TM_k^{-2}\})} \leq \epsilon ,
\end{equation}
\noindent for $k$ large enough, where 
\[W_k(x,t)=e^{-i\eta|\xi_k|^2} e^{ix\xi_k}V_{M_k}(x-x_k-2\xi_k \eta,y,\eta),\quad \eta=t-t_k
\]
\noindent with $V_k$ defined as before. Moreover,
\begin{equation}\label{euqation}
||U_k(t)-\Pi_{t_k-t,x_k}e^{ix\xi_k}T^{ls}_{M_k}\underline{\vec{v}_0^{\pm \infty}}||_{X^1(\{\pm(t-t_k) \geq TM_k^{-2}\})} \leq \epsilon,
\end{equation}
\noindent for $k$ large enough (depending on $\psi,\epsilon,T$). 
\end{theorem}
\noindent \emph{Proof:} Without loss of generality, we may assume that $x_k=0$. Using a Galilean transform and the fact that $\xi_k$ is bounded, we may assume that $\xi_k=0$ for all $k$.\vspace{3mm}

\noindent First we can consider the case when $t_k=0$ for all $k$ and we let $\vec{v}$ be the solution of (1.6) with initial data $\vec{\psi}$. Then by using Theorem 5.5, we see that there exists $T_0=T_0(\psi,\epsilon)$ such that
\begin{equation}\label{euqation}
\sup\limits_{t \geq T_0} ||\vec{v}(t)-e^{it\Delta_x}\vec{v}_0^{+\infty}||_{h^1L^2}+||e^{it\Delta_x}\vec{v}_0^{+\infty}||_{\vec{W}(\{t\geq T_0)\}} \leq \epsilon,
\end{equation}
\[\sup\limits_{t \leq -T_0} ||\vec{v}(t)-e^{it\Delta_x}\vec{v}_0^{-\infty}||_{h^1L^2}+||e^{it\Delta_x}\vec{v}_0^{-\infty}||_{\vec{W}(\{t\leq T_0)\}} \leq \epsilon,
\]
fix $T\geq T_0$. Applying Lemma 5.6, we see that, if $k$ is large enough,
\[||U_k-V_{M_k}||_{X^1(\{ |t|\leq T M_k^{-2}\})} \leq \epsilon.
\]
By using Stricharz estimates, (5.39) and Lemma 5.8, we see 
\[||e^{it\Delta_{\mathbb{R}^2\times \mathbb{T}^2}}U_k(\pm TM_k^{-2})||_{Z(\pm t \geq M_k^{-2}T)} \leq \epsilon.
\]
\noindent Now, Theorem 4.3 implies $U_k$ extends to a global solution $U_k \in X^1_c(\mathbb{R})$ satisfying (5.38).\vspace{3mm}

\noindent Now we consider the other case when $M_k^{2}|t_k| \rightarrow \infty$. For definiteness, we assume that  $M_k^{2} t_k \rightarrow +\infty$ and let $\vec{v}$ be the solutions to (1.6) satisfying
\[ ||\vec{v}(t)-e^{it\Delta_x}\vec{\psi}||_{h^1L^2} \rightarrow 0
\]
\noindent as $t \rightarrow -\infty$. Let $\psi^{'}=\underline{\vec{v}(0)} \in H^{0,1}(\mathbb{R}^2\times \mathbb{T}^2)$ and let $V_k$ be the solution of (1.1) with initial data $T^{ls}_{M_k}\psi^{'}$. Applying the first case of the proof to the frame $(N_k,0,0,0)$ and the family $V_k$ we conclude that:
\[ ||V_k(-t_k)-\Pi_{t_k,0}T^{ls}_{M_k}\psi||_{H^1(\mathbb{R}^2 \times \mathbb{T}^2)} \rightarrow 0
\]
\noindent as $k \rightarrow \infty$. The conclusion of the proof now follows from Theorem 4.5 and the behavior of $V_k$.
\section{Profile Decomposition}
\noindent Then we can define three different kind of profiles corresponding to different frames. We use frames to make the different profiles written in the same form as in [13, 17, 18, 19].
\begin{definition} [Frames and Profiles](1) We define a frame to be sequence $(N_k,t_k,p_k,\xi_{k})_{k} \in 2^{\mathbb{Z}} \times \mathbb{R} \times (\mathbb{R}^2 \times \mathbb{T}^2) \times \mathbb{R} $ which contains 4 elements. And we can distinguish three types of profiles as follows.

\noindent a) A Euclidean frame is a sequence $\mathcal{F}_e =(N_k,t_k,p_k,0)$ with $N_k \geq 1, N_k \rightarrow \infty, t_k \in \mathbb{R} , p_k \in \mathbb{R}^2 \times \mathbb{T}^2$.

\noindent b) A Large-scale frame is a sequence $\mathcal{F}_{ls} =(M_k,t_k,p_k,\xi_k)$ with $M_k \leq 1,M_k \rightarrow 0 ,t_k \in \mathbb{R}, p_k \in \mathbb{R}^2\times \mathbb{T}^2$.

\noindent c) A  Scale-one frame is a sequence $\mathcal{F}_1 =(1,t_k,p_k,0)$ with $t_k \in \mathbb{R}, p_k \in \mathbb{R}^2 \times \mathbb{T}^2$.

\noindent (2) We say that two frames $(N_k,t_k,p_k,\xi_{k})_{k}$ and $(M_k,s_k,q_k,\eta_{k})_{k}$ are orthogonal if
\[ \lim_{k \rightarrow +\infty} (|ln \frac{N_k}{M_k}|+N_k^2|t_k-s_k|+N_k^{-1}|\xi_k-\eta_k|+N_k|(p_k-q_k)-2(t_k-s_k)\xi_k|)=+\infty.
\]
\noindent(3) We associate a profile defined as: 

\noindent a) If $\mathcal{O}=(N_k,t_k,p_k,0)_k $ is a Euclidean frame and for $\phi \in \dot{H}^1(\mathbb{R}^4)$ we define the Euclidean profile associated to $(\phi,\mathcal{O})$ as the sequence $\tilde{\phi}_{\mathcal{O},k}$ with
\[\tilde{\phi}_{\mathcal{O},k}=\Pi_{t_k,p_k}(T^{e}_{N_k}\phi)(x,y).
\]

\noindent b) If $\mathcal{O}=(M_k,t_k,p_k,\xi_k)_k$ is a large scale frame, if $p_k=(x_k,0)$ and if $\psi \in H^{0,1}(\mathbb{R}^2 \times \mathbb{T}^2)$, we define the large scale profile associated to $(\psi,\mathcal{O})$ as the sequence $\tilde{\psi}_{\mathcal{O},k}$ with
\[\tilde{\psi}_{\mathcal{O},k}=\Pi_{t_k,p_k}[e^{i\xi_k x}T^{ls}_{M_k}\psi(x,y)].
\]

\noindent c) If $\mathcal{O}=(1,t_k,p_k,0)_k$  is a scale one frame, if $W \in H^1(\mathbb{R}^2 \times \mathbb{T}^2)$, we define the scale one profile associated to $(W,\mathcal{O})$ as $\tilde{W}_{\mathcal{O},k}$ with
\[\tilde{W}_{\mathcal{O},k}=\Pi_{t_k,p_k} W.
\]

\noindent (4) Finally, we say that a sequence of functions $\{f_k\}_k \subset H^1(\mathbb{R}^2 \times \mathbb{T}^2)$ is absent from a frame $\mathcal{O}$ if, up to a subsequence: 

\noindent $\langle f_k,\tilde{\psi}_{\mathcal{O},k} \rangle_{H^1 \times H^1}\rightarrow 0$ as $k \rightarrow \infty$ for any profile $\tilde{\psi}_{\mathcal{O},k}$ associated with $\mathcal{O}$.
\end{definition}
\begin{lemma}[Frame equivalence and orthogonality]
\noindent (1) Suppose that $\mathcal{O}$ and $\mathcal{O}^{'}$ are equivalent Euclidean (respectively large-scale or scale-one) frames, then there exists an isometry $L$ of $\dot{H}^1(\mathbb{R}^4)$ (resp. $H^{0,1}(\mathbb{R}^2 \times \mathbb{T}^2)$, $H^1(\mathbb{R}^2 \times \mathbb{T}^2)$) such that, for any profile generator $\psi \in $ $\dot{H}^1(\mathbb{R}^4)$ (resp. $H^{0,1}(\mathbb{R}^2 \times \mathbb{T}^2)$, $H^1(\mathbb{R}^2 \times \mathbb{T}^2)$), it holds that, up to a subsequence:
\begin{equation}\label{equation}
\limsup\limits_{k \rightarrow +\infty} ||\widetilde{L\psi}_{\mathcal{O},k}-\tilde{\psi}_{\mathcal{O}^{'},k}  ||_{H^1(\mathbb{R}^2 \times \mathbb{T}^2)}=0 .
\end{equation}
\noindent (2) Suppose that $\mathcal{O}$ and $\mathcal{O}^{'}$ are orthogonal frames and $\tilde{\psi}_{\mathcal{O},k}$ and $\tilde{\phi}_{\mathcal{O}^{'},k}$ are two profiles associated with $\mathcal{O}$ and $\mathcal{O}^{'}$ respectively. Then
\[\lim\limits_{k \rightarrow +\infty} \langle \tilde{\psi}_{\mathcal{O},k},\tilde{\phi}_{\mathcal{O}^{'},k} \rangle_{H^1\times H^1(\mathbb{R}^2\times \mathbb{T}^2)} =0,
\]
\[\lim\limits_{k \rightarrow +\infty} \langle |\tilde{\psi}_{\mathcal{O},k}|^2,|\tilde{\phi}_{\mathcal{O}^{'},k}|^2 \rangle=0.
\]
\noindent (3) If $\mathcal{O}$ is a Euclidean frame and $\tilde{\psi}_{\mathcal{O},k}$, and $\tilde{\phi}_{\mathcal{O^{'}},k}$are two profiles associated to $\mathcal{O}$, then:
\[ \lim\limits_{k \rightarrow \infty} \langle \tilde{\psi}_{\mathcal{O},k},\tilde{\phi}_{\mathcal{O},k} \rangle_{H^1\times H^1(\mathbb{R}^2 \times \mathbb{T}^2)}=\langle \phi,\psi \rangle_{\dot{H}^{1}\times \dot{H}^{1}(\mathbb{R}^4)},
\]
\[ \lim\limits_{k \rightarrow \infty} ||\tilde{\psi}_{\mathcal{O},k}||_{L^2}+||\tilde{\phi}_{\mathcal{O},k}||_{L^2}=0.
\]
\noindent (4) If $\mathcal{O}$ is a scale-one frame and $\tilde{\psi}_{\mathcal{O},k}$, $\tilde{\phi}_{\mathcal{O},k}$ are two profiles associated to $\mathcal{O}$, then:
\[ \lim\limits_{k \rightarrow \infty} \langle \tilde{\psi}_{\mathcal{O},k},\tilde{\phi}_{\mathcal{O},k}\rangle_{H^1\times H^1(\mathbb{R}^2 \times \mathbb{T}^2)}=\langle \phi,\psi \rangle_{{H^1}\times {H^1}(\mathbb{R}^2\times \mathbb{T}^2)}.
\]
\noindent (5) If $\mathcal{O}$ is a large-scale frame and $\tilde{\psi}_{\mathcal{O},k}$, $\tilde{\phi}_{\mathcal{O},k}$ are two profiles associated to $\mathcal{O}$, then:
\[ \lim\limits_{k \rightarrow +\infty} ||\tilde{\psi}_{\mathcal{O},k}||_{L_{x,y}^4(\mathbb{R}^2\times \mathbb{T}^2)}=0,
\]
\begin{align*} 
\lim\limits_{k \rightarrow +\infty} \langle \tilde{\psi}_{\mathcal{O},k},\tilde{\phi}_{\mathcal{O},k} \rangle_{H^1 \times H^1(\mathbb{R}^2 \times \mathbb{T}^2))}&=\langle \psi,\phi \rangle_{H^{0,1} \times H^{0,1}(\mathbb{R}^2 \times \mathbb{T}^2)}+|\xi_{\infty}|^2 \langle \psi,\phi \rangle _{L^2\times L^2(\mathbb{R}^2 \times \mathbb{T}^2)}\\ &\simeq \langle \psi,\phi \rangle_{H^{0,1} \times H^{0,1}(\mathbb{R}^2 \times \mathbb{T}^2)}.
\end{align*}
\end{lemma}
\noindent The proof is straightforward. (see [17, 18, 19])\vspace{3mm}

\noindent The following theorem is the main theorem of this section, i.e. profile decomposition theorem, which will be used in the energy induction method.
\begin{theorem}[Profile Decomposition] Assume $\{ \phi_k\}_{k}$ is a sequence of functions satisfying
\begin{equation}\label{equation}   
\sup\limits_{k \geq 0} ||\phi_k||_{L^2(\mathbb{R}^2 \times \mathbb{T}^2)}+||\nabla_{x,y} \phi_k||_{L^2(\mathbb{R}^2 \times \mathbb{T}^2)} < E \leq \infty .
\end{equation}
\noindent Then there exists a subsequence (for convenience which we also denote by $\phi_k$)
, a family of Euclidean profiles $\tilde{\varphi}_{\mathcal{O}^\alpha,k}$, a family of large scale profiles $\tilde{\psi}_{\mathcal{S}^\beta,k}$ and a family of scale-one profiles $\tilde{W}_{\mathcal{O}^\gamma,k}$ such that, for any $A \geq 1$ and any $k \geq 0$ in the subsequence 
\begin{equation}\label{equation}
\phi_k(x,y)=\sum_{1 \leq \alpha \leq A}\tilde{\varphi}^{\alpha}_{\mathcal{O}^\alpha,k} +\sum_{1 \leq \beta \leq A}\tilde{\psi}^{\beta}_{\mathcal{S}^\beta,k} +\sum_{1 \leq \gamma \leq A}\tilde{W}^{\gamma}_{\mathcal{O}^\gamma,k} +R_k^A(x,y),
\end{equation}
\noindent with 
\begin{equation}\label{equation}
\lim\limits_{A \rightarrow \infty } \limsup\limits_{k \rightarrow \infty} ||e^{it \Delta_{\mathbb{R}^2 \times  \mathbb{T}^2}}R_k^A||_{Z(\mathbb{R})}=0. 
\end{equation}
\noindent In addition, all the frames are pairwise orthogonal and we have the following orthogonality property:
\[M(\phi_{k})=\sum_{1 \leq \beta \leq A} M(\psi^{\beta}) +\sum_{1 \leq \gamma \leq A} M(W^{\gamma}) +M(R_k^A)+o_{A,k \rightarrow +\infty}(1),
\]
\begin{equation}\label{equation}
\aligned
||\nabla_{x,y} \phi_k||^2_{L^2(\mathbb{R}^2\times \mathbb{T}^2)} &=\sum_{1 \leq \alpha \leq A}||\varphi^{\alpha}||_{\dot{H}_1(\mathbb{R}^4)}+\sum_{1 \leq \beta \leq A}[|\xi_{\infty}^{\beta}|^2 M(\psi^{\beta})+||\nabla_y \psi^{\beta}||^2_{L^2}]\\ 
&+\sum_{1 \leq \gamma \leq A}||\nabla_{x,y}W^{\gamma}||^2_{L^2(\mathbb{R}^2\times \mathbb{T}^2)}+||\nabla_{x,y}R^A_k||^2_{L^2(\mathbb{R}^2\times \mathbb{T}^2)}+o_{A,k \rightarrow +\infty}(1),
\endaligned
\end{equation}
\[ ||\phi_{k}||_{L^4(\mathbb{R}^2\times \mathbb{T}^2)}^{4}=\sum_{1 \leq \alpha \leq A} ||\varphi^{\alpha}||_{L^4}^4 +\sum_{1 \leq \gamma \leq A} ||W^{\gamma}||_{L^4}^{4} +o_{A \rightarrow +\infty,k \rightarrow +\infty}(1),
\]\noindent where $\xi_{\infty}^{\beta}=\lim_{k \rightarrow +\infty}\xi_k^{\beta}$, $o_{A,k \rightarrow +\infty}(1) \rightarrow 0$ as $k \rightarrow +\infty$ for each fixed $A$, and $o_{A\rightarrow +\infty,k\rightarrow +\infty}(1) \rightarrow 0$ in the ordered limit $\lim\limits_{A \rightarrow +\infty}\lim\limits_{k \rightarrow +\infty}$.
\end{theorem}
\noindent As in [13, 17, 18], this follows from iteration of the following statement,
\begin{lemma} Let $\delta > 0$. Assume that $\phi_k$ is a sequence satisfying:
\begin{equation}\label{equation}   \sup\limits_{k \geq 0} ||\phi_k||_{L^2(\mathbb{R}^2 \times \mathbb{T}^2)}+||\nabla_{x,y} \phi_k||_{L^2(\mathbb{R}^2 \times \mathbb{T}^2)} < E \leq \infty .
\end{equation}
\noindent then there exists a subsequence (for convenience, which we also denote by $\phi_k$), $A=A(E,\delta)$ Euclidean profiles $\tilde{\varphi}^{\alpha}_{\mathcal{O}^\alpha,k}$, $A$ large scale profiles $\tilde{\psi}^{\beta}_{\mathcal{S}^\beta,k} $, and $A$ scale 1 profiles $\tilde{W}^{\gamma}_{\mathcal{O}^\gamma,k} $ such that, for any $k \geq 0$ in the subsequence
\begin{equation}\label{equation}
\phi_k(x,y)=\sum_{1 \leq \alpha \leq A}\tilde{\varphi}^{\alpha}_{\mathcal{O}^\alpha,k} +\sum_{1 \leq \beta \leq A}\tilde{\psi}^{\beta}_{\mathcal{S}^\beta,k} +\sum_{1 \leq \gamma \leq A}\tilde{W}^{\gamma}_{\mathcal{O}^\gamma,k} +R_k^A(x,y) \end{equation}
\noindent with
\[  \limsup\limits_{k \rightarrow \infty} [||e^{it \Delta_{\mathbb{R}^2 \times  \mathbb{T}^2}}R_k^A||_{Z(\mathbb{R})}+\sup\limits_{t} ||e^{it \Delta_{\mathbb{R}^2 \times  \mathbb{T}^2}}R_k^A||_{L_{x,t}^4(\mathbb{R}^2 \times \mathbb{T}^2)}]  \leq \delta.
\]\noindent Also, the frames are pairwise orthogonal and we have the following orthogonality property:
\[M(\phi_{k})=\sum_{1 \leq \beta \leq A} M(\psi^{\beta}) +\sum_{1 \leq \gamma \leq A} M(W^{\gamma}) +M(R_k^A)+o_{k \rightarrow +\infty}(1),
\]\begin{align*}
||\nabla_{x,y} \phi_k||^2_{L^2(\mathbb{R}^2\times \mathbb{T}^2)} &=\sum_{1 \leq \alpha \leq A}||\varphi^{\alpha}||_{\dot{H}_1(\mathbb{R}^4)}+\sum_{1 \leq \beta \leq A}[|\xi_{\infty}^{\beta}|^2 M(\psi^{\beta})+||\nabla_y \psi^{\beta}||^2_{L^2}]\\
&+\sum_{1 \leq \gamma \leq A}||\nabla_{x,y}W^{\gamma}||^2_{L^2(\mathbb{R}^2\times \mathbb{T}^2)}+||\nabla_{x,y}R^A_k||^2_{L^2(\mathbb{R}^2\times \mathbb{T}^2)}+o_{k \rightarrow +\infty}(1),
\end{align*}
\begin{align*}
 ||\phi_{k}||_{L^4(\mathbb{R}^2\times \mathbb{T}^2)}^{4}&=\sum_{1 \leq \alpha \leq A} ||\varphi^{\alpha}||_{L^4(\mathbb{R}^4)}^4 +\sum_{1 \leq \gamma \leq A} ||W^{\gamma}||_{L^4(\mathbb{R}^2\times \mathbb{T}^2)}^{4} \\
 &+||R^A_k||^4_{L^4(\mathbb{R}^2\times \mathbb{T}^2)}+o_{k \rightarrow +\infty}(1),
\end{align*}\noindent where $o_{k \rightarrow +\infty}(1) \rightarrow 0$ as $k \rightarrow +\infty$.
\end{lemma}
\noindent The proof will be completed in two steps: first, we extract the Euclidean and scale-one profiles by studying the defects of compactness of the Stricharz estimate. This extraction leaves only sequences whose linear flow has small critical Besov norm but large $Z(\mathbb{R})$ norm, from which we extract the scale-one profiles and the large scale profiles and finish the proof. For a sequence of functions $\{f_k\}$ in $H^1(\mathbb{R}^2 \times \mathbb{T}^2)$, we define the following functional:
\begin{equation}\label{equation} 
\Lambda_{\infty}(\{f_k\})=\limsup\limits_{k \rightarrow \infty}||e^{it\Delta}f_k||_{L_t^{\infty}B_{\infty,\infty}^{-1}}=\limsup\limits_{k \rightarrow \infty} \sup\limits_{\{N,t,x,y\}} N^{-1} |(e^{it\Delta}P_N f_k)(x,y)| .
\end{equation}
\noindent where the supremum is taken over all scales $N\geq 1$, times $t\in \mathbb{R}$ and $(x,y)\in \mathbb{R}^2\times \mathbb{T}^2$.
\begin{lemma}Let $v > 0$. Assume that $\phi_k$ is a sequence satisfying (6.6) (bounded in $H^1(\mathbb{R}^2\times \mathbb{T}^2)$), then there exists a subsequence of $\phi_k$, $A$ Euclidean profiles $\tilde{\varphi^{\alpha}}_{\mathcal{O}^\alpha,k}$, and $A$ scale-one profiles $\tilde{W^{\gamma}}_{\mathcal{O}^\gamma,k}$ such that, for any $k \geq 0$ in the subsequence
\begin{equation}\label{equation} 
\phi_k^{'}(x,y)=\phi_k(x,y)-\sum_{1\leq \alpha \leq A} \tilde{\varphi}^{\alpha}_{\mathcal{O}^\alpha,k}-\sum_{1\leq \gamma \leq A}\tilde{W}^{\gamma}_{\mathcal{O}^\gamma,k} \end{equation} \noindent satisfies
\begin{equation}\label{equation} 
\Lambda_{\infty}(\{ \phi_k^{'}\}) < v.
\end{equation}
\noindent Besides, all the frames involved are pairwise orthogonal and $\phi_k^{'}$ is absent from all these frames.
\end{lemma}
\noindent \emph{Proof:} We first claim that if $\Lambda_{\infty}(\{f_k\}) \geq v$, then there exists a frame $\mathcal{O}$ and an associated profile $\tilde{\psi}_{\mathcal{O},k}$ satisfying
\begin{equation}\label{equation} 
\limsup\limits_{k \rightarrow \infty}  ||\tilde{\psi}_{\mathcal{O},k}||_{H^1} \lesssim 1 ,
\end{equation}
\noindent and
\begin{equation}\label{equation}
\limsup\limits_{k \rightarrow \infty} |\langle f_k,\tilde{\psi}_{\mathcal{O},k} \rangle_{ H^1 \times H^1 }| \gtrsim v.
\end{equation}
\noindent In addition, if $f_k$ was absent from a family of frames $\mathcal{O}^\alpha$, then $\mathcal{O}$ is orthogonal to all the previous frames.\vspace{3mm}

\noindent Let us prove the claim above first. By assumption, up to extracting a subsequence, there exists a sequence $\{N_k,t_k,(x_k,y_k)\}_k$ such that, for all $k$ 
\[ \frac{1}{2} v \leq N_k^{-1}|(e^{it_k\Delta}P_{N_k}f_k)(x_k,y_k)| \leq |\langle f_k, N_k^{-1}(e^{it_k\Delta}P_{N_k})\delta_{(x_k,y_k)} \rangle_{H^1 \times H^{-1}} |.
\]
\noindent Let us consider two situations.\vspace{3mm}

\noindent First, assume that $N_k$ remains bounded, then up to a subsequence, we can assume that $N_k \rightarrow N_{\infty}$ and since $N_k$ is dyadic, we may assume that $N_k=N_{\infty}$ for all $k$. In this case, we define the scale-one profile $\mathcal{O}=(1,t_k,(x_k,y_k),0)$ and
\[ \psi= (1-\Delta)^{-1} N_{\infty}^{-1} P_{N_{\infty}}\delta_{(0,0)}.
\]
\noindent Now assume that $N_k \rightarrow +\infty$ and we define the Euclidean frame $\mathcal{O}=(N_k,t_k,(x_k,y_k),0)_k $ and the function: 
\[ \psi= \mathcal{F}_{\mathbb{R}^4}^{-1}(|\xi|^{-2}[\eta^4(\xi)-\eta^4(2\xi)]) \in H^1(\mathbb{R}^4).
\]
\noindent By using Poisson Summation Formula, we can prove that
\[\lim\limits_{k \rightarrow +\infty} ||(1-\Delta)T_{N_k}\psi-N_k^{-1}P_{N_k}\delta_0||_{L^{\frac{4}{3}}}=0.
\]
\noindent Thus, by definition, we have $||(1-\Delta)T_{N_k}\psi-N_k^{-1}P_{N_k}\delta_0 ||_{H^{-1}}\rightarrow 0$ and then we conclude,
\[ \frac{1}{2}v \lesssim |\langle f_k,N_k^{-1}(e^{it_k\Delta}P_{N_k})\delta_{(x_k,y_k)}  \rangle| \lesssim |\langle f_k,(1-\Delta)\tilde{\psi}_{\mathcal{O},k} \rangle_{ H^1 \times H^{-1}} \rangle|.
\]
\noindent The last claim about orthogonality $\mathcal{O}$ with $\mathcal{O}^{\alpha}$ follows from Lemma 6.2 and the existence of a nonzero scalar product in (6.12).\vspace{3mm}

\noindent Now continuing with the sequence $\{ f_k\}_k$ as above, if the frame selected was a scale-one frame, we consider 
\[ g_k(x,y):=e^{it_k \Delta} f_k((x,y)+p_k)=\Pi_{-(t_k,p_k)} f_k.
\]
\noindent This is a bounded sequence in $H^1$, up to a subsequence, we assume it converges weakly to $\varphi \in H^1$. We then define the profile corresponding to $\mathcal{O}$ as $\tilde{\varphi}_{\mathcal{O},k}$. By its definition and (6.6), $\varphi$ has norm smaller than $E$. Also, we have,
\[\frac{v}{2} \lesssim \lim\limits_{k \rightarrow +\infty} \langle f_k, \tilde{\psi}_{\mathcal{O},k} \rangle_{H^1 \times H^1} \lesssim \lim\limits_{k \rightarrow +\infty} \langle g_k, \psi \rangle_{H^1 \times H^1}=\langle \varphi, \psi \rangle_{H^1 \times H^1}.
\]
\noindent Consequently, we get that 
\begin{equation}\label{equation} ||\varphi||_{H^1}  \gtrsim v. \end{equation}
\noindent We also observe that since $g_k-\psi$ weakly converges to $0$ in $H^1$, there holds that
\begin{equation}\label{equation} ||Af_k||_{L^2}^2 =||Ag_k||_{L^2}^2=||A g_k-\varphi||^2_{L^2}+||A\varphi||_{L^2}^2+o_k(1)=||A(f_k-\tilde{\varphi}_{\mathcal{O},k})||^2_{L^2}+||A\varphi||_{L^2}^2+o_k(1) \end{equation}
for $ A=1$ or $A=\nabla_{x,y}$.
\noindent The situation for Euclidean frame is similar.\vspace{3mm}

\noindent For the Euclidean case, for $k$ large enough, we consider
\[\varphi_k(y)=N_k^{-1}\eta^4(y/N_k)(\Pi_{(-t_k,-x_k)}f_k)(\Psi(y/N_k)),\quad y\in \mathbb{R}^4.
\]
\noindent This is a bounded sequence in $\dot{H}^1(\mathbb{R}^4)$. We can extract a subsequence that converges weakly to a function $\varphi \in \dot{H}^1(\mathbb{R}^4)$ satisfying
\[ ||\varphi||_{\dot{H}^1(\mathbb{R}^4)} \lesssim 1.
\]
\noindent Now, let $\gamma \in C_0^{\infty}(\mathbb{R}^4)$; for $k$ large enough,
\begin{align*}
\langle f_k,\tilde{\gamma}_{\mathcal{O},k} \rangle _{H^1\times H^1(\mathbb{R}^2\times \mathbb{T}^2)}&=\langle \Pi_{-(t_k,x_k)}f_k,T_{N_k}^e \gamma \rangle_{H^1\times H^1(\mathbb{R}^2\times \mathbb{T}^2)}\\
&=\langle \varphi,\gamma \rangle_{ \dot{H}^1 \times \dot{H}^1(\mathbb{R}^4)} +o_k(1).
\end{align*}
\noindent We conclude that 
\begin{equation}\label{equation} ||\varphi||_{\dot{H}^1(\mathbb{R}^4)} \gtrsim v \end{equation}
\noindent and that 
\[h_k=f_k-\tilde{\varphi}_{\mathcal{O},k}
\]
\noindent is absent from the frame $\mathcal{O}$. Using Lemma 6.2, we see that
\[ ||h_k||_{L^2}^2=||f_k||_{L^2}^2+o_k(1),
\]
\begin{align*}
||\nabla_{x,y} h_k||_{L^2}^2&=||\nabla_{x,y} f_k||_{L^2}^2+||\nabla \varphi||_{L^2}^2-2\langle \nabla f_k, \nabla \tilde{\varphi}_{\mathcal{O},k}\rangle_{L^2 \times L^2} \\
&=||\nabla_{x,y} f_k||_{L^2}^2-||\nabla \varphi||_{L^2}^2+o_k(1).
\end{align*}
\noindent Defining $f^0_k=\phi_k$ and for each $\alpha$, $f_k^{\alpha+1}=f_k^{\alpha}-\tilde{\phi}_{\mathcal{O}^\alpha,k}$ where is the profile given based on the considerations above; iterating this claim at most $\mathcal{O}(v^{-2})$ times and replacing $\phi_k$ by 
\[ \phi_k^{'}=\phi_k-\sum_{\alpha} \tilde{\varphi}_{\mathcal{O}^\alpha,k}.
\]
\noindent We obtain that $\{\phi_k^{'}\}_k$ satisfies
\begin{equation}\label{equation} \limsup\limits_{k \rightarrow +\infty} ||\phi_k^{'}||_{H^1} \leq E < +\infty \end{equation}
\noindent and
\begin{equation}\label{equation} \limsup\limits_{k \rightarrow +\infty}\sup\limits_{N\geq 1,t,x,y}N^{-1}|(e^{it\Delta}P_N\phi_k^{'})(x,y)| < v. \end{equation}
\noindent This proves the Lemma 6.4.\vspace{3mm}

\noindent \emph{Proof of Lemma 6.3:} First, for $v=v(\delta,E)$ to be decided later we use Lemma 6.5 and extract some profiles. Then, we replace $\phi_k$ by $\phi_k^{'}$, thus ensuring that (6.10) holds for the sequence $\{\phi_k^{'}\}_k$. We now consider 
\[ \Lambda_{0}(\{\phi_k^{'}\}) =\limsup\limits_{k \rightarrow +\infty} ||e^{it\Delta}\phi_k^{'}||_{Z(\mathbb{R})}.
\]\noindent If $\Lambda_{0}(\{\phi_k^{'}\}) <\delta$, we may set $R^A_k=\phi_k^{'}$ for all $k$ and we get Lemma 6.4. \vspace{3mm}

\noindent Now we claim that if $\{\phi_k^{'}\}_k$ satisfies $\Lambda_0(\{\phi_k^{'}\}_k) \geq \delta$ and $\{\phi_k^{'}\}_k$ is orthogonal to a family of frames $\mathcal{O}^{\alpha}$, $1 \leq \alpha \leq A$, then there exists a frame $\mathcal{O}$ orthogonal to $\mathcal{O}^{\alpha}$ and an associated profile $\tilde{\varphi}_{\mathcal{O},k}$ such that, after passing to a subsequence, we have that 
  \begin{equation}\label{equation} \limsup\limits_{k \rightarrow +\infty} ||\tilde{\varphi}_{\mathcal{O},k}||_{H^1(\mathbb{R}^2\times \mathbb{T}^2)} \gtrsim_{\delta} 1 \qquad \phi_k^{'}-\tilde{\varphi}_{\mathcal{O},k} \quad \textmd{is absent from} \quad \mathcal{O}. \end{equation}
\noindent Once the claim is established then the end of the proof follows by iterating the extraction process as in Lemma 6.4. Thus, now we will focus on this claim.\vspace{3mm}

\noindent Since $\Lambda_0(\{\phi_k^{'}\}) \geq \delta$, by H$\ddot{o}$lder's inequality and Strichartz estimates (3.2), we have that, 
\[ c_N^k=N^{\frac{1}{2}}||e^{it\Delta}P_N \phi_k^{'}||_{L^4(\mathbb{R}^2\times \mathbb{T}^2 \times \mathbb{R})},
\]
\[ ||c_N^k||_{l_N^{4}} \leq ||c_N^k||_{l_N^{2}}^{\frac{1}{2}}||c_N^k||_{l_N^{\infty}}^{\frac{1}{2}} \leq (\sum_{N \geq 1}N^2||P_N \phi_k^{'}||^2_{L^2})^{\frac{1}{2}}(\sup\limits_{N}c_N^k)^{\frac{1}{2}}.
\] 
\noindent Using (6.16), we obtain that there exists a sequence of scales $N_k \geq 1$ such that, 
\[ (\frac{\delta}{2})^{2}< \Lambda_0(\{\phi_k^{'}\})^{2} \leq E^{\frac{1}{2}}N_k^{\frac{1}{2}}||e^{it\Delta}P_{N_k} \phi_k^{'}||_{L^4(\mathbb{R}^2\times \mathbb{T}^2 \times \mathbb{R}).}
\]
\noindent We conclude that there exists a sequence $h_k \in C_c^{\infty}(\mathbb{R}^2 \times \mathbb{T}^2 \times \mathbb{R})$ such that 
\[ 1 \leq ||h_k||_{L^{\frac{4}{3}}}  \leq 2,
\]
\[(\frac{\delta}{2})^{2}E^{-\frac{1}{2}}N_k^{-\frac{1}{2}} \leq \langle h_k, e^{it\Delta}P_{N_k}\phi_k^{'} \rangle_{L^2 \times L^2}.
\]
\noindent Now for a given threshold $B$, we introduce the partition function $ \chi_B(\gamma) $ satisfying \vspace{3mm}

\noindent $ \chi_B(\gamma)=1 $  if  $||h_k||_{L^{\frac{4}{3}}(\mathbb{R}^2\times \mathbb{T}^2 \times I_{\gamma})} \geq B,$  $\quad \chi_B(\gamma)=0$  otherwise. \vspace{3mm}

\noindent And we decompose as follows, 
\[ h_k(x,y,t)=h_k^{>B}+h_k^{<B}=h_k(x,y,t) \chi_B([\frac{t}{2\pi}])+h_k(x,y,t)(1- \chi_B([\frac{t}{2\pi}]))
\]
\noindent so we have 
\[||h_k||_{L^{\frac{4}{3}}(\mathbb{R}^2\times \mathbb{T}^2 \times \mathbb{R})} \leq ||h_k^{>B}||_{L^{\frac{4}{3}}(\mathbb{R}^2\times \mathbb{T}^2 \times \mathbb{R})}+||h_k^{<B}||_{L^{\frac{4}{3}}(\mathbb{R}^2\times \mathbb{T}^2 \times \mathbb{R})},
\]
\[ \sup\limits_{\gamma} ||h_k^{<B}||_{L^{\frac{4}{3}}(\mathbb{R}^2\times \mathbb{T}^2 \times I_{\gamma})} \leq B.
\]
\noindent Using Strichartz estimates, we have that for $\frac{10}{3}<p_1 < 4$:
\[ \limsup\limits_{k \rightarrow +\infty}N_k^{\frac{6}{p_1}-1}||e^{it\Delta}P_{N_k}\phi_k^{'}||_{l_{\gamma}^{\frac{2p_1}{p_1-2}} L^{p_1}} \lesssim \limsup\limits_{k \rightarrow +\infty} ||\phi_k^{'}||_{H^1} \lesssim E.
\]
\noindent Interpolating with (6.17), we obtain that 
\[\limsup\limits_{k \rightarrow +\infty}N_k^{\frac{1}{2}}||e^{it\Delta}P_{N_k}\phi_k^{'}||_{l^{\frac{8}{(p_1-2)}}_{\gamma}L^4} \lesssim E^{\frac{p_1}{4}}v^{\frac{4-p_1}{4}}.
\]
\noindent An observation is that: 
\[ |supp_{\gamma}(h_k^{>B})| \leq (\frac{2}{B})^{\frac{4}{3}}.
\]
\noindent Thus by H$\ddot{o}$lder's inequality in $\gamma$,
\begin{align*}
\langle e^{it\Delta}P_{N_k}\phi_k^{'},h_k^{>B} \rangle &\leq ||e^{it\Delta}P_{N_k} \phi_k^{'}||_{l^{\frac{8}{(p_1-2)}}_{\gamma}L^4}||h_k^{>B}||_{L^{\frac{4}{3}}} [(\frac{2}{B})^{\frac{4}{3}}]^{\frac{4-p_1}{8}} \\ 
&\lesssim (\frac{1}{B})^{\frac{4-p_1}{6}}N_k^{-\frac{1}{2}}E^{\frac{p_1}{4}}v^{\frac{4-p_1}{4}}.
\end{align*}
\noindent Eventually, for any fixed $B>0$, we can choose $v=v(B,\delta)$ such that
\[v^{\frac{4-p_1}{4}}=c\delta^2 E^{-\frac{1}{2}-\frac{p_1}{4}}B^{\frac{4-p_1}{6}},
\]
\noindent for $c>0$ a constant sufficiently small so that,
\begin{equation}\label{equation}
\aligned
(\frac{\delta}{2})^{2}E^{-\frac{1}{2}} & \leq 2\langle h_k^{<B},e^{it\Delta} P_{N_k} \phi_k^{'}\rangle_{L^2\times L^2}\\
& \leq  2 EN_k^{-1} ||\int_{\mathbb{R}}e^{-is\Delta}P_{N_k}h_k^{<B}(x,y,s) ds||_{L^2_{x,y}(\mathbb{R}^2 \times \mathbb{T}^2)}.\\
\endaligned
\end{equation}
\noindent Using Strichartz estimate (3.3), we obtain that 
\begin{align*}
||\int_{\mathbb{R}} e^{-is\Delta}P_{N_k}h_k^{<B}(x,y,s)ds||_{L^2_{x,y}}  &\lesssim N_k^{\frac{1}{2}}||h_k^{<B}||_{l^2_{\gamma}L^{\frac{4}{3}}_{x,y,t}(I_{\gamma})}+N_k^{\frac{1}{3}}||h_k^{<B}||_{L^{\frac{4}{3}}_{x,y,t}(I_{\gamma})}\\
&\lesssim N_{k}^{\frac{1}{2}}||h_k^{<B}||_{L^{\frac{4}{3}}_{x,y,t}(I_{\gamma})}^{\frac{2}{3}}B^{\frac{1}{3}}+N_k^{\frac{1}{3}}\\
&\lesssim N_k^{\frac{1}{2}}B^{\frac{1}{3}}+N_k^{\frac{1}{3}}.
\end{align*}
\noindent Choosing $B$ such that 
\[B^{\frac{1}{3}}=\epsilon E^{-\frac{3}{2}} \delta^{2}.
\]
\noindent for some constant $\epsilon>0$ small enough and plugging into (6.19), we obtain that
\[ \delta^{2} \lesssim \epsilon \delta^{2}+E^{\frac{3}{2}} N_k^{-\frac{1}{6}}.
\]
\noindent If $\epsilon>0$ is small enough, from the estimate above we can obtain a uniform bound for $N_k$ and the bound only relies on $E$ and $\delta$
\begin{equation}\label{equation} N_k \lesssim_{E,\delta} 1 .
\end{equation}
\noindent To sum  up, we have showed that if $\Lambda_0(\{\phi_k^{'}\})>\delta$, then there exists a sequence of scales $N_k$ satisfying (6.20) and such that 
\[ ||P_{N_k}e^{it\Delta}\phi_k^{'}||_{L^4_{x,y,t}(\mathbb{R}^2\times \mathbb{T}^2 \times \mathbb{R})}> c(\delta,E)
\]\noindent for some $c(\delta,E)>0$ which denotes a small positive constant depending only on $\delta$ and $E$. \vspace{3mm}

\noindent Now, we write
\[ P_{N_k}e^{it\Delta_{\mathbb{R}^2\times \mathbb{T}^2}}\phi_k^{'}(x,y)=\sum_{z=(z_1,z_2)\in \mathbb{Z}^2 ;|z_i| \leq N_k} e^{-it|z|^2} e^{i\langle z,y\rangle} \eta^2_{N_k}(z)[e^{it\Delta_{\mathbb{R}^2}}\phi^{''}_{z,k}(x)]
\]
\noindent where
\[\phi_{z,k}^{''}(x)=\frac{1}{(2\pi)^2} \int_{\mathbb{T}^2} P_{N_k}\phi^{'}_k(x,y) e^{-i\langle z,y \rangle} dy.
\]
\noindent Extracting a subsequence, we conclude that there exists $z$ such that, for all k,
\[supp(\mathcal{F}_x\phi_{z,k}^{''}) \subset [-3N_k,3N_k]^2,
\]
\begin{equation}\label{equation}
||e^{it\Delta_{\mathbb{R}^2}}\phi^{''}_{z,k}||_{L^4_{x,t}(\mathbb{R}^2\times \mathbb{R})} >c(\delta,E), 
\end{equation}
\[||\phi_{z,k}^{''} ||_{H^1(\mathbb{R}^2)} \leq M.
\]
\noindent Now we need to apply the following result from [1, 5, 29]. (We use the version as [29, Corollary A.3])
\begin{theorem}\label{theorem}For any $M$, $c(\delta, E) > 0$, there exists a finite set $\mathcal{C} \subset L^2(\mathbb{R}^2)$ of functions satisfying,
\begin{equation}\label{equation}
||v||_{L^2(\mathbb{R}^2)}=1 ,\quad ||v||_{L^1(\mathbb{R}^2)} \leq S(E,M,\delta) ,  \quad \forall v \in \mathcal{C}.\end{equation}
\noindent and $\kappa(M, c(\delta, E))>0$ such that whenever $u \in  L^{2}(\mathbb{R}^2)$ obeys the bounds
\[ ||u||_{L^2(\mathbb{R}^2)} \leq M,\quad ||e^{it\Delta_{\mathbb{R}^2}} u||_{L^4_{x,t}(\mathbb{R}^2 \times \mathbb{R})}  \geq c(\delta, E).
\]
\noindent Then there exists $v \in \mathcal{C}$ and $(\lambda, \xi_{0},t_0,x_0) \in \mathbb{R}_{+} \times \mathbb{R}^2 \times \mathbb{R} \times \mathbb{R}^2$ such that
\[\langle u,v^{'}\rangle_{L^2 \times L^2(\mathbb{R}^2)} \geq \kappa ,\quad  v^{'}(x)= \lambda e^{ix\xi_{0}}[e^{it_{0}\Delta_{\mathbb{R}^2}} v](\lambda (x-x_{0})).
\]\end{theorem}
\noindent Using this Theorem, after extraction, there exists a function $v \in L^2(\mathbb{R}^2)$ satisfying (6.22) and a sequence $(\lambda_k, \xi_{k},t_k,x_k) \in \mathbb{R}_{+} \times \mathbb{R}^2 \times \mathbb{R} \times \mathbb{R}^2$ such that 
\[\langle \phi_k^{''},v_k\rangle \geq \kappa, \quad  v_k(x)= \lambda_k e^{ix\xi_{k}}[e^{it_{k}\Delta_{\mathbb{R}^2}} v](\lambda_k (x-x_{k})) \tag{6.23}.\]
\noindent We may assume $v$ has a compactly supported Fourier transform. Also we claim that $\lambda_k$ and $|\xi_k|$ remain bounded. By computation,
\[ \mathcal{F}_{\mathbb{R}^2}v_k(\xi)=\lambda_k^{-1} e^{ix_k \xi_k}e^{-ix_k \xi}[\mathcal{F}_{\mathbb{R}^2} e^{it_k \Delta_{\mathbb{R}^2}} v](\frac{\xi-\xi_k}{\lambda_k}),
\]
\[ ||\langle \xi \rangle ^{-1}\mathcal{F}_{\mathbb{R}^2} v_k(\xi)||_{L^2([-4N_k,4N_k]^2)} \lesssim \lambda_k^{-1}S(M,E,\delta) ,
\]
\[ ||v_k||^2_{H^{-1}(\mathbb{R}^2)} \sim \int_{\mathbb{R}^2} [\frac{1}{1+|\xi_k+\lambda_k \eta|}]^2 |\mathcal{F}_{\mathbb{R}^2} v|^2(\eta) d\eta.
\]\noindent Together with (6.21) and (6.23), we can show $\lambda_k$ and $|\xi_k|$ remain bounded. \vspace{3mm}

\noindent Assume first that $\lambda_k$ remains bounded below. Then, up to extracting a subsequence, we may assume that $\lambda_k \rightarrow \lambda_{\infty} \in (0,\infty)$. Similarly, we may assume that  $\xi_{k} \rightarrow \xi_{\infty}$. Then, setting
\[\tilde{\psi}(x,y)=e^{i\langle z,y \rangle}e^{ix\xi_{\infty}}\lambda_{\infty}v(\lambda_{\infty}x) ,\qquad t_k^{'}=-\lambda_k^{-2}t_k,\qquad x^{'}_k=x_k+2t_k^{'}\xi_k\]
\noindent and defining the frame $\mathcal{O}=\{(1,t_k^{'},(x_k^{'},0),0)_k\}$, we see from (6.23) that 
\begin{align*} 
\kappa &\leq \langle P_{N_k}\phi_k^{'}, e^{i\langle z,y \rangle}e^{ix\xi_{k}}\lambda_{k}[e^{-it_k \Delta_{\mathbb{R}^2}}v](\lambda_k(x-x_k)) \rangle \\ &\leq \langle \phi_k^{'},e^{i(x_k\xi_{\infty}-t^{'}_k(|z|^2-|\xi_{\infty}|^2))}P_{N_k}e^{-it_k^{'}\Delta_{\mathbb{R}^2\times \mathbb{T}^2}}\tilde{\psi}(x-x_k^{'}) \rangle+o_k(1) .
\end{align*}
\noindent Since $1 \leq N_k \lesssim 1$, up to a subsequence, we may assume that $N_k=N_{\infty}$. As a result, setting $\psi=P_{N_{\infty}}\tilde{\psi}$, we see the scale-one profile $\tilde{\psi}_{\mathcal{O},k}$ satisfies (6.11) and (6.12). We also conclude that $\mathcal{O}$ is orthogonal to $\mathcal{O}^{\alpha}$, $1\leq \alpha \leq A$. As in the proof of Lemma 6.4, we find $\varphi$ satisfying (6.18).\vspace{3mm}

\noindent Assume now that $\lambda_k \rightarrow 0$. Let $M_k$ be a dyadic number such that $1\leq \lambda_k^{-1} M_k \leq 2$ and consider the sequence
\[ \Phi_k(x,y)=M_k^{-1} e^{-it_k\Delta_x} [e^{i\xi_k x/M_k}\phi^{'}(x_k+\frac{x}{M_k},y)].
\]
\noindent We have that 
\[ ||\Phi_k||^2_{L^2(\mathbb{R}^2\times \mathbb{T}^2)}+||\nabla_y \Phi_k||^2_{L^2(\mathbb{R}^2\times \mathbb{T}^2)} \leq M^2+E^2.
\]
\noindent Up to a subsequence, we may assume that $\Phi_k \rightharpoonup \Phi$ in $H^{0,1}(\mathbb{R}^2\times \mathbb{T}^2)$. We define 
\[\qquad \qquad \qquad t_k^{'}=-M_k^{-2}t_k,\qquad \xi_k^{'}=-\xi_k,\qquad x_k^{'}=x_k+2t_k^{'}\xi_k^{'},\]
\noindent and $\mathcal{O}=(M_k,t^{'}_k,(x^{'}_k,0),\xi_k^{'})$. Then we obtain that $\mathcal{O}$ is orthogonal to $\mathcal{O}^{\alpha}$ from (6.23) and we see from the definition of $\Phi_k$ that (6.18) holds with $\varphi=e^{i\theta_{\infty}} \Phi$ for some $\theta_{\infty} \in \mathbb{R}/\mathbb{Z}$.\vspace{3mm}

\noindent This finishes the proof of Lemma 6.3.
\section{Induction on Energy}
\noindent Then we are now ready to prove the main theorem. We follow an induction on energy method formalized in [20, 21]. Define
\[ \Lambda(L)=\textmd{sup}\{||u||_{Z(I)}:u \in X^1_{loc}(I),E(u)+\frac{1}{2} M(u) \leq L\}
\]
\noindent where the supremum is taken over all strong solutions of full energy less than $L$. By the local theory, this is sublinear in $L$ and finite for $L$ sufficiently small. We also define
\[L_{max}= \textmd{sup}\{ L: \Lambda(L) < +\infty \}.
\]\noindent We want to show that $ L_{max}=+\infty.$ That is our goal. If $ L_{max}=+\infty $ holds, we can extend the small data result to our main theorem, i.e. Theorem 1.2 by using the Extension Theorem (Theorem 4.4). \vspace{3mm}

\noindent The key proposition is: 
\begin{theorem}\label{theorem} Assume that $L_{max} < +\infty $ and the Conjecture 1.2 holds for $L_{max}$. Let $\{ t_{k}\}_{k} ,\{ a_{k}\}_{k} ,\{ b_{k}\}_{k} $ be arbitrary sequences of real numbers and $\{u_k\}$ be a sequence of solutions to (1.1) such that $u_{k} \in X^{1}_{c,loc}(t_k-a_k,t_k+b_k)$ and satisfying
\begin{equation}\label{equation} L(u_k) \rightarrow L_{max},\quad ||u_k||_{Z(t_k-a_k,t_k)} \rightarrow +\infty,\quad ||u_k||_{Z(t_k,t_k+b_k)} \rightarrow +\infty. \end{equation}
Then passing to a subsequence, there exists a sequence $x_{k} \in \mathbb{R}^2 $ and $\omega \in H^{1}(\mathbb{R}^2 \times \mathbb{T}^2)$ such that
\begin{equation}\label{equation} \omega_k(x,y) =u_k(x-x_k,y,t_k) \rightarrow \omega \end{equation}
strongly in $H^{1}(\mathbb{R}^2 \times \mathbb{T}^2)$.
\end{theorem}
\noindent We will give the proof later by using Theorem 6.2 (profile decomposition). Based on Theorem 7.1, we can prove:
\begin{corollary} Assume that $L_{max} < +\infty $ and the Conjecture 1.2 holds for some $E_{max}^{ls} \geq L_{max}$. Then there exists $u \in X^1_{c,loc}(\mathbb{R})$ solving (1.1) and a Lipschitz function $\underline{x}:\mathbb{R} \rightarrow \mathbb{R}^2$ such that $L(u)=L_{max}$ and
\begin{equation}\label{equation} \sup\limits_{t \in \mathbb{R}} |\underline{x}^{'}(t)| \lesssim 1 \end{equation} 
\[ (u(x-\underline{x}(t),y,t) :t \in \mathbb{R})\quad \textmd{is precompact in} \quad H^1(\mathbb{R}^2 \times \mathbb{T}^2).
\]
\end{corollary}
\noindent \emph{Proof:} Assuming that $L_{max} < +\infty$, we can find a sequence of solutions of (1.1) $u_k$ satisfying (7.1). Apply Theorem 7.1 we can extract a subsequence and obtain a sequence $x_k$ such that Corollary 7.2 holds for some $\omega \in H^1(\mathbb{R}^2 \times \mathbb{T}^2)$.
\noindent Thus $L(\omega)=L_{max}<+\infty$. Let $U \in C(I:H^1)$ be the maximal strong solution of (1.1) with initial data $\omega$, defined on $I=(-a_{\infty},b_{\infty})$. According to local theory and (7.1),
\begin{equation}\label{equation} ||U||_{Z(-a_{\infty},0)}=||U||_{Z(0,b_{\infty})}=+\infty. \end{equation}
\noindent We claim that there exists $\kappa >0$ such that for all $t\in I$,
\begin{equation}\label{equation} ||U||_{Z(t-2\kappa,t+2\kappa)} \leq 2. \end{equation}
\noindent If $U$ is global, $a_{\infty}=b_{\infty}=\infty$.\vspace{3mm}

\noindent Assume if (7.5) is not true. Then there exists a sequence $t_k \in I$ such that 
\[ ||U||_{(t_k-\frac{1}{k},t_k+\frac{1}{k})}\geq 2.
\]\noindent We can apply Theorem 7.1 to the sequence $U(t_k)$ and obtain that, up to a subsequence, there exists $x_k$ such that $U_k(x,y)=U(t_k,x-x_k,y) \rightarrow \omega^{'}$ in $H^1$. Let $W$ be the nonlinear solution of (1.1) with initial data $\omega^{'}$. By the local theory, we know there exists $\kappa^*>0$ such that
\[ ||W||_{Z(-k^*,k^*)} \leq 1
\]
\noindent and by the stability theory, we obtain that, for $k$ large enough
\[ ||U||_{Z(t_k-\kappa^*,t_k+\kappa^*)} \leq 2,
\]
\noindent which gives a contradiction for $k$ large enough.\vspace{3mm}

\noindent Now we prove (7.3). We define the sequence of times $t_k=k\kappa$ and for each $t_k$, we define $\underline{x}_k$ and $R_k$ such that 
\begin{equation}\label{equation} \frac{1}{2}\int_{\{|x-\underline{x}_k|\leq R_k\}}\int_{\mathbb{T}^2}[|u(t_k,x,y)|^2+|\nabla u(t_k,x,y)|^2+\frac{1}{2}|u(t_k,x,y)|^4]dxdy=\frac{99}{100} L_{max} \end{equation}
\noindent and $R_k$ is the minimal with this property. While $\underline{x}_k$ is not necessarily unique, we claim that there exists $D$ such that for all $k$
\begin{equation}\label{equation} R_k \leq D, \quad |\underline{x}_{k}-\underline{x}_{k+1}| \leq D \end{equation}
\noindent and that
\begin{equation}\label{equation} \{ u(t_k+s,x-\underline{x}_k),k \in \mathbb{Z},s \in (-\kappa,\kappa)\} \quad \textmd{is precompact in}\quad H^1(\mathbb{R}^2\times \mathbb{T}^2). \end{equation}
\noindent The fact that the $R_k$ are uniformly bounded comes from the compactness up to translations of $\{u(t_k)\}_{k}$. Assume that $\{v_k(x,y)=u(x-\underline{x}_k,y,t_k)\}$ was not precompact in $H^1(\mathbb{R}^2\times \mathbb{T}^2)$. In that case, there exists $\epsilon >0$ and a subsequence $k^{'}$such that for all $k_1^{'}$, $k_2^{'}$,
\begin{equation}\label{equation} ||v_{k_1^{'}}-v_{k_2^{'}}||_{H^1(\mathbb{R}^2\times \mathbb{T}^2)}>\epsilon.
\end{equation}
\noindent Now apply Theorem 7.1, we see that there exists a sequence $\overline{x}^k$ and a subsequence of $k^{'}$ such that
\[ \qquad  v_{k^{''}}(x-\overline{x}^{k^{''}},y) \rightarrow \omega(x,y)  \quad \textmd{strongly in} \quad H^1(\mathbb{R}^2\times \mathbb{T}^2).\]
From (7.6), $\{\overline{x}^k \}_k$ remains bounded, so that the convergence of $v_{k^{''}}$ contradicts (7.9). Using (7.5) and the precompactness of $\{v_k\}_k$, we obtain (7.8). Similarly, this implies the second statement in (7.7). Choose $ \underline{x}(t)$ to be a Lipschitz function satisfying $ \underline{x}(t_k)=\underline{x}_k$, we obtain (7.3). This completes the proof.\vspace{3mm}

\noindent Now we can finish the proof of the main theorem by using a Morawetz-type argument as follows.
\begin{theorem}\label{theorem}Assume that $u$ satisfies the conclusion of Corollary 7.2, then $u=0$. In particular, $L_{max}\geq E^{ls}_{max}$.
\end{theorem}
\noindent \emph{Proof:} Assume $u \neq 0$. Then, from the compactness property, we see that there exists $\rho >0$ such that 
\begin{equation}\label{equation} \inf\limits_{t \in \mathbb{R}}min (||u(t)||_{L^4_{x,y}(\mathbb{R}^2 \times \mathbb{T}^2)},||u(t)||_{L^2_{x,y}(\mathbb{R}^2 \times \mathbb{T}^2)}) \geq \rho.
\end{equation}
\noindent Now we consider the conserved momentum 
\[ P(u)=Im \int_{\mathbb{R}^2\times \mathbb{T}^2}\bar{u}(x,y,t) \partial_{x_1} u(x,y,t) dxdy.
\]
\noindent Considering the Galilean transform 
\[v(z,t)=e^{-i|\xi_0|^2t+i \langle z,\xi_0 \rangle} u(z-2\xi_0 t,t),
\]
\noindent let
\[\xi_0=-\frac{P(u)}{M(u)},
\]
\noindent without loss of generality, we can assume that 
\begin{equation}\label{equation} P(u)=0.
\end{equation} \noindent Then we define the Virial action by 
\[ A_{R}(t)=\int_{\mathbb{R}^2\times \mathbb{T}^2}  \chi_{R}(x_1-\underline{x}_1(t))(x_1-\underline{x}_1(t))Im[\bar{u}(x,y,t)\partial_{x_1} u(x,y,t)]dxdy
\]
\noindent for $\chi_R(x)=\chi(x/R)$ where $\chi$ satisfies $\chi(x)=1$ when $|x| \leq 1$ and $\chi(x)=0$ when $|x| \geq 2$ ($x \in \mathbb{R}$).\vspace{3mm}

\noindent On one hand, clearly 
\begin{equation}\label{equation} \sup\limits_{t}|A_{R}(t)| \lesssim R .
\end{equation}
\noindent On the other hand, we compute that
\begin{align*}
\frac{d}{dt}A_R &=- \underline{x}_1^{'}(t) Im \int_{\mathbb{R}^2\times \mathbb{T}^2}\bar{u}(x,y,t) \partial_{x_1}  u(x,y,t) dxdy \\
& - \underline{x}_1^{'}(t)\int_{\mathbb{R}^2\times \mathbb{T}^2} \{(\chi ^{'})_R(x_1-\underline{x}_1(t))\frac{x-\underline{x}_1(t)}{R} -(1- \chi_{R}(x_1-\underline{x}_1(t)))\} Im [\bar{u}(x,y,t) \partial_{x_1}  u(x,y,t)] dxdy\\
&+  \int_{\mathbb{R}^2\times \mathbb{T}^2} \chi_{R}(x_1-\underline{x}_1(t))(x_1-\underline{x}_1(t)) \partial_t Im[\bar{u}(x,y,t)\partial_{x_1} u(x,y,t)]dxdy.
\end{align*}

\noindent The first term will vanish automatically based on the assumption (7.11) and the second term can be bounded by  
\[ \int_{\{|x-\underline{x}(t)| \geq R\}} \int_{\mathbb{T}^2}[|u(x,y,t)|^2+|\nabla u(x,y,t)|^2]dxdy=O_{R}(t),
\]
 $\qquad \qquad \qquad  \qquad  \qquad  \qquad  \sup\limits_t O_{R}(t) \rightarrow 0 \qquad$ as $\qquad  \quad R  \rightarrow +\infty $ .

\noindent Notice that
\[  \partial_t Im[\bar{u}(x,y,t) \partial_{x_1} u(x,y,t)]=\partial_{x_1} \Delta \frac{|u|^2}{2}-2div\{Re[\partial_{x_1}\bar{u} \nabla u]\}-\frac{1}{4}\partial_{x_1} |u|^4.
\]
\noindent For the last term, we have 
\begin{align*}
\frac{d}{dt}A_R &=\int_{\mathbb{R}^2\times \mathbb{T}^2}\chi_{R}(x_1-\underline{x}_1(t))[\frac{1}{4}|u(x,y,t)|^4+\frac{1}{2}|\partial_{x_1} u(x,y,t)|^2]dxdy \\
&+ \int_{\mathbb{R}^2\times \mathbb{T}^2}\chi_{R}^{'}(x_1-\underline{x}_1(t)) \frac{x_1-\underline{x}_1(t)}{R}[\frac{1}{4}|u(x,y,t)|^4+\frac{1}{2}|\partial_{x_1} u(x,y,t)|^2]dxdy\\
&-  \int_{\mathbb{R}^2\times \mathbb{T}^2} \frac{|u(x,y,t)|^2}{2} \partial_{x_1}^3[ \chi_{R}(x_1-\underline{x}_1(t))(x_1-\underline{x}_1(t)) ]dxdy +O_R(t)\\
&= \int_{\mathbb{R}^2\times \mathbb{T}^2} [\frac{1}{4}|u(x,y,t)|^4+\frac{1}{2}|\partial_{x_1} u(x,y,t)|^2]dxdy+\tilde{O}_{R}(t).
\end{align*}
\noindent Integrating this equality, we obtain 
\[|A_R(t)-A_R(0)| \geq Ct\rho-t \sup\limits_t \tilde{O}_R(t).
\]
\noindent Taking $R$ sufficiently large enough, we obtain, when $t$ is sufficiently large, there is a contradiction. This finishes the proof of Theorem 7.3.\vspace{3mm}

\noindent \textbf{Proof of Theorem 7.1:} Without loss of generality, we may assume that $t_k=0$. We apply the profile decomposition, i.e. Theorem 6.3 to the sequence $\{u_k(0)\}$ which is  bounded in $H^1(\mathbb{R}^2 \times \mathbb{T}^2)$. Then we get:
\[ u_k(0)=\sum_{1 \leq \alpha \leq J}\tilde{\varphi}^{\alpha}_{\mathcal{O}^\alpha,k} +\sum_{1 \leq \beta \leq J}\tilde{\psi}^{\beta}_{\mathcal{S}^\beta,k} +\sum_{1 \leq \gamma \leq J}\tilde{W^{\gamma}}_{\mathcal{O}^\gamma,k} +R_k^J(x,y).
\]
\noindent There are three cases to be discussed: no profile, one profile and multiple profiles. \vspace{3mm}

\noindent \textbf{Case 1:} There is no profiles. So if we take $J$ sufficiently large, we will have:
\[ ||e^{it\Delta} u_{k}(0)||_{Z(\mathbb{R})} = ||e^{it\Delta} R^{J}_{k}||_{Z(\mathbb{R})}\leq  \delta_{0}/2
\]\noindent for k sufficiently large, where $\delta_{0}$ is given in Theorem 4.3. Then we know that $u_k$ can be extended on $\mathbb{R}$ and that
\[ \lim\limits_{k \rightarrow +\infty} ||u_k||_{Z(\mathbb{R})} \leq \delta_{0}.
\]
\noindent It is a contradiction. Hence, we consider the situation when there are at least one profile. \vspace{3mm}

\noindent Now for every linear profile, we define the associated nonlinear profile as the maximal solution of (1.1) with the same initial data as in [13].\vspace{3mm}

\noindent From Section 5 and Section 6, we can define:
\[ L_E(\alpha):=\lim\limits_{k \rightarrow +\infty}(E(\tilde{\varphi}^{\alpha}_{\mathcal{O}^\alpha,k})+\frac{1}{2}M(\tilde{\varphi}^{\alpha}_{\mathcal{O}^\alpha,k}))=E_{\mathbb{R}^4}(\phi^{\alpha}) \in (0,L_{max}],
\]
\[L_{ls}(\beta):=\lim\limits_{k \rightarrow +\infty}(E(\tilde{\psi}^{\beta}_{\mathcal{S}^\beta,k} )+\frac{1}{2}M(\tilde{\psi}^{\beta}_{\mathcal{S}^\beta,k} ))=||\psi^{\beta}||^2_{H^{0,1}(\mathbb{R}^2\times \mathbb{T}^2)} \in (0,L_{max}],
\]
\[L_1(\gamma):=\lim\limits_{k \rightarrow +\infty}(E(\tilde{W^{\gamma}}_{\mathcal{O}^\gamma,k} )+\frac{1}{2}M(\tilde{W^{\gamma}}_{\mathcal{O}^\gamma,k} ))=E(W)+\frac{1}{2}M(W^{\gamma}) \in (0,L_{max}].
\]
\noindent Noticing that:
\begin{equation}\label{equation} \lim\limits_{J \rightarrow +\infty}[\sum_{1 \leq \alpha,\beta,\gamma \leq J}[L_{E}(\alpha)+L_{ls}(\beta)+L_{1}(\gamma)]+\lim\limits_{k \rightarrow +\infty} L(R_k^J)] \leq L_{max}. \end{equation}
\textbf{Case 2a:} $L_{E}(1)=L_{max}$, there is only one Euclidean profile, that is
\[ u_k(0)=\tilde{\phi}_{\varepsilon,k}+o_k(1)
\]
\noindent in $H^1(\mathbb{R}^2 \times \mathbb{T}^2)$, where $\varepsilon$ is a Euclidean frame. In this case, since the corresponding nonlinear profile $U_k$ satisfies $||U_k||_{Z(\mathbb{R})} \lesssim_{E_{\mathbb{R}^4}(\phi)} 1 $ and $  \lim\limits_{k \rightarrow +\infty } ||U_k(0)-u_k(0)||_{H^1(\mathbb{R}^2 \times \mathbb{T}^2)} \rightarrow 0 $.\vspace{3mm}

\noindent We may use Theorem 4.5 to deduce that 
\[ ||u_k||_{Z(\mathbb{R})} \lesssim ||u_k||_{X^1(\mathbb{R})} \lesssim_{L_{max}} 1,
\]
\noindent which contradicts (7.1).\vspace{3mm}

\noindent \textbf{Case 2b:} $L_{ls}(1)=L_{max}$, there is only one large scale profile, that is 
\[ u_k(0)=\tilde{\psi}_{S,k}+o_k(1)
\]
\noindent in $H^1$, where $S$ is a large-scale frame. Thus, the corresponding nonlinear profile $U_k$ satisfies $||U_k||_{Z(\mathbb{R})} \lesssim_{||\psi||_{H^{0,1}}} 1 $ and $ \lim\limits_{k \rightarrow +\infty } ||U_k(0)-u_k(0)||_{H^1(\mathbb{R}^2 \times \mathbb{T}^2)} \rightarrow 0. $ \vspace{3mm}

\noindent We may use Theorem 4.5 to deduce that
\[ ||u_k||_{Z(\mathbb{R})} \lesssim ||u_k||_{X^1(\mathbb{R})} \lesssim_{L_{max}} 1
\]
\noindent which contradicts (7.1). \vspace{3mm}

\noindent \textbf{Case 2c:} $L_{1}(1)=L_{max}$, there is only one scale-one profile, we have that 
\[ u_k(0)=\tilde{\omega}_{\mathcal{O},k}+o_k(1)
\]
in $H^1(\mathbb{R}^2\times \mathbb{T}^2)$, where $\mathcal{O}=\{1,t_k,x_k,0\}$ is a scale-one frame. If $t_k \equiv 0$, this is precisely the conclusion(7.2). \vspace{3mm}

\noindent If $t_k \rightarrow +\infty$, then 
\[||e^{it\Delta_{\mathbb{R}^2\times \mathbb{T}^2}}\tilde{\omega}_{\mathcal{O},k}||_{Z(a_k,0)} \leq ||e^{it\Delta_{\mathbb{R}^2\times \mathbb{T}^2}}\tilde{\omega}_{\mathcal{O},k}||_{Z(-\infty,0)}=||e^{it\Delta_{\mathbb{R}^2\times \mathbb{T}^2}}\omega||_{Z(-\infty,-t_k)}
\]
\noindent which goes to 0 as $t_k \rightarrow +\infty$. Using Theorem 4.3, we see that, for $k$ large enough, 
\[ ||u_k||_{Z(-\infty,0)} \leq \delta_0.
\]
\noindent It contradicts (7.1). The case $t_k \rightarrow  -\infty$ is similar.\vspace{3mm}

\noindent \textbf{Case 3:} $L_{\mu}(1)< L_{max}$ for all $\mu \in \{ E,ls,1\}$. In this case, we construct an approximate solution of (1.1) with initial data $u_k(0)$ whose $Z$-norm is finite and derive a contradiction by using Theorem 4.5. For this case, first, there exists $\eta >0$ such that for all $\alpha \geq 1$, $\mu \in \{ E,ls,1\}$ and $L_{\mu}(\alpha)<L_{max}-\eta$, we have that all nonlinear profiles are global and satisfy, for any $k$, $\alpha \geq 1$ and $\mu \in \{E,ls,1\}$ (after extracting a subsequence)
\[||U^{\mu,\alpha}_{k}||_{Z(\mathbb{R})} \leq 2\Lambda(L_{max}-\eta/2) \lesssim 1.
\]
One thing we need to mention is that from now on all implicit constants are allowed to depend on $\Lambda(L_{max}-\eta/2)$. Using Theorem 4.5, it follows that 
\begin{equation}\label{equation}  ||U^{\mu,\alpha}_{k}||_{X^1(\mathbb{R})} \lesssim 1.
\end{equation}

\noindent For $J,k \geq 1$, we define
\begin{equation}\label{equation} U^J_{prof,k} = \sum_{1 \leq \alpha \leq J} \sum_{\mu \in \{ E,ls,1\}} U^{\mu,\alpha}_{k}. 
\end{equation} \noindent We can prove that $||U^J_{prof,k}||_{X^1(\mathbb{R})} \lesssim 1.$ \vspace{3mm}

\noindent More precisely, we can show that there exists a constant $Q\lesssim 1$ such that 
\begin{equation}\label{equation} ||U^J_{prof,k}||_{X^1}^2+\sum_{1 \leq \alpha \leq J} \sum_{\mu \in \{ E,ls,1\}}||U^{\mu,\alpha}_{k}||_{X^1}^2 \leq Q^2 \end{equation}
\noindent uniformly in $J$, for all $k$ sufficiently large. Let $\delta_0(2L_{max})$ defined in Theorem 4.3. We know that there are only finitely many profiles such that $L(\alpha)\geq \frac{\delta_0}{2}$. Without loss of generality, we may assume that for all $\alpha \geq A, L(\alpha) \leq \delta_0$. Noticing (7.13),
\begin{align*}
||U^J_{prof,k}||_{X^1(\mathbb{R})}&=||\sum_{1 \leq \alpha \leq J} \sum_{\mu \in \{ E,ls,1\}} U^{\mu,\alpha}_{k}||_{X^1(\mathbb{R})}\\
&\leq ||\sum_{1 \leq \alpha \leq A} \sum_{\mu \in \{ E,ls,1\}} U^{\mu,\alpha}_{k}||_{X^1(\mathbb{R})}+||\sum_{A \leq \alpha \leq J} \sum_{\mu \in \{ E,ls,1\}} (U^{\mu,\alpha}_{k}-e^{it\Delta}U^{\mu,\alpha}_{k}(0))||_{X^1(\mathbb{R})}\\
&+||e^{it\Delta} \sum_{A \leq \alpha \leq J} \sum_{\mu \in \{ E,ls,1\}} U^{\mu,\alpha}_{k}(0)||_{X^1(\mathbb{R})}\\
&\leq ||\sum_{1 \leq \alpha \leq A} \sum_{\mu \in \{ E,ls,1\}} U^{\mu,\alpha}_{k}||_{X^1(\mathbb{R})}+\sum_{A \leq \alpha \leq J} \sum_{\mu \in \{ E,ls,1\}} L_{\mu}(\alpha)^{\frac{3}{2}}\\
&+|| \sum_{A \leq \alpha \leq J} \sum_{\mu \in \{ E,ls,1\}} U^{\mu,\alpha}_{k}(0)||_{H^1(\mathbb{R}^2\times \mathbb{T}^2)}\\
&\lesssim A+\sum_{A \leq \alpha \leq J} \sum_{\mu \in \{ E,ls,1\}} L_{\mu}(\alpha)^{\frac{3}{2}}+|| \sum_{A \leq \alpha \leq J} \sum_{\mu \in \{ E,ls,1\}} U^{\mu,\alpha}_{k}(0)||_{H^1(\mathbb{R}^2\times \mathbb{T}^2)} \lesssim 1.
\end{align*}
\noindent The bound on $\sum\limits_{1 \leq \alpha \leq J} \sum\limits_{\mu \in \{ E,ls,1\}}||U^{\mu,\alpha}_{k}||_{X^1}^2$ is similar.\vspace{3mm}

\noindent  Then we are now ready to construct the approximate solution. Let $F(z)=z|z|^2$ and also we have
\begin{equation}\label{equation}  F^{'}(G)u=2G\bar{G}u+G^2 \bar{u} .
\end{equation}
\noindent For each $B$ and $J$, we define $g_k^{B,J}$ to be the solution of the initial value problem:
\begin{equation}\label{equation} i\partial_t g+\Delta g-F^{'}(U^B_{prof,k})g=0 ,\quad g(0)=R^{J}_{k}.
\end{equation}
\noindent The solution $g_{k}^{B,J}$ is well defined on $\mathbb{R}$ for $k>k_0(B,J)$ and satisfies:
 \begin{equation}\label{equation} ||g_k^{B,J}||_{X^1(\mathbb{R})} \leq Q^{'} .
\end{equation}
\noindent For some $Q^{'}$ independent of $J$ and $B$. This follows by splitting $\mathbb{R}$ into $O(Q)$ intervals $I_j$ over which $||U_{prof,k}^{B}||_{Z(I_j)}$ is small and applying the local theory on each subinterval.\vspace{3mm}

\noindent Now we can define the approximate solution. We let ($A$ will be chosen shortly)\vspace{3mm}

$U_k^{app,J}=U_{prof,J}^{A}+g_{k}^{A,J}+U^{>A}_{prof,k} \quad$ where $\quad U_{prof,k}^{>A}=\sum\limits_{A<\alpha \leq J}  \sum\limits_{\mu} U_{k}^{\mu,\alpha}$ \vspace{5mm}

\noindent which has $u_k(0)$ as its initial data and satisfies, for any $1\leq A \leq J$, the bound: 
\[||U_k^{app,J}||_{X^1(\mathbb{R})} \leq 3(Q+Q^{'})
\]\noindent for all $k \geq k_0(J)$. According to Theorem 4.5 with $M=6(1+Q+Q^{'})$ gives us an $\epsilon_1=\epsilon_1(M)\leq \frac{1}{K(1+Q+Q^{'})}$ for some $K$ sufficiently large, such that if the error term defined below whose $N$-norm is bounded by $\epsilon_1$, then we can upgrade the uniform $X^1(\mathbb{R})$ bounds on into a uniform bound on $||u_k||_{Z(\mathbb{R})}$, thus deriving a contradiction. First we choose $A$ such that:
\begin{equation}\label{equation} 
||U^{>A}_{prof,k}||^2_{X^1(\mathbb{R})}+\sum_{A<\alpha \leq J}  \sum_{\mu} ||U_{k}^{\mu,\alpha}||^2_{X^1(\mathbb{R})} \leq \epsilon_{1}^{10}.
\end{equation}
\noindent for any $J \geq A$ and $k$ sufficiently large.\vspace{3mm}

\noindent After fixing $A$ we can bound the error term:
\begin{align*}
e^{J}_{k}&=(i\partial_{t}+\Delta)U_k^{app,J}-F(U_k^{app,J})  \tag{7.21}\\
&=-F(U^A_{prof,k}+g_k^{A,J}+U^{>A}_{prof,k})+\sum_{1\leq \alpha \leq J,\mu} F(U_k^{\mu,\alpha})+F^{'}(U_{prof,k}^{A})g_k^{A,J}  \tag{7.22} \\
&=- F(U^{A}_{prof,k}+g_k^{A,J}+U^{>A}_{prof,k})+F(U^A_{prof,k}+g^{A,J}_{k})+F(U^{>A}_{prof,k})  \tag{7.23}\\
&-F(U^A_{prof,k}+g_k^{A,J})+F(U^{A}_{prof,k})+F^{'}(U_{prof,k}^{A})g_k^{A,J}    \tag{7.24} \\
&-F(U^A_{prof,k})+\sum_{1 \leq \alpha \leq A} F(U_k^{\mu,\alpha})   \tag{7.25} \\
&-F(U^{>A}_{prof,k})+\sum_{A+1 \leq \alpha \leq J,\mu} F(U^{\mu,\alpha}_k)    \tag{7.26}.
\end{align*}
\noindent We will estimate the four terms separately.\vspace{3mm}

\noindent By Lemma 4.2, (7.16) and (7.20), we estimate:
\[||(7.23)||_{N(\mathbb{R})} \lesssim (||U^{A}_{prof,k}+g^{A,J}_k||_{X^1(\mathbb{R})}+||U^{>A}_{prof,k}||_{X^1(\mathbb{R})})^2 ||U^{>A}_{prof,k}||_{X^1(\mathbb{R})}< \frac{\epsilon_1}{4}
\]
\noindent for $k$ large enough. By Lemma 4.2 and Lemma 7.5, we estimate:
\[||(7.24)||_{N(\mathbb{R})}\lesssim  (||U^{A}_{prof,k}||_{X^1(\mathbb{R})}+||g^{A,J}_k||_{X^1(\mathbb{R})})^2 ||g^{A,J}_k||_{Z^{'}(\mathbb{R})} \lesssim (Q+Q^{'})^2||g^{A,J}_k||_{Z^{'}(\mathbb{R})}< \frac{\epsilon_1}{4}
\]
\noindent if $J$ is big enough and for $k>k_0(J)$. By Lemma 7.4, we estimate:
\[||(7.25)||_{N(\mathbb{R})}\lesssim_A \sum_{(\alpha_1,\mu_1)\neq (\alpha_2,\mu_2)}||\tilde{U}^{\mu_1,\alpha_1}_{k} \tilde{U}^{\mu_2,\alpha_2}_{k} \tilde{U}^{\mu_3,\alpha_3}_{k}||_{N(\mathbb{R})}<\frac{\epsilon_1}{4}
\]
\noindent if $k$ is big enough. By (7.20), we estimate:
\[||(7.26)||_{N(\mathbb{R})} \lesssim ||U^{>A}_{prof,k}||^3_{X^1(\mathbb{R})}+\sum_{A<\alpha \leq J}||U^{\mu,\alpha}_{prof,k}||^3_{X^1(\mathbb{R})} \leq \frac{\epsilon_1}{4}.
\]\noindent By using Theorem 4.5, we get that $u_k$ extends as a solution in $X_c^{1}(\mathbb{R})$ satisfying:
\[||u_k||_{Z(\mathbb{R})}<+\infty
\]
\noindent which contradicts (7.1). \vspace{3mm}

\noindent There are two more lemmas which are used in the estimates above. 
\begin{lemma} Assume that $U^{\alpha}_{k}, U^{\beta}_{k}, U^{\gamma}_{k}$ are three nonlinear profiles from the set $\{U^{\mu,\alpha}_k:1\leq \alpha \leq A,\mu \in \{E,ls,1 \} \}$ such that $U^{\alpha}_{k}$ and $U^{\beta}_{k}$ correspond to orthogonal frames. Then for these nonlinear profiles:
\[ \limsup\limits_{k \rightarrow +\infty} ||\tilde{U}^{\alpha}_{k} \tilde{U}^{\beta}_{k} \tilde{U}^{\gamma}_{k}||_{N(\mathbb{R})} =0
\]\noindent where for $\delta \in \{\alpha,\beta,\gamma\}$ , $\tilde{U}^{\delta}_{k} \in \{U^{\delta}_{k},\bar{U}^{\delta}_{k}\}$.
\end{lemma}
\begin{lemma}For any fixed $A$, it holds that: 
\[ \limsup\limits_{J \rightarrow \infty} \limsup\limits_{k \rightarrow \infty} ||g_k^{A,J}||_{Z(\mathbb{R})}=0.
\]
\end{lemma}
\noindent The proofs of Lemma 7.4 and Lemma 7.5 are similar to the proofs in [13, 19].
\section{Local Theory of the Resonant System}
\begin{theorem}[Local well-posedness and small-data scattering for (1.6)]. Let $\vec{u}(0)=\{u_p(0)\}_p \in h^1L^2$ satisfies $||\vec{u}_{0}||_{h^1L^2} \leq E$, then:\vspace{3mm}

\noindent (1)There exists an open interval $I$ which contains $0$ and a unique solution $\vec{u}(t)$ of (1.6) in $C^0_{t}(I:h^1L^2) \cap \vec{W}(I)$.\vspace{3mm}

\noindent (2)There exists $E_0$ such that $E \leq E_0$, $\vec{u}(t)$ is global and scatters in both directions.\vspace{3mm}

\noindent (3)Persistence of regularity: if $\vec{u}(0)\in h^{\eta}H^k$ for some $\eta \geq 1$ and $k \geq 0$, then $\vec{u} \in C^0_{t}(I:h^{\eta}H^k)$.
\end{theorem}
\noindent \emph{Proof:} The proof follows from a simple fixed point theorem(and classical arguments), once we have established the nonlinear estimate. By Strichartz estimate, we have
\[||u_j||_{L^{4}_{x,t}(\mathbb{R}^2\times I)} \lesssim ||u_j(0)||_{L^2}+\sum_{R(j)}||u_{p_1}\bar{u}_{p_2}u_{p_3}||_{L^{\frac{4}{3}}_{x,t}(\mathbb{R}^2\times I)}
\]
\noindent where $R(j)$ was defined in (1.6). Multiplying by $\langle j \rangle$ and square-summing, the first term on the right-handed side is bounded by the square of the $h^1L^2_x-$norm.\vspace{3mm}

\noindent For the second term, we compute as follows
\begin{align*}
&\sum_{j \in \mathbb{Z}^2} \langle j \rangle^2 [\sum_{R(j)} ||u_{p_1}\bar{u}_{p_2}u_{p_3}||_{L^{\frac{4}{3}}_{x,t}(\mathbb{R}^2 \times I)}]^2 \\
&\lesssim \sum_{j \in \mathbb{Z}^2} \langle j \rangle^2 [\sum_{R(j)}\Pi_{k=1}^{3} ||u_{p_k}||_{L^{4}_{x,t}(\mathbb{R}^2 \times I)}]^2\\
&\lesssim \sum_{j \in \mathbb{Z}^2} \{\sum_{R(j)}\Pi_{k=1}^{3} \langle p_k \rangle^2 ||u_{p_k}||_{L^{4}_{x,t}(\mathbb{R}^2 \times I)} \times \langle j \rangle^2 \sum_{R(j)}\langle p_1 \rangle^{-2}\langle p_2 \rangle^{-2}\langle p_3 \rangle^{-2} \} \\
&\lesssim \sum_{j \in \mathbb{Z}^2} \sum_{R(j)}\Pi_{k=1}^{3} \langle p_k \rangle^2 ||u_{p_k}||_{L^{4}_{x,t}(\mathbb{R}^2 \times I)} \lesssim ||\vec{u}||^6_{\vec{W}(I)}.
\end{align*}
\noindent We obtain
\[||\vec{u}||_{\vec{W}(I)} \lesssim ||e^{it\Delta_x} \vec{u}_0||_{\vec{W}(I)}+||\vec{u}||^3_{\vec{W}(I)}.
\]
\noindent Also we know
\[||e^{it\Delta_x} \vec{u}_0||_{\vec{W}(\mathbb{R})} \lesssim ||\vec{u}_0||_{h^1L^2} \lesssim E .
\]
\noindent We can now run a classical fixed-point argument in $W(I) \cap C_t(I:h^1L^2)$ provided $I$ or $E$ is small enough. The rest of the theorem follows from standard arguments.
\begin{lemma}There holds that 
\[\sup\limits_{j \in \mathbb{Z}^2}\{\langle j \rangle^2 \sum_{R(j)}\langle p_1 \rangle^{-2}\langle p_2 \rangle^{-2}\langle p_3 \rangle^{-2}\} \lesssim 1.
\]
\end{lemma}
\noindent \emph{Proof:} Without loss of generality, we may assume that
\[|p_1|\leq |p_3|, \quad max(|j|,|p_2|)\sim |p_3|.
\]
\noindent Also we can see that $p_1$ is on a specific circle $\mathcal{C}$,
\[|p_1-\frac{p_2-j}{2}|^2=(\frac{p_2-j}{2})^2.
\]
\begin{align*}
S_1 &=\sum_{(p_1,p_2,p_3)\in R(j);|p_1|\leq|p_3|;|p_2|\leq|p_1|} \langle p_1 \rangle^{-2} \langle p_2 \rangle^{-2} \frac{\langle j \rangle^{2}}{\langle p_3 \rangle^{2}}\\
&\lesssim \sum_{(p_1,p_2,j+p_2-p_1)\in R(j);|p_2|\leq|p_1|} \langle p_1 \rangle^{-2} \langle p_2 \rangle^{-2} [\frac{\langle j \rangle}{\langle max(|j|,|p_2|) \rangle}]^2\\
&\lesssim \sum_{p_2} \langle p_2 \rangle^{-2} \sum_{p_1} \langle p_1 \rangle^{-2}\\
&\lesssim \sum_{p_2} \langle p_2 \rangle^{-2} \langle |p_2| \rangle^{-1} \lesssim 1.
\end{align*}

\noindent The sum when $|p_1|\leq|p_2|$ is bounded similarly, using the following lemma to bound the sum over $p_2$ instead of the bound over $p_1$.

\begin{lemma}For any $P\in \mathbb{R}^2$, $R>0$ and $A>1$ there hold that:
\[\sum_{|p|\geq A,p\in \mathbb{Z}^2\cap C(P,R)}\frac{1}{\langle p \rangle^2} \lesssim A^{-1}
\]
where $C(P,R)$ denotes the circle of radius $R$ centered at $P$.
\end{lemma}
\noindent \emph{Proof:} It is exactly as same as [13, Lemma 8.3].
\begin{lemma}Assume the conclusion of Theorem 1.2 holds for all initial data $u_0\in H^1(\mathbb{R}^2\times \mathbb{T}^2)$ with full energy $L(u_0)<E_{max}$. Then Conjecture 1.2 holds true for all initial data $\vec{v}_0 \in h^1L^2_x$ satisfying $E_{ls}(\vec{v}_0)<E_{max}$. Furthermore, if all finite full-energy solutions scatter for (1.1), then the same thing holds for finite $E_{ls}$-energy solutions of (1.6).
\end{lemma}
\noindent \emph{Proof:} It follows as in [13, Lemma 8.4].\vspace{3mm}

\noindent At last, we discuss the proof of Lemma 3.3 (Local-in-time $L^p$ estimate). The idea of the proof is similar to [19, Proposition 2.1]. The main ingredient is the following distributional inequality:
\begin{lemma}Assume $p_0>\frac{10}{3}$, $N\geq 1$, $\lambda\in [N^{\frac{2p_0-6}{p_0-2}},2^{10}N^2]$, $||m||_{L^2(\mathbb{R}^2\times \mathbb{Z}^2)}\leq 1$, and $m(\xi)=0$ for $|\xi|>N$, then 
\begin{equation}
|\{(x,t)\in \mathbb{R}^2\times \mathbb{T}^2 \times [-2^{-10},2^{-10}]:|\int_{\mathbb{R}^2\times \mathbb{Z}^2} m(\xi)e^{-it|\xi|^2}e^{ix\cdot \xi}d\xi|\geq \lambda\}|\lesssim N^{2p_0-6} \lambda^{-p_0}.
\end{equation}
\end{lemma}
\noindent First of all, Lemma 3.3. follows from the lemma above (Lemma 8.5). \vspace{3mm}

\noindent \emph{Proof of Lemma 3.3:} We let
\[F(x,t)=\int_{\mathbb{R}^2\times \mathbb{Z}^2}m(\xi)e^{-it|\xi|^2}e^{ix\cdot \xi}d\xi,
\]
\noindent where $m$ is as in Lemma 8.5, it suffices to prove that if $p>\frac{10}{3}$ and $N\geq 1$, then
\begin{equation}
||1_{[-2^{-10},2^{-10}]}(t)F||_{L^p(\mathbb{R}^2\times \mathbb{T}^2 \times \mathbb{R})}\lesssim_{p} N^{2-\frac{6}{p}}.
\end{equation}
\noindent We may assume $p\in (\frac{10}{3},4]$ and $N\gg 1$. Then
\[||1_{[-2^{-10},2^{-10}]}(t)F||^p_{L^p(\mathbb{R}^2\times \mathbb{T}^2\times \mathbb{R})}\leq \sum_{2^l\leq 2^{10}N^2}2^{pl}|\{(x,t)\in \mathbb{R}^2\times \mathbb{T}^2\times [-2^{-10},2^{-10}]:|F(x,t)|\geq 2^l\}|.
\]
\noindent If $2^l\geq N^{\frac{2p_0-6}{p_0-2}}$, $p_0\in (\frac{10}{3},p)$, we use the distributional inequality (8.1). If $2^l \leq N^{\frac{2p_0-6}{p_0-2}}$, we use the following bound:
\[2^{2l}|\{(x,t)\in \mathbb{R}^2\times \mathbb{T}^2\times [-2^{-10},2^{-10}]:|F(x,t)|\geq 2^l\}|\leq ||F||^2_{L^2(\mathbb{R}^2\times \mathbb{T}^2\times \mathbb{R})} \lesssim 1
\]
\noindent Therefore
\begin{equation}
\aligned
||1_{[-2^{-10},2^{-10}]}(t)F||^p_{L^p(\mathbb{R}^2\times \mathbb{T}^2\times \mathbb{R})} &\lesssim \sum_{2^l\leq N^{\frac{2p_0-6}{p_0-2}}}2^{(p-2)l}+\sum_{N^{\frac{2p_0-6}{p_0-2}} \leq 2^l\leq 2^{10}N^2}2^{pl}\cdot N^{2p_0-6}2^{-p_0 l}\\
&\lesssim N^{2p-6},
\endaligned
\end{equation}
\noindent which gives (8.2). It suffices to prove Lemma 8.5. Now we focus on Lemma 8.5.\vspace{3mm}

\noindent \emph{Proof of Lemma 8.5: }We may assume that $N\gg 1$ and consider the kernel $K_N:\mathbb{R}^2\times \mathbb{T}^2 \times \mathbb{R}\rightarrow \mathbb{C}$,
\begin{equation}
K_N(x,t)=\eta^1(2^5t/(2\pi))\int_{\mathbb{R}^2\times \mathbb{Z}^2} e^{-it|\xi|^2} e^{ix\cdot \xi} \eta^4(\xi/N) d\xi
\end{equation}
\noindent Let
\[S_{\lambda}:=\{(x,t)\in \mathbb{R}^2\times \mathbb{T}^2 \times [-2^{-10},2^{-10}]:|\int_{\mathbb{R}^2\times \mathbb{T}^2} m(\xi)e^{-it|\xi|^2}e^{ix\cdot \xi} d\xi|\geq \lambda\}
\]
\noindent and fix a function $f:\mathbb{R}^2\times \mathbb{T}^2 \times [-2^{-10},2^{-10}]\rightarrow \mathbb{C}$ such that
\begin{equation}
|f|\leq 1_{S_{\lambda}}
\end{equation}
\noindent and 
\begin{equation}
\lambda|S_{\lambda}|\leq |\int_{\mathbb{R}^2\times \mathbb{T}^2\times \mathbb{R}} f(x,t)\cdot [\int_{\mathbb{R}^2\times \mathbb{T}^2} m(\xi)e^{-it|\xi|^2}e^{ix\cdot \xi}d\xi]dxdt|.
\end{equation}
\noindent Using the assumption on $m$ we estimate the right-hand side of the inequality (8.6) above by 
\[||\prod\limits_{j=1}^{4} \eta^1(\xi_j/N) \cdot \int_{\mathbb{R}^2\times \mathbb{T}^2\times \mathbb{R}}f(x,t)e^{-it|\xi|^2}e^{ix\cdot \xi}dxdt||_{L^2_{\xi}}.
\]
\noindent Thus
\begin{equation}
\lambda^2|S_{\lambda}|^2 \leq \int_{\mathbb{R}^2\times \mathbb{T}^2 \times \mathbb{R}} \int_{\mathbb{R}^2\times \mathbb{T}^2 \times \mathbb{R}} f(x,t)\overline{f(y,s)}K_N(t-s,x-y)dtdxdsdy.
\end{equation}
\noindent Using Lemma 8.6 below, we estimate the right-hand-side in (8.7) as follows
\begin{equation}
\aligned
&\int_{\mathbb{R}^2\times \mathbb{T}^2 \times \mathbb{R}} \int_{\mathbb{R}^2\times \mathbb{T}^2 \times \mathbb{R}} f(x,t)\overline{f(y,s)}K_N(t-s,x-y)dtdxdsdy\\
&\leq \sum_{\mu \in \{1,2,3\}}|\int_{\mathbb{R}^2\times \mathbb{T}^2 \times \mathbb{R}} \int_{\mathbb{R}^2\times \mathbb{T}^2 \times \mathbb{R}} f(x,t)\overline{f(y,s)}K^{\mu,\lambda}_N(t-s,x-y)dtdxdsdy| \\
&\leq (\lambda^2/2)||f||^2_{L^1}+C\lambda^2(N^{2p_0-6}\lambda^{-p_0})||f||^2_{L^2}+C\lambda^2(N^{2p_0-6}\lambda^{-p_0})^{\frac{r-1}{r}}||f||^2_{L^{\frac{2r}{r+1}}}\\
&\leq (\lambda^2/2)|S_{\lambda}|^2+C\lambda^2(N^{2p_0-6}\lambda^{-p_0})|S_{\lambda}|+C\lambda^2(N^{2p_0-6}\lambda^{-p_0})^{\frac{r-1}{r}}|S_{\lambda}|^{\frac{r+1}{r}}.
\endaligned
\end{equation}
\noindent Using (8.7), it follows that
\[|S_{\lambda}|\lesssim N^{2p_0-6}\lambda^{-p_0}+(N^{2p_0-6}\lambda^{-p_0})^{\frac{r-1}{r}}|S_{\lambda}|^{\frac{1}{r}}.
\]
\noindent which gives (8.1) and finishes the proof of Lemma 8.5. Now, we only need to prove the following lemma, which is a crucial step in the proof of Lemma 8.5.
\begin{lemma}Assume $\lambda\in [N^{\frac{2p_0-6}{p_0-2}},2^{10}N^2]$ as in Lemma 8.5 and $r\in [2,4]$, there is a decomposition
\[K_N=K^{1,\lambda}_N+K^{2,\lambda}_N+K^{3,\lambda}_N
\] 
\noindent such that
\[||K_N^{1,\lambda}||_{L^{\infty}(\mathbb{R}^2\times \mathbb{T}^2\times \mathbb{R})} \leq  \frac{\lambda^2}{2},
\]
\begin{equation}
||\widehat{K_N^{2,\lambda}}||_{L^{\infty}(\mathbb{R}^2\times \mathbb{T}^2\times \mathbb{R})}\lesssim \lambda^2(N^{2p_0-6}\lambda^{-p_0}),
\end{equation}
\[||\widehat{K_N^{3,\lambda}}||_{L^{r}(\mathbb{R}^2\times \mathbb{T}^2\times \mathbb{R})}\lesssim \lambda^2(N^{2p_0-6}\lambda^{-p_0})^{\frac{r-1}{r}}.
\]
\end{lemma}

\noindent \emph{Proof of Lemma 8.6:} For a continuous function $h:\mathbb{R}\rightarrow \mathbb{C}$ and any $(\xi,\tau) \in \mathbb{R}^2\times \mathbb{Z}^2\times \mathbb{R}$
\begin{equation}
\mathcal{F}[K_N(x,t)\cdot h(t)](\xi,\tau)=C\eta^4(\xi/N)\int_{\mathbb{R}}h(t)\eta^1(2^5t/(2\pi))e^{-it(\tau+|\xi|^2)}dt.
\end{equation}
\noindent It is shown in [2, Lemma 3.18] that
\begin{equation}
|\sum_{n\in\mathbb{Z}}e^{-it|n|^2}e^{ixn}\eta^1(\xi/N)^2|\lesssim \frac{N}{\sqrt[]{q}(1+N|t/(2\pi)-a/q|^{1/2})}
\end{equation}
\noindent if
\begin{equation}
t/(2\pi)=a/q+\beta,\quad q\in \{1,...,N\},a\in \mathbb{Z},(a,q)=1,|\beta|\leq (Nq)^{-1}.
\end{equation}
\noindent Also for $t$ as in (8.12), we have,
\begin{equation}
|K_N(x,t)|\lesssim \frac{N^2}{q(1+N|t/(2\pi)-a/q|^{1/2})^2} \cdot (\frac{N}{1+N|t/(2\pi)|^{\frac{1}{2}}})^2.
\end{equation}
\noindent For $j\in \mathbb{Z}$ we define $\eta_j$, $\eta \geq 0:\mathbb{R}\rightarrow [0,1],$
\[\eta_j(s):=\eta^1(2^j s)-\eta^1(2^{j+1}s), \quad \eta_{\geq j}(s):=\sum_{k\geq j}\eta_k(s).
\]
\noindent Fix integers $K$, $L$, satisfying
\begin{equation}
K\in \mathbb{Z}_{+},\quad N\in [2^{K+4},2^{K+5}),\quad L\in \mathbb{Z}\cap[0,2K+20],\quad \lambda^{p_0-2}N^{6-2p_0}\in [2^L,2^{L+1}). 
\end{equation}

\noindent Now we start with the important decomposition as follows:
\[1=[\sum_{k=0}^{K-1} \sum_{j=0}^{K-k} p_{k,j}(s)]+e(s),
\]
\begin{equation}
p_{k,j}(s):=\sum_{q=2^k}^{2^{k+1}-1} \sum_{a\in \mathbb{Z},(a,q)=1} \eta_{j+K+k+10}(s/(2\pi)-a/q) \quad \textmd{if }j\leq K-k-1,
\end{equation}
\[p_{k,K-k}(s):=\sum_{q=2^k}^{2^{k+1}-1} \sum_{a\in \mathbb{Z},(a,q)=1} \eta_{\geq 2K+10}(s/(2\pi)-a/q).
\]
\noindent Let $T_K=\{(k,j)\in \{0,...,K-1\}\times \{0,...,K\}:k+j\leq K\}$. In view of Dirichlet's lemma, we observe that
\begin{equation}
\textmd{if }t\in supp(e) \textmd{ satisfies }(8.12), \textmd{ then either }N\lesssim q \textmd{ or } (Nq)^{-1}\approx |t/(2\pi)-a/q|. 
\end{equation}
\noindent We define the first component of the kernel $K_N^{2,\lambda}$,
\begin{equation}
K^{2,\lambda}_{N,1}(x,t)=K_N(x,t)\cdot \eta^1(2^{L-40}t/(2\pi)).
\end{equation}
\noindent It follows from (8.10) that
\begin{equation}
||\widehat{K^{2,\lambda}_{N,1}}||_{L^{\infty}(\mathbb{R}^2\times \mathbb{T}^2 \times \mathbb{R})} \lesssim 2^{-L} \lesssim N^{2p_0-6} \lambda^{2-p_0},
\end{equation}
\noindent which agrees with (8.9).\vspace{3mm}

\noindent Therefore we may assume that $L\geq 45$ and write 
\begin{equation}
K_N(x,t)-K^{2,\lambda}_{N,1}(x,t)=\sum_{l=4}^{L-41} K_N(x,t)\cdot \eta_{l}(t/(2\pi)).
\end{equation}
\noindent Using (8.13) and (8.16), for any $(k,j)\in T_K$ and $l \in [4,L-41] \cap \mathbb{Z}$,
\begin{equation}
\sup\limits_{x,t}|K_N(x,t)\cdot \eta_{l}(t/(2\pi))p_{k,j}(t)| \lesssim 2^{l}2^{K+j},
\end{equation}
\begin{equation}
\sup\limits_{x,t}|K_N(x,t)\cdot \eta_{l}(t/(2\pi))e(t)| \lesssim 2^{l}2^{K}.
\end{equation}
\noindent We analyze two cases.\vspace{3mm}

\noindent \emph{Case 1: }$L \leq 2K-\delta K,\delta=\frac{1}{100}$. In this case we set
\begin{equation}
K^{1,\lambda}_N(x,t):=K_N(x,t)\cdot \big[\sum_{l=4}^{L-41}\eta_l(t/(2\pi))[e(t)+\sum_{k,j\in T_K,2j\leq L}p_{k,j}(t)+\sum_{k,j\in T_K,2j> L}\rho_{k,j} p_{k,0}(t)]\big],
\end{equation}
\begin{equation}
K^{2,\lambda}_N(x,t):=K^{2,\lambda}_{N,1}(x,t)+K_N(x,t)\cdot  \big[\sum_{l=4}^{L-41}\eta_l(t/(2\pi))[\sum_{k,j\in T_K,2j> L}(p_{k,j}(t)-\rho_{k,j} p_{k,0}(t))]\big],
\end{equation}
\begin{equation}
K^{3,\lambda}_N(x,t):=0,
\end{equation}
\noindent where $K_{N,1}^{2,\lambda}$ is define as in (8.17) and 
\begin{equation}
\rho_{k,j}:=2^{-j} \textmd{ if } j\leq K-k-1 \textmd{ and } \rho_{k,N-k}:=2^{-K+k+1}.
\end{equation}
\noindent The bound on $K_N^{1,\lambda}$ is trivial. Using (8.20) and (8.21), for any $(x,t)\in \mathbb{R}^2\times \mathbb{T}^2 \times \mathbb{R}$, we have
\begin{equation}
\aligned
|K_N^{1,\lambda}(x,t)| &\leq \sum_{l=4}^{L-41}|K_N(x,t) \cdot \eta_l(t/(2\pi)) e(t)| \\
&+\sum_{l=4}^{L-41} \sum_{k,j\in T_K,2j\leq L}|K_N(x,t) \cdot \eta_l(t/(2\pi))p_{k,j}(t)|\\
&+\sum_{l=4}^{L-41} \sum_{k,j\in T_K,2j> L}|K_N(x,t) \cdot \eta_l(t/(2\pi))\rho_{k,j}p_{k,0}(t)| \\
&\lesssim 2^{L}2^{K}+K2^{L}2^{K+\frac{L}{2}}+K2^{L}2^{K}\\
&\lesssim K2^{(L-2K)(3p_0-10)/2(p_0-2)}\cdot 2^{2L/(p_0-2)} 2^{K(4p_0-12)/(p_0-2)}.
\endaligned
\end{equation}
\noindent Recall that $(2K-L)\geq \delta K,N\gg 1$, and $p_0>\frac{10}{3}$. Notice that $\lambda^2 \approx 2^{2L/(p_0-2)}2^{K(4p_0-12)/(p_0-2)}$ the desired bound $|K_N^{1,\lambda}(x,t)| \leq \lambda^2/2$ follows.\vspace{3mm}

\noindent It remains to bound (8.9) on the kernel $K_N^{2,\lambda}$ which follows as in [19].\vspace{3mm}

\noindent \emph{Case 2: }$L \geq 2K-\delta K,\delta=\frac{1}{100}$. For $b\in \mathbb{Z}_{+}$ sufficiently large, in this case we set
\begin{equation}
K_N^{1,\lambda}(x,t):=K_N(x,t)\cdot [\sum_{l=4}^{L-41}\eta_l(t/(2\pi))[e(t)+\sum_{k,j\in T_K,2j\leq L-b}p_{k,j}(t)]].
\end{equation}
\noindent Using the bound (8.20) and (8.21), for $(x,t)\in \mathbb{R}^2\times \mathbb{T}^2 \times \mathbb{R}$
\begin{equation}
\aligned
|K_N^{1,\gamma}(x,t)| &\lesssim \sum_{l=4}^{L-41}|K_N(x,t)\cdot \eta_l(t/(2\pi)) e(t)|+\sum_{l=4}^{L-41}\sup\limits_{k,j\in T_K,2j\leq L-b}|K_N(x,t)\cdot \eta_l(t/(2\pi)) p_{k,j}(t)|\\
&\lesssim 2^{L}2^{K}+2^{L}2^{K+L/2}2^{-b/2}\\
&\lesssim 2^{-b/2}2^{(L-2K)(3p_0-10)/2(p_0-2)}\cdot 2^{2L/(p_0-2)} 2^{K(4p_0-12)/(p_0-2)}.
\endaligned
\end{equation}
\noindent Since $\lambda^2 \approx 2^{2L/(p_0-2)}2^{K(4p_0-12)/(p_0-2)}$, it follows that $|K_N^{1,\gamma}(x,t)|\leq \lambda^2/2$ if $b$ is fixed sufficiently large.\vspace{3mm}

\noindent Now, let
\[L_N(x,t):=K_N(x,t)-K^{2,\lambda}_{N,1}(x,t)-K^{1,\lambda}_{N}(x,t)=\sum_{l=4}^{L-41} \sum_{k,j\in T_K,2j>L-b}K_N(x,t)\cdot \eta_l(t)p_{k,j}(t).
\]
\noindent It remains to prove that one can decompose $L_N=L_N^{2,\lambda}+L_N^{3,\lambda}$ satisfying
\begin{equation}
||\widehat{L_N^{2,\lambda}}||\lesssim 2^{-L},\quad ||\widehat{L_N^{3,\lambda}}||\lesssim \lambda^{2/r}2^{-L(r-1)/r}.
\end{equation}
\noindent We let $\tilde{\eta}_{L}(s):=\sum_{l=4}^{L-41}\eta_l(s)=\eta^1(2^4s)-\eta^1(2^{L-40}s)$. And we have
\begin{equation}
\widehat{L_N}(\xi,\tau)=C\sum_{k,j\in T_K,2j>L-b}\eta^4(\xi/N)\int_{\mathbb{R}}\tilde{\eta}_L(t)\eta^1(2^5 t)p_{k,j}(2\pi t)e^{-2\pi it(\tau+|\xi|^2)} dt.
\end{equation}
\noindent The cardinality of the set $T_{K,L}:=\{k,j\in T_K:2j>L-b\}$ is bounded by $C(1+|2K-L|)^2$. Let $f_{k,j}:\mathbb{R}\rightarrow \mathbb{C}$,
\begin{equation}
f_{k,j}(\mu):=\int_{\mathbb{R}} \tilde{\eta}_L(t)\eta^1(2^5 t)p_{k,j}(2\pi t)e^{-2\pi it\mu} dt.
\end{equation}
\noindent It suffices to prove that for $(k,j)\in T_{K,L}$, one can decompose
\begin{equation}
f_{k,j}=f^2_{k,j}+f^3_{k,j},
\end{equation}
\noindent satisfying
\begin{equation}
||f^2_{k,j}||_{L^{\infty}(\mathbb{R})}\lesssim 2^{-L}(1+|2K-L|)^{-2}
\end{equation}
\noindent and 
\begin{equation}
||f^3_{k,j}||_{L^{r}(\mathbb{R})}\lesssim (\lambda/N^2)^{2/r} 2^{-L(r-1)/r}(1+|2K-L|)^{-2}.
\end{equation}
\noindent The rest of the proof follows as in [19, Proposition 2.1]. \vspace{3mm}

\noindent \textbf{Acknowledgments.} I would like to express my deepest thanks to my advisor Professor Benjamin Dodson for many useful suggestions and comments. He has always been nice and has helped me a lot during my three years' study at Johns Hopkins University. And I really appreciate Professor Zaher Hani and Professor Benoit Pausader for their kind help and valuable comments. I would also like to thank C. Fan, C. Luo, T. Ren, Q. Su, M. Wang, L. Zhao and some other friends for insightful discussions and encouragement. Moreover, part of this work was progressed when the author was at AIM workshop (Wave Turbulence) and I would like to thank the organizers for hosting.
\\

\hfill \linebreak
\noindent \author{Zehua Zhao}

\noindent \address{Johns Hopkins University, Department of Mathematics, 3400 N. Charles Street, Baltimore, MD 21218, U.S.}

\noindent \email{zzhao25@jhu.edu}\\

\end{document}